%% file: main.tex
\documentclass[3p, preprint, number, sort&compress]{elsarticle}

\makeatletter
\def\ps@pprintTitle{%
 \let\@oddhead\@empty
 \let\@evenhead\@empty
 \def\@oddfoot{}%
 \let\@evenfoot\@oddfoot}
\makeatother

\usepackage{amsmath,amssymb}
\usepackage[T1]{fontenc}
\usepackage{color}
\usepackage{lmodern}
\usepackage{lineno,hyperref}
\usepackage[toc,page]{appendix}
\usepackage{stmaryrd}
\modulolinenumbers[5]

\hypersetup{colorlinks=false}
\hypersetup{citecolor=green}
\hypersetup{linkcolor=red}
\hypersetup{urlcolor=magenta}
\usepackage{graphicx}
\usepackage{subfigure}
\graphicspath{{figures/}}
\usepackage{array}
\newcolumntype{P}[1]{>{\centering\arraybackslash}p{#1}}
\newcolumntype{M}[1]{>{\centering\arraybackslash}m{#1}}
\usepackage{multirow}
\usepackage[dvipsnames]{xcolor}
\usepackage{defjs3}
\usepackage[section]{placeins} 
\usepackage{hhline}
\usepackage{comment}
\usepackage{cleveref}
\usepackage{xcolor}
\usepackage{float}

\usepackage[]{natbib}
\setcitestyle{numbers}

\usepackage{lineno}

\begin{document}

\large

\begin{frontmatter}
\nolinenumbers
\title{Multi-phase-field elasticity model based on partial rank-one energy relaxation on pairwise interfaces}

\author[Dor,IFM]{Mohammad Sarhil\corref{cor1}}
 \cortext[cor1]{Corresponding author}
 \ead{mohammad.sarhil@tu-dortmund.de}
\author[ICAMS]{ Oleg Shchyglo}
\author[ICAMS]{ Hesham Salama} 
\author[IFM]{ Dominik Brands}
\author[ICAMS]{ \\ Ingo Steinbach}
\author[IFM]{ J\"org Schr\"oder}

\address[Dor]{Institut für Baumechanik, Statik und Dynamik, Faculty of Architecture and Civil Engineering, TU Dortmund, August-Schmidt-Str. 8, 44227 Dortmund, Germany}
\address[IFM]{Institut f\"ur Mechanik, Abteilung Bauwissenschaften, Fakult\"at f\"ur Ingenieurwissenschaften, Universit\"at Duisburg-Essen, \\ Universit\"atsstr. 15, 45141 Essen, Germany}
\address[ICAMS]{Interdisciplinary Centre for Advanced Materials Simulation (ICAMS), Ruhr-Universit\"at Bochum, \\ Universit\"atsstr. 150, 44801 Bochum, Germany}

\begin{abstract}
To model mechanically-driven phase transformations using the phase-field theory, suitable models are needed for describing the mechanical fields related to individual  phase-fields in the interfacial regions. They play a crucial role in obtaining the mechanical driving forces of phase-field evolution. Quantitative modeling requires satisfying the interfacial static equilibrium and kinematic compatibility conditions. To the best of our knowledge, no existing multi-phase-field elasticity model has been able to satisfy the jump conditions between all the locally-active phase-fields associated to their pairwise normals, except in the dual-phase-field regions. In this work, we introduce a novel multi-phase-field elasticity model based on the partial rank-one relaxation of the elastic energy density defined on the pairwise interfaces as a function of pairwise strains. These ad hoc pairwise definitions enable us to satisfy the static equilibrium and kinematic compatibility conditions between all the locally-active phase-fields.  Different numerical examples are presented, which compare the developed model against the equal-strain and equal-stress limiting cases. 
\end{abstract}

\begin{keyword}
constitutive modeling, phase-field modeling, multi-phase-field, rank-one relaxation
\end{keyword}

\end{frontmatter}

\nolinenumbers
\input{Introduction.tex}

\input{Section2.tex}
\input{Section3.tex}

\input{Section4.tex}

\input{Section5.tex}

\input{Conclusion.tex}



\bibliography{References.bib}
\bibliographystyle{plainnat}


\end{document}

%% file: Introduction.tex
\section{Introduction}
\label{sec:intro}
Modeling of the evolution of material microstructure is essential in many fields such as biology, hydrodynamics, chemistry, engineering, and materials science. It provides a deeper understanding of natural phenomena and enables the enhancement of manufacturing processes toward better materials. The evolution of the microstructure is generally determined by minimizing the total energy under consideration, which may include chemical energy, interfacial energy, mechanical energy, magnetic energy, electrostatic energy, and others. Conventionally, interfaces are treated mathematically as sharp boundaries separating the phases. Subsequently, the interface velocity is determined as an additional boundary condition, which requires explicit tracking of the interfaces. While this is possible for one-dimensional problems, it quickly becomes impractical and complex for higher-dimensional cases. The need to address moving interface problems in real three-dimensional systems has led to the development of phase-field theory. A general review of phase-field theory is available in \cite{Che:2002:pfm,Ste:2009:mas}. 
		
The phase-field method describes the microstructure's geometry using indicator field variables and requires no explicit tracking of the interfaces. The field variables used in the phase-field method are uniform within the associated phase and vary continuously across the interfaces, gradually approaching zero in the adjacent phase. Here, the term "phase" does not necessarily refer to a thermodynamic phase (e.g. austenite and martensite) and can indicate different orientations of one thermodynamic phase (e.g. martensite variants). Thus, phase-field models have diffuse interfaces with a certain thickness selected for numerical convenience. 
Historically, the phase-field method was used to model the microstructure evolution in solidifying systems with a focus on the prediction of dendritic patterns during solidification without explicitly tracking the liquid-solid interfaces. This was first established in \citep{Kob:1993:man,Kob:1994:ana}. Later, the phase-field method emerged as an advantageous and powerful computational approach to model and predict the microstructure evolution in various types of materials. In general, the applications of the phase-field method include solidification, grain growth, solid-state phase transformation, the modeling of thin films and surfaces, dislocation dynamics, crack propagation and ductile fracture, and electromigration \citep{SteSonHar:2010:pfm, SteApe:2006:mfp, Voi:2011:3pf, KosCuiOrt:2002:apf, MieTeiAld:2016:pfm, BarGraNue:2007:apf}.

Mainly, there are two types of phase-field models that have been developed independently by two communities \citep{MoeBlaWol:2008:ait,Che:2002:pfm}. The first type uses phase-field variables, which are related to microscopic parameters such as the local composition and long-range order parameter (physical order parameter models), e.g. \citep{Kha:1983:tos}. The second type of phase-field models uses an indicator field variable as a numerical tool to avoid tracking of the interfaces (mesoscopic models), e.g. \citep{StePezNesSeePriSchRez:1996:apf}. The field variables can also be classified as conserved or non-conserved variables. Conserved variables have to satisfy a local conservation condition or law, e.g., concentrations and temperature. Non-conserved variables do not have to fulfill any condition or conservation law, e.g.,  long-range order parameter. In this work, the indicator variable $\phi \in [0,1]$ is non-conserved and is called phase-field. The evolution of the conserved field variables is obtained by solving the Cahn-Hilliard equation \citep{Cah:1961:osd} while for the evolution of the non-conserved field variables,  the Allen-Cahn equation \citep{CahAll:1977:amt} should be used instead.

 For modeling the microstructure evolution in solid systems such as martensitic transformation and  Widmanst\"atten pattern formation, the mechanical response caused by differences in lattice parameter, orientation, or crystal symmetry between phases, variants, and grains dominates during the evolution. In the context of the phase-field method, the strain and stress fields of the individual phase-fields in the diffuse interfaces, which are not uniquely defined, are necessary to obtain the corresponding mechanical driving forces. Therefore, mechanical models are needed focusing only on elasticity in this work. 
These elastic models are classified into two categories \citep{AmmAppCaiSam:2009:cpf}: homogenization models that satisfy micro-macro energy conservation, known as the Hill-Mandel condition \citep{Hil:1963:epo} and interpolation models that violate it. Equal-strain or equal-stress assumptions \citep{Reu:1929:bdf,Voi:1887:tsu} combined with the energy conservation condition lead to the classical unrealistic Voigt/Taylor and Reuss/Sachs homogenization models \citep{AmmAppCaiSam:2009:cpf,MosShcMon:2014:anh,HilMie:2011:cot,HilMie:2012:apf,SchSchSelBoeNes:2015:sse} corresponding to the limiting upper and lower bounds of elasticity. Interpolation models are the most common in the phase-field community in which an energy formulation based on bulk quantities is derived without defining the mechanical responses of the active phase-fields (localization), thereby seeking simplification. A bulk energy formulation merged with a linearly interpolated transformation deformation gradient and effective compliance tensors is available in \citep{SteApe:2006:mfp,BorDuStrBoeShcHarSte:2016:mdo,SteSonHar:2010:pfm,BorEngBoeSchSte:2014:lse,BorEngMoslSchStei:2015:ldf,ShcDuEngSte:2019:pfs,Ste:2009:mas,Ste:2013:pfm} or with linearly interpolated effective stiffness tensor in \citep{Kha:1983:tos,IdeLevPreCho:2005:fes,SPaMueBreNes:2007:pfm,EsfChaLevCol:2018:mpf,LevIdesPre:2004:mso,ArtJinKha:2001:tdp,SeoHuLiChrOh:2002:cso,BasLev:2017:isw,ArtSluRoy:2005:pfm,BabLev:2020:fss,esfGhaLev:sim:sim}. 
 Higher-order interpolation functions are needed at the nanoscale, as imposed by the laws of thermodynamic equilibrium, where the first derivative of the interpolation function must vanish in the bulk \cite{LevPre:2002:tdl1,LevPre:2002:tdl2}. A nonlinear interpolated transformation deformation gradient was assumed as a linear combination of Bain tensors, which describe lattice distortions associated to phase transformations, multiplied with a second-order interpolation polynomial in \cite{ArtWanKha:2000:tdp,MoeBlaWol:2008:ait} or with a higher-order interpolation polynomial in \cite{LevRoyPre:2013:mts,BasLev:2019:fep,BasLev:2018:nmp,IdeChoLev:2008:fem,Lev:2014:pfa,Lev:2013:tcp,BabBasLev:2017:aaa,Lev:2013:pft,Lev:2018:pfa,BabLev:2018:pfa,BasLev:2020:mpi,LevRoy:015:mpf,BabLev:2020:smd}.  For preserving the specific volume of the variant along the path of the martensitic variant-variant transformation, an isochoric transformation deformation gradient as an exponential-logarithmic combination was presented in \cite{BasLev:2018:nmp,BasLev:2019:fep,BasLev:2017:isw,BasLev:2020:aef,TumaStu:2016:sei,RezStu:2020:pfm,TumaStup:2016:pfs,Lev:2018:pfa,BabLev:2018:pfa}.  However, all these interpolation models were unable to achieve stress-free variant-variant twin boundaries and the isochoric transformation deformation gradient even significantly increases these interfacial stresses compared to the nonlinear interpolation polynomials. 
 
Interpolation models achieved stress-free twin boundaries \cite{HilMie:2011:cot,HilMie:2012:apf,BasLev:2018:nmp,BasLev:2017:isw,ClaKna:2011:apm} using a parameterized transformation deformation gradient along the rank-one connection between two chosen variants. This requires knowing the solution of the twinning equation. However, one of two solutions actually has to be picked which imposes the preferred orientation of the interface. Moreover, the solution is not available for all kinds of transformation and this model cannot be expanded to the general multi-phase-field case.  A detailed analysis by comparing the diffuse and sharp interface solutions was conducted in \cite{DurWolMoe:2013:eoi,DurWolMoe:2015:aqp} to eliminate the interfacial excess energy and stresses for the diffuse interface case which ended with introducing a quantitative homogenization model applying equal strains parallel to the interface and equal stresses in the interface's normal direction. These are actually the kinematic compatibility and static equilibrium conditions of the interface (jump conditions) used to develop the dual-phase model in \cite{SchnTschChouSelzerBoehNest:2015:pfe,TscSchNes:2019:aec} by enforcing them explicitly
which implies transforming  the stresses and strains to a new orthogonal coordinate system where the normal of the interface is the new first coordinate axis and then splitting the mechanical fields additively into a part in the direction of the interface's normal and a part perpendicular to the normal (in the interface's plane). This model was generalized to the multi-phase-field case in \cite{SchSchSchReiHerSelBoeNes:2017:ots,HerSchSchSchReiSelNes:2018:mfm,SchReiHerrSchNes:2020:smd} by defining a favorable phase-field (the one with the highest value) where the jump conditions for the pairs involving this phase-field only were satisfied and the rest were ignored. Another generalization scheme was introduced in \cite{SchSchTscReiHerSchSelBri:2018:ssm,SchHerStrSelSchNes:2019:otm,SchSchStrMitSelNes:2018:mpf,SchHerSchHoeNes:2019:mpf,AmoSchSchNes:2018:cep,KubSchStrSchNes:2019:pfa,SchAmoSchNes:2020:iosf} by satisfying all the jump conditions but on one homogenized normal vector instead of multiple pairwise normal vectors which means all phase-fields share equal stress in the direction of this homogenized normal and equal strains in the perpendicular directions. Unfortunately, this homogenized normal follows the phase-field pair with the highest values and lacks accuracy for the remaining phase-field pairs. Recently, a partial rank-one homogenization scheme in the multi-phase-field framework has been introduced in \cite{ChaSchMoe:ace:2023} by satisfying the jump condition on $n-1$ pairs from the available $n(n-1)$/2 pairs in the junctions. 

In this work, we introduce a novel elasticity homogenization model for the multi-phase-field case based on the partial rank-one relaxation of the elastic energy following the previous works for the dual-phase materials \cite{BorEngMoslSchStei:2015:ldf,MosShcMon:2014:anh,KieFurMos:2017:anc,Sar:2018:cra,Sar:2021:mti}.  The models we introduce implicitly impose kinematic compatibility and static through controlled energy relaxation along the rank-one connection. This saves the complexity that arises from to rotating the coordinate system twice and splitting the mechanical fields into components corresponding to each jump condition when the jump conditions are fulfilled explicitly.  In contrast to the previous multi-phase-field models, the model derived here enforces all the interfacial jump conditions of all locally active pairs on their pairwise normals. A new phase-field energy definition is made by linearly interpolating the pairwise phase-field energies resulting from interactions with the other locally active phase-fields. This pairwise phase-field energy is assumed to be a function of a pairwise phase-field deformation gradient that accounts for the pairwise Hadamard jump condition (kinematic compatibility), and then the relaxation of the elastic energy defined on the pairwise interface with respect to the related pairwise  jump vector satisfies the pairwise static equilibrium. This jump vector represents the projection of the deformation gradient difference onto the normal of the interface.  In other words, the (bulk) effective energy is defined by the homogenization theory as the "volumetric" average of the active phase-field energies, and the phase-field energy is obtained by linearly interpolating the pairwise phase-field energies resulting from satisfying the pairwise jump conditions with the other active phase-fields. The model treats junction regions pairwise motivated by sharp interface literature \cite{Bha:2003:mom,BalJam:1992:pet,BalJam:1987:fpm,Bha:1993:cot,Bha:1992:sai,ZhaJamMue:2009:eba,KohMue:1992:bow,Mue:1999:vmf,JamHan:2000:mta,Dol:2004:vmf,BalCar:1999:ccf} since these regions arise due to the diffuse interface approximation and do not exist physically. 

The outline of the paper is as follows: In \Cref{sect:var} we derive the constitutive relation of the combined mechanical-phase transformation problem by evaluating the entropy inequality. We introduce in \Cref{sec:mpf} the multi-phase-field model the homogenization theory within the framework of multi-phase-field theory. The novel elasticity model is derived in \Cref{sec:mpfr1} before we compare the presented model with equal-strain and equal-stress assumptions in five numerical examples. We conclude our work in \Cref{sect:fin}. 
 All the numerical examples in this work were performed using the OpenPhase software \citep{Ope:2018:Pha,MedShcDarGlaZenPouSpaSte:2011:oao}. OpenPhase is a phase-field simulation library based on the multi-phase-field (MPF) model \citep{StePezNesSeePriSchRez:1996:apf,TiaNesDieSte:1998:tmf,StePez:1999:agf,NesBriGarHarStiBjo:2005:mas,Ste:2009:mas}.

%% file: Section2.tex
\section{The mechanical-phase transformation problem}
\label{sect:var}
The problem will first be separated into two subproblems:  the mechanical and the phase transformation subproblems. The constitutive laws of the combined problem will then be derived. Consider a homogeneous (single-phase) hyperelastic material under mechanical loading, where a material point within the reference configuration $\bX \in \Bnull$ is mapped to the deformed configuration $\Btime$ via the deformation mapping $\Bvarphi$. The deformation gradient is given as
$
\bF = \frac{\partial \Bvarphi}{\partial \bX} = \nabla_{\bX} \Bvarphi = \nabla \Bvarphi\,,
$
which maps an infinitesimal line element from the reference configuration to the deformed configuration \cite{Gur:1981:ait}. The local statement of entropy inequality (the Clausius-Duhem inequality) for isothermal processes is \cite{GurFriAna:2010:tma}
\begin{equation}
- \rho_0 \dot{\hat{\psi}} + J \boldsymbol{\sigma} : \bL  \geq 0\,,
\end{equation} 
where $\rho_0$ is the reference density, $\hat{\psi}$ is the specific Helmholtz free energy per unit mass, $J$ is the Jacobian (the determinant of deformation gradient), $\Bsigma$ is the Cauchy stress tensor, and $\bL$ is the spatial velocity gradient tensor. This entropy inequality can be reformulated in terms of the first Piola-Kirchhoff  stress tensor~$\bP$ by using the relation between the first Piola-Kirchhoff and Cauchy stress tensors $ \Bsigma= \frac{1}{J} \bP \cdot \bF^\textrm{T}\,$  along with  $\bL  = \frac{\partial \dot{\bx}}{\partial \bx} = \dot{\bF}  \cdot \bF^\textrm{-1}$. We obtain
\begin{equation}
- \rho_0 \dot{\hat{\psi}} + \bP    : \left[ \dot{\bF}  \cdot \bF^\textrm{-1} \cdot \bF \right] = - \rho_0  \dot{\hat{\psi}}+ \bP    : \dot{\bF} \geq 0\,. 
\end{equation} 
By introducing the volume-specific Helmholtz free energy $ \psi = \rho_0 {\hat{\psi}}$, we reach the final form of the entropy inequality
\begin{equation}
\label{eq:a}
 - \dot{\psi} + \bP  : \dot{\bF} =  \dot{\mathcal{W}} - \dot{\mathit{\psi}}  \geq 0\,,
\end{equation} 
where the term $ \dot{\mathcal{W}} = \bP   : \dot{\bF}$ represents the external power. \\

We extend the problem to a dual-phase material with a single order parameter (phase-field) $\phi \in \left[ 0,1\right]$. The micro-forces system associated with the phase transformation in an arbitrary control volume  ${\cal R} \subset \Bnull$  is introduced following \citep{Gur:1996:ggl,FriGur:1993:cto,FriGur:1994:dsst} by a stress vector $\Bxi$ and a scalar internal force $\Pi$. Here, this volume unit represents a material point with two phases having volume fractions $\phi$ and $1-\phi$.  The power associated with the stress vector is given by  $\int_{\partial {\cal R}} \dot{\phi} \; \Bxi \cdot \bn_{\partial {\cal R}} \; \textrm{d} A$ which represents the power transferred across the surface $\partial {\cal R}$ from and to neighboring volume units. The scalar internal force power equals  $\int_{ \cal R}  \Pi \dot{\phi} \; \textrm{d} V$ spent on ordering and moving  atoms inside the volume unit itself (changing the crystal structure). Here,  $\bn_{\partial {\cal R}}$ is the outward normal vector of the surface $\partial { \cal R}$. The local statement of the microstructure balance reads \cite{FriGur:1993:cto}
\begin{equation}
\label{eq:microbalance}
\Div \boldsymbol{\xi} +  \Pi = 0\,.
\end{equation} 
The Helmholtz free energy increase cannot exceed the external power for the micro-forces system. We get the following  local dissipation inequality  
\begin{equation}
\label{eq:b} 
  -  \dot{{\psi}} + \nabla\dot{\phi} \cdot \Bxi   - \Pi \dot{\phi}   \geqslant 0\,.   
\end{equation}
The local statement of entropy inequality for the coupled mechanical-phase transformation  problem is obtained by combining \Cref{eq:a,eq:b} as \cite{FriGur:1994:dsst,Gur:1996:ggl} 
\begin{equation} 
\label{eq:loc}
  -  \dot{{\psi}} + \bP    : \dot{\bF} + \nabla\dot{\phi} \cdot \Bxi   - \Pi \dot{\phi}   \geqslant 0\,.   
\end{equation}

The evaluation of the aforementioned entropy inequality will be done using the Coleman and Noll procedure \citep{ColNol:1963:tho}. In the following, we assume first that the Helmholtz free energy density depends on the state variables $\bF,\phi,\nabla\phi$ and $\dot{\phi}$. Substituting the time derivative of the Helmholtz free energy density in terms of the rates of its state variables, the dissipation inequality becomes
\begin{equation}
\label{eq:h}
{\cal D} = \left[ \bP - \frac{\partial \psi}{\partial \bF} \right]   : \dot{\bF}   + \left[ \Bxi - \frac{\partial \psi}{\partial \nabla \phi}  \right]  \cdot {{\nabla  \dot{\phi}}} - \left[ \Pi  + \frac{\partial \psi}{\partial \phi} \right]   \dot{\phi} -  \frac{\partial \psi}{\partial \dot{\phi}}   \ddot{\phi} \geq 0\,.
\end{equation} 
The entropy inequality is linear in $\dot{\bF}, {{\nabla \dot{\phi}}}$ and $\ddot{\phi}$ and the  accompanying terms must vanish to ensure that the inequality holds for all independent processes. Therefore, we get the following constitutive relations 
\begin{equation}
\label{eq:xi}
\bP = \frac{\partial \psi}{\partial \bF}\,, \qquad \Bxi = \frac{\partial \psi}{\partial \nabla \phi}\,, \qquad   \frac{\partial \psi}{\partial \dot{\phi}} = 0\,,
\end{equation} 
which means that the free energy is independent of interfacial velocity, i.e,  $\psi = \psi(\bF,\phi,\nabla\phi)$. The reduced dissipation inequality reads 
\begin{equation}
{\cal D}_{red} =\Delta G {\dot \phi} \geq 0\,,  \quad \textrm{with} \quad
\Delta G =   - \left[  \Pi  + \frac{\partial \psi}{\partial \phi} \right] \,,
\end{equation}
and to make the dissipation inequality fulfilled, we assume the dissipative driving force $\Delta G$ is proportional to the interface velocity
\begin{equation}
\Delta G =  m \dot{\phi}\,,
\end{equation}
where $m \geq 0$ is a constant related to the mobility of the interface (kinetic parameter). This dissipation is associated with rearranging the atoms and moving the interface.  Substituting the micro-stress balance form \Cref{eq:microbalance} and the constitutive micro-stress relation from \Cref{eq:xi} into the previous dissipative driving force, we obtain the following driving force using the functional derivative of the total energy $E$ 
\begin{equation}
\label{eq:i}
\Delta G = - \left[ \frac{\partial \psi}{\partial \phi}   - \Div(\Bxi) \right] = - \left[ \frac{\partial \psi}{\partial \phi} - \Div\left(\frac{\partial \psi}{\partial \nabla \phi}\right) \right] = - \frac{\delta E}{\delta \phi}\,.
\end{equation} 
 For the multi-phase-field case with $n$ locally non-vanishing phase-fields, we introduce a Lagrange multiplier $\lambda$ in the diffusion equations with the constraint $\sum_{\beta=1}^n \phi_\beta = 1$.  The equation of motion of phase-field $\alpha$ reads 
\begin{equation}
m \dot{\phi}_\alpha =  - \frac{\delta }{\delta \phi_\alpha}  \left( E - \lambda  \left[\sum_{\beta=1}^n \phi_\beta - 1 \right] \right) = - \left[ \frac{\delta E }{\delta \phi_\alpha} - \lambda   \right] , 
\end{equation} 
where the Lagrange multiplier is found to be $\lambda = \frac{1}{n} \sum_{\beta=1}^n  \frac{\delta E}{\delta \phi_\beta}    $    in \cite{StePez:1999:agf} leading to the following thermodynamically consistent diffusion equation
\begin{equation}
\label{eq:stepez}
m\dot{\phi}_\alpha =   \frac{1}{n}  \sum_{\beta=1}^n \left[\frac{\delta E}{\delta \phi_\beta}  - \frac{\delta E}{\delta \phi_\alpha} \right].
\end{equation}

%% file: Section3.tex
\section{Multi-phase-field elasticity model based on pairwise relaxation} 
\label{sec:mpf}
\subsection{Multi-phase-field model}
For most applications in materials science, dual-phase modeling is insufficient to describe the complexity of the microstructure.  The multi-phase-field model (MPF) \cite{ShcDuEngSte:2019:pfs,Ste:2009:mas,Ste:2013:pfm,SteApe:2006:mfp,BorDuStrBoeShcHarSte:2016:mdo,StePezNesSeePriSchRez:1996:apf,TiaNesDieSte:1998:tmf,StePez:1999:agf,SteSonHar:2010:pfm,BorEngBoeSchSte:2014:lse,BorEngMoslSchStei:2015:ldf,MedShcDarGlaZenPouSpaSte:2011:oao} utilizes an arbitrary number of phase-fields to indicate  the different phases or orientations of the same phase (variants and grains).  The diffusion equation reads as follows, in accordance with \Cref{eq:stepez}, 
\begin{equation}
\label{Eq:d} 
\dot{\phi}_\alpha =  \sum_{\beta=1}^{n} \frac{\pi^2 M_{\alpha\beta}}{4\eta n} \Delta G_{\alpha\beta}\,,
\end{equation}  
where $\eta \textrm{ and } M_{\alpha\beta}$ are, respectively, the interface width and the interfacial mobility  of $\alpha$-$\beta$ pairs \cite{Ste:2009:mas}. The pairwise driving force density $\Delta G_{\alpha\beta}$  accounts for the operative physical phenomena under consideration (interfacial, elastic, chemical, ... etc)
\begin{equation}
\Delta G_{\alpha\beta} = \Delta G^\textrm{int}_{\alpha\beta}  + \Delta G^\textrm{elas}_{\alpha\beta} +  \Delta G^\textrm{chem}_{\alpha\beta}  + .... = \frac{\delta E}{\delta \phi_\beta} -  \frac{\delta E}{\delta \phi_\alpha}\,,
\end{equation}
where the total free energy is obtained  by integrating the energy densities of the active phenomena over the entire considered domain  
\begin{equation} 
\label{Eq:b} 
 E = \int_{{\cal B}_0} \psi \;\textrm{d} V\,, \quad \textrm{with} \quad \psi = \psi^\textrm{int} + \psi^\textrm{elas} +  \psi^\textrm{chem}+... \;.
\end{equation}
Utilizing the functional derivative of the total energy $E$,  the pairwise driving force becomes
\begin{equation} 
\label{eq:gdf}
\Delta G_{\alpha\beta} = \frac{\partial \psi}{\partial \phi_\beta} - \frac{\partial \psi}{\partial \phi_\alpha} - \Div \left(\frac{\partial \psi}{\partial \nabla \phi_\beta} - \frac{\partial \psi}{\partial \nabla \phi_\alpha} \right)\,.
\end{equation}
The interfacial free energy for the MPF model is assumed as 
\begin{equation} 
\label{Eq:c} 
\psi^{\textrm{int}} = \sum_{\alpha=1}^{n-1}  \; \sum_{\beta=\alpha+1}^{n}   \frac{4\gamma_{\alpha\beta}}{\eta} \left[ - \frac{\eta^{2}}{\pi^{2}} \nabla \phi_{\alpha} \cdot \nabla \phi_{\beta} +   \phi_{\alpha}  \phi_{\beta} \right]\,,\\
\end{equation} 
where $\gamma_{\alpha\beta}$ is the interfacial energy per unit reference area between $\alpha$ and $\beta$ leading to the following diffusion equation for a dual-phase system ($\phi_\beta = 1 - \phi_\alpha$)  
\begin{equation} 
\label{eq:DeltaGtotal}
\dot{\phi}_\alpha =  M_{\alpha\beta}\left[ \gamma_{\alpha\beta} \left[ \nabla^2\phi_\alpha + \frac{\pi^2}{\eta^2} (\phi_\alpha - \frac{1}{2}) \right] + \frac{\pi^2 }{8 \eta} \left[ \Delta G^\textrm{elas}_{\alpha\beta} + \Delta G^\textrm{chem}_{\alpha\beta}  \right] \right] \,, 
\end{equation} 

Here, the Laplace operator $\nabla^2\phi_\alpha$ tries to widen (spread) the interface while the term  
$(\phi_\alpha - \frac{1}{2})$ acts in the opposite way and sharpens the interface. The mechanical and chemical driving force moves the interface, however, even for a constant driving force profile, it disturbs the traveling wave solution, i.e.  the phase-field steady state contour deforms for a moving interface, see  \cite{Ste:2009:mas}. We will return to this point later. 

\subsection{Homogenization theory within the multi-phase-field framework} 

We consider a material point $\cal R$ in the diffuse interface region with $n$ locally active phase-fields.   The mechanical behavior of the individual phase-fields is described by the elastic energy density of each individual phase-field  $\psi^\textrm{elas}_{i}$. The effective (bulk) elastic free energy density of the  material point $\cal R$ is obtained by the  Hill-Mandel condition \citep{Hil:1963:epo}, which states that the macroscopic strain energy density should be equal to the volumetric average of the microscopic strain energy density
\begin{equation}
\label{eq:hill}
\psi^\textrm{elas} = \frac{1}{V_{\cal R}} \int_{{\cal B}_{\cal R}} \psi^\textrm{elas}_i\, \textrm{d}V \, = \sum_{i=1}^n  \frac{V_i}{V_{\cal R}} \psi^\textrm{elas}_i  = \sum_{i=1}^n  \phi_i \psi^\textrm{elas}_i  \,,
\end{equation}
and similarly, the effective deformation gradient reads 
\begin{equation}
\bF = \frac{1}{V_{\cal R}} \int_{{\cal B}_{\cal R}} \bF_i\, \textrm{d}V \, =  \sum_{i=1}^n  \phi_i \bF_i \,.
\end{equation}
The effective first Piola-Kirchhoff stress is governed by the following constitutive relation
 \begin{equation}
\label{Eq:j} 
\bP = \frac{\partial \psi^\textrm{elas}}{\partial \bF}  =   \sum_{\alpha=1}^{n} \phi_\alpha \bP_\alpha \,, \quad \textrm{with} \quad \bP_\alpha =  \frac{\partial \psi^\textrm{elas}_\alpha}{\partial \bF_\alpha} \,.
\end{equation}

The remainder of this work focuses on small deformations and rotations. The effective infinitesimal strain and true stress tensors, both symmetric, are expressed as 

\begin{equation}
\label{eq:df}
\Bvarepsilon =  \sum_{i=1}^n  \phi_i \Bvarepsilon_i  \,, \quad  \Bsigma =  \sum_{i=1}^n  \phi_i \Bsigma_i \,, \quad \textrm{with} \quad  \Bsigma_i = \CC_i : (\Bvarepsilon_i - \Bvarepsilon_i^\textrm{B})  \,, 
\end{equation}

where $\CC_i$ and $\Bvarepsilon_i^\textrm{B}$ are the standard fourth-order elasticity tensor and the Bain strain (eigenstrain) of phase-field $i$, respectively. 

\begin{figure}[ht]
	\centering
	\unitlength=1mm
\begin{picture}(120,60)
	\graphicspath{{figures/}}
	 \put(0,5){\def\svgwidth{0.7\textwidth}{\small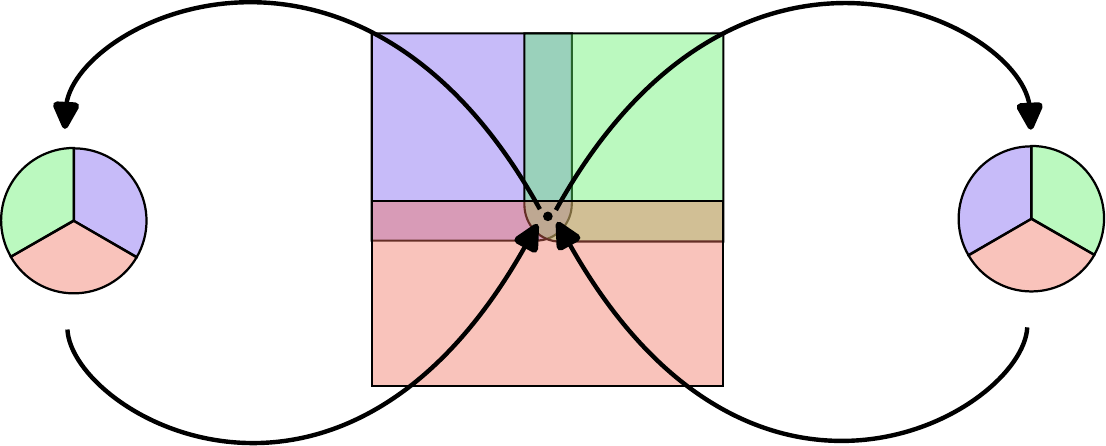}}
	\end{picture}
	\caption{The mechanical-phase field interactions; in order to solve the balance equation $\Div \Bsigma = \bzero $, the effective properties are needed in the interfacial regions (left). Once the mechanical problem is solved ($\Bvarepsilon,\Bsigma$ are known), the pairwise mechanical driving forces are calculated (right) in order to update the phase-fields.   }
	\label{fig:macromicro}
\end{figure}

The mechanical-phase field interaction occurs in two stages. This is illustrated in Figure \ref{fig:macromicro}. In the first interaction, we solve the mechanical problem, where the effective properties $\CC, \Bvarepsilon^B$ are needed in the interfacial regions. The influence of the effective properties in the interfacial regions is typically very small in the overall mechanical problem within $\cal B$. This aspect is beyond the scope of this paper. Therefore, for the numerical examples, we consider phase-fields with equal stiffness. The second interaction takes place after solving the "macroscopic" mechanical problem, where the mechanical driving forces are required, which contribute to the evolution of the microstructure. This, in turn, requires the mechanical response of the individual phase-field (localization). This is the central focus of this work.

\subsection{Equal-strain and equal-stress assumptions}

The equal-strain and equal-stress assumptions correspond to Voigt/Taylor and Reuss/Sachs limiting bounds of elasticity, cf. \cite{AmmAppCaiSam:2009:cpf,MosShcMon:2014:anh}. For the equal-strain assumption, we have $\Bvarepsilon_i = \Bvarepsilon$ for $i=1,...,n$, and the elastic energy reads
\begin{equation}
\psi^\textrm{elas,iso-strain} = \sum_{i=1}^n  \phi_i \psi^\textrm{elas}_i (\Bvarepsilon_i = \Bvarepsilon) \,.
\end{equation} 
This leads to the mechanical driving force 
  \begin{equation}
 \label{Eq:m}
 \Delta  G^\textrm{elas,iso-strain}_{\alpha\beta}  = \psi^\textrm{elas}_\beta -  \psi^\textrm{elas}_\alpha\,. 
 \end{equation}

For the equal stress assumption, we have $\Bsigma_i = \Bsigma$ for $i=1,...,n$. The strain of the phase-field $\alpha$ can be reformulated in terms of pairwise strain jumps with respect to all other locally active phase-fields $\llbracket \Bvarepsilon \rrbracket_{\alpha\beta} = \Bvarepsilon_{\beta} - \Bvarepsilon_ {\alpha}$.  \Cref{eq:df} can be rewritten as follows
\begin{equation}
\label{eq:DGJs}
\Bvarepsilon = \phi_\alpha \Bvarepsilon_\alpha + \sum_{\beta=1,\beta \neq \alpha}^n \phi_\beta \left[ \llbracket \Bvarepsilon \rrbracket_{\alpha\beta} + \Bvarepsilon_ {\alpha} \right]  \longrightarrow \Bvarepsilon_\alpha = \Bvarepsilon - \sum_{\beta=1,\beta \neq \alpha}^n \phi_\beta \llbracket \Bvarepsilon \rrbracket_{\alpha\beta}.  
\end{equation} 
The elastic energy density is given by 
\begin{equation}
\psi^\textrm{elas,iso-stress} = \sum_{i=1}^n  \phi_i \psi^\textrm{elas}_i (\Bvarepsilon_i) \,, 
\end{equation} 
and the corresponding mechanical driving force has the form 
   \begin{equation}
 \Delta  G^\textrm{elas,iso-stress}_{\alpha\beta}  = \psi^\textrm{elas}_\beta -  \psi^\textrm{elas}_\alpha- (\phi_\alpha + \phi_\beta) \, \Bsigma : \llbracket \Bvarepsilon \rrbracket_{\alpha\beta} \,, 
 \end{equation}
where the strain jump is derived from the assumption $\Bsigma_\alpha = \Bsigma_\beta$ as 
\begin{equation}
 \llbracket \Bvarepsilon \rrbracket_{\alpha\beta} =  [\CC_\beta^{-1} - \CC_\alpha^{-1}] : \Bsigma + [\Bvarepsilon_\beta^B -  \Bvarepsilon_\alpha^B ] \,. 
\end{equation}

%% file: figures/micro_macro.pdf_tex
\begingroup%
  \makeatletter%
  \providecommand\color[2][]{%
    \errmessage{(Inkscape) Color is used for the text in Inkscape, but the package 'color.sty' is not loaded}%
    \renewcommand\color[2][]{}%
  }%
  \providecommand\transparent[1]{%
    \errmessage{(Inkscape) Transparency is used (non-zero) for the text in Inkscape, but the package 'transparent.sty' is not loaded}%
    \renewcommand\transparent[1]{}%
  }%
  \providecommand\rotatebox[2]{#2}%
  \newcommand*\fsize{\dimexpr\f@size pt\relax}%
  \newcommand*\lineheight[1]{\fontsize{\fsize}{#1\fsize}\selectfont}%
  \ifx\svgwidth\undefined%
    \setlength{\unitlength}{530.48493153bp}%
    \ifx\svgscale\undefined%
      \relax%
    \else%
      \setlength{\unitlength}{\unitlength * \real{\svgscale}}%
    \fi%
  \else%
    \setlength{\unitlength}{\svgwidth}%
  \fi%
  \global\let\svgwidth\undefined%
  \global\let\svgscale\undefined%
  \makeatother%
  \begin{picture}(1,0.40290423)%
    \lineheight{1}%
    \setlength\tabcolsep{0pt}%
    \put(0,0){\includegraphics[width=\unitlength,page=1]{micro_macro.pdf}}%
    \put(-0.04,0.2){\color[rgb]{0,0,0}\makebox(0,0)[lb]{\smash{\large $\cal R$}}}%
    \put(1.01,0.2){\color[rgb]{0,0,0}\makebox(0,0)[lb]{\smash{\large $\cal R$}}}%
    \put(0.15,0.22){\color[rgb]{0,0,0}\makebox(0,0)[lb]{\smash{\large $\{\CC_i\}$}}}%
    \put(0.15,0.17){\color[rgb]{0,0,0}\makebox(0,0)[lb]{\smash{\large $\{\Bvarepsilon_i^\textrm{B}\}$}}}%
     \put(0.78,0.22){\color[rgb]{0,0,0}\makebox(0,0)[lb]{\smash{\large $\{\CC_i\}$}}}%
    \put(0.78,0.17){\color[rgb]{0,0,0}\makebox(0,0)[lb]{\smash{\large $\{\Bvarepsilon_i^\textrm{B}\}$}}}%
    \put(0.2,0.43){\color[rgb]{0,0,0}\makebox(0,0)[lb]{\smash{\large $\{\phi_i\}$}}}%
    \put(0.7,0.431){\color[rgb]{0,0,0}\makebox(0,0)[lb]{\smash{\large $\{\phi_i\}, \Bsigma, \Bvarepsilon$}}}%
     \put(0.2,-0.04){\color[rgb]{0,0,0}\makebox(0,0)[lb]{\smash{\large $\CC,\varepsilon^\textrm{B}$}}}%
    \put(0.2,-0.04){\color[rgb]{0,0,0}\makebox(0,0)[lb]{\smash{\large $\CC,\varepsilon^\textrm{B}$}}}%
      \put(0.02,0.22){\color[rgb]{0,0,0}\makebox(0,0)[lb]{\smash{\large $\phi_\alpha$}}}%
      \put(0.08,0.22){\color[rgb]{0,0,0}\makebox(0,0)[lb]{\smash{\large $\phi_\beta$}}}%
      \put(0.05,0.16){\color[rgb]{0,0,0}\makebox(0,0)[lb]{\smash{\large $\phi_\gamma$}}}%
        \put(0.89,0.22){\color[rgb]{0,0,0}\makebox(0,0)[lb]{\smash{\large $\phi_\alpha$}}}%
      \put(0.95,0.22){\color[rgb]{0,0,0}\makebox(0,0)[lb]{\smash{\large $\phi_\beta$}}}%
      \put(0.92,0.16){\color[rgb]{0,0,0}\makebox(0,0)[lb]{\smash{\large $\phi_\gamma$}}}%
     \put(0.73,-0.04){\color[rgb]{0,0,0}\makebox(0,0)[lb]{\smash{\large $\{\Delta G^\textrm{elas}_{ij}\}$}}}%
      \put(0.35,0.24){\color[rgb]{0,0,0}\makebox(0,0)[lb]{\smash{\large $ \alpha$}}}%
      \put(0.6,0.24){\color[rgb]{0,0,0}\makebox(0,0)[lb]{\smash{\large $ \beta$}}}%
      \put(0.48,0.07){\color[rgb]{0,0,0}\makebox(0,0)[lb]{\smash{\large $ \gamma$}}}%
  \end{picture}%
\endgroup%

%% file: Section4.tex
\section{Multi-phase-field elasticity model based on pairwise  relaxation (MPFR1 model)}   
\label{sec:mpfr1}
Let us assume we have a material point in a triple junction region with three locally active phase-fields with order parameters $\phi_1,\phi_2$, and $\phi_3$, the strains of the phase-fields $1$ and $2$ are given in \Cref{eq:DGJs} as
\begin{equation}
 \begin{aligned}
 \Bvarepsilon_1 &= \Bvarepsilon - \phi_2  \llbracket \Bvarepsilon \rrbracket_{12} - \phi_3 \llbracket \Bvarepsilon \rrbracket_{13} \\
  \Bvarepsilon_2 &= \Bvarepsilon + \phi_1  \llbracket \Bvarepsilon \rrbracket_{12} - \phi_3 \llbracket \Bvarepsilon \rrbracket_{23} \,.
 \end{aligned}
 \end{equation}
The  strain jump across the interface between the phase-field pair  ($1,2$)  can be obtained by subtracting    
 \begin{equation}
 \Bvarepsilon_2 - \Bvarepsilon_1 = \left[ \phi_1 + \phi_2 \right]  \llbracket \Bvarepsilon \rrbracket_{12} - \phi_3 \left[   \llbracket \Bvarepsilon \rrbracket_{23}  -  \llbracket \Bvarepsilon \rrbracket_{13} \right] \,,
\end{equation}
but we also have  $ \Bvarepsilon_2 - \Bvarepsilon_1 - \llbracket \Bvarepsilon \rrbracket_{12} = \bzero$ which leads to the necessary condition $   \llbracket \Bvarepsilon \rrbracket_{23}  -  \llbracket \Bvarepsilon \rrbracket_{13}  = -  \llbracket \Bvarepsilon \rrbracket_{12}$.  This means two jumps between two pairs of phase-fields can be freely chosen (e.g. relaxation along rank-one connection) but the third jump cannot be independently chosen (i.e. overdetermined system of equations).   Generalizing to $n$ locally-active phase-fields, the kinematic compatibility and static equilibrium conditions can be satisfied  for $n-1$ interfaces simultaneously from the existing $n(n-1)/2$ pairs. The idea to overcome this issue is motivated by the structure of the MPF  model which deals with pairwise driving forces. Thus, we will define pairwise elastic energies $\psi_{\alpha \beta}$ as a function of the pairwise strains $\Bvarepsilon_{\alpha \beta}$  which are calculated by enforcing the rank-one relaxation on the related normal vector $\bn_{\alpha \beta}$.  Finally, we obtain the elastic pairwise driving forces $\Delta G^\textrm{elas}_{\alpha\beta}$. 

We define an elastic energy for each individual phase-field in the diffuse interface region utilizing the pairwise phase-field energies resulting from the interactions with all the other locally non-vanishing phase-fields. The elastic energy of the phase-field $\alpha$ is assumed as  
\begin{equation}
\label{eq:pfenergy}
\psi^\textrm{elas}_\alpha := \frac{\sum\limits_{\beta=1,\beta\neq\alpha}^n   \phi_\beta \psi^\textrm{elas}_{\alpha\beta}}{\sum\limits_{\beta=1,\beta\neq\alpha}^n  \phi_\beta} = \frac{\sum\limits_{\beta=1,\beta\neq\alpha}^n   \phi_\beta \psi^\textrm{elas}_{\alpha\beta}}{1 - \phi_\alpha},
\end{equation}
where $\psi^\textrm{elas}_{\alpha\beta}$ is the energy of phase-field $\alpha$ with respect to $\beta$ (pairwise phase-field elastic energy). This pairwise phase-field energy is assumed to be a function of the associated pairwise phase-field strain $\psi_{\alpha\beta}^\textrm{elas} (\Bvarepsilon_{\alpha\beta})$ with $\Bvarepsilon_{\alpha\beta}$ is the strain of the phase-field $\alpha$  with respect to phase-field $\beta$.  The effective energy in \Cref{eq:hill} becomes then 
\begin{equation}
\label{eq:ee}
\psi^\textrm{elas} = \sum_{\alpha=1}^n  \frac{\phi_\alpha}{1-\phi_\alpha} \sum_{\beta=1,\beta\neq\alpha}^n   \phi_\beta \psi^\textrm{elas}_{\alpha\beta} \,.
\end{equation}
Note that, the previous treatment is active for $n \ge 3$ and has no effect for dual-phase-field regions. 

The rank-one relaxation (or convexification) of the free energy has been conventionally used to predict the formation of the sequential laminates  \cite{AubFagOrt:2003:acs,BarKieBukMen:2015:ake,BarCarHacHop:2004:erf,CarHacMie:2002:ncp,Lus:1996:otc,MieLam:2003:ats,OrtRep:1999;nem,KruLus:2003:tco,AraPed:2001:otc}. In the context of phase-field modeling, a rank-one convexification model \cite{BorEngMoslSchStei:2015:ldf,MosShcMon:2014:anh,KieFurMos:2017:anc} was introduced for dual-phase-field materials by enforcing rank-one connection (Hadamard condition) between the strains (or deformation gradients) of the two phase-fields. Then, a "partial" energy relaxation with respect to the  jump vector  (in normal direction) has been shown to satisfy the static equilibrium of the interface. Alternatively, a "full" rank-one convexification was introduced to obtain the interface normal vector as an additional field independent of the phase-field gradient, which is reasonable in sharp interface models but not suitable for the phase-field method. 

For the multi-phase-field case, the previously defined phase-field energy in \Cref{eq:pfenergy} enables us to enforce the partial rank-one relaxation pairwise on the interfaces of the active pairs. Starting from \Cref{eq:DGJs}, the deformation gradient of phase-field $\alpha$ can be rewritten as
\begin{equation}
\Bvarepsilon_\alpha = \Bvarepsilon - \phi_\beta \llbracket \Bvarepsilon \rrbracket_{\alpha\beta} - \sum_{i\ne\alpha,\beta}^n \phi_i \llbracket \Bvarepsilon \rrbracket_{\alpha i}\,,
\end{equation}
while retaining the strain jump of the phase-field $\alpha$ with respect to phase-field $\beta$. We eliminate the strain jumps  of all pairs except the pair ($\alpha,\beta$) by substituting the related strains 
\begin{equation}
 \Bvarepsilon_\alpha =   \Bvarepsilon - \phi_\beta \llbracket \Bvarepsilon \rrbracket_{\alpha\beta} - \sum_{i\ne\alpha,\beta}^n \phi_i \left[ \Bvarepsilon_i - \Bvarepsilon_\alpha \right] \longrightarrow \Bvarepsilon_{\alpha \beta} := \frac{1}{\phi_\alpha + \phi_\beta} \left[\Bvarepsilon  - \phi_\beta \llbracket \Bvarepsilon \rrbracket_{\alpha\beta} - \sum_{i=1,i\neq\alpha,\beta}^n \phi_i \Bvarepsilon_i \right],
 \end{equation}
which defines the phase-field $\alpha$ strain interacting with phase-field $\beta$ (pairwise phase-field strain) as a function of the strain jump between this pair.  The pairwise strain of the phase-field $\alpha$ with respect to $\beta$ and vice versa can be reformulated as
 \begin{equation}
  \begin{aligned}
 \label{eq:bdg}
 \Bvarepsilon_{\alpha \beta} = \frac{1}{\phi_\alpha + \phi_\beta} \left[\Bvarepsilon  - \phi_\beta \left[ \ba_{\alpha\beta} \otimes \bn_{\alpha\beta} \right]^\textrm{sym}  - \sum_{i=1,i\neq\alpha,\beta}^n \phi_i \Bvarepsilon_i \right],  \\
 \quad \Bvarepsilon_{\beta \alpha} = \frac{1}{\phi_\alpha + \phi_\beta} \left[\Bvarepsilon  + \phi_\alpha \left[ \ba_{\alpha\beta} \otimes \bn_{\alpha\beta} \right]^\textrm{sym}  - \sum_{i=1,i\neq\alpha,\beta}^n \phi_i \Bvarepsilon_i \right],
  \end{aligned}
 \end{equation} 
by introducing a rank-one connection (Hadamard jump condition)
 \begin{equation}
\llbracket \Bvarepsilon \rrbracket_{\alpha\beta} := \Bvarepsilon_{\beta\alpha} -  \Bvarepsilon_{\alpha\beta} = \left[ \ba_{\alpha\beta} \otimes \bn_{\alpha\beta} \right]^\textrm{sym} \,,
 \end{equation} 
 where $\ba_{\alpha\beta}$ and $\bn_{\alpha\beta}$ are, respectively, the strain jump and normal vectors of the phase-field pair $\alpha$ and $\beta$. The  pairwise normal vector is assumed in the form  \cite{SchSchSchReiHerSelBoeNes:2017:ots,MoeBlaWol:2008:qao}
\begin{equation}
\bn_{\alpha\beta} := \frac{ \nabla\phi_\beta -  \nabla \phi_\alpha}{ \| \nabla\phi_\beta -  \nabla \phi_\alpha \| } \,.
\end{equation}
 We define the elastic energy of the pair $\alpha$ and $\beta$ using the mixture rule of the relevant pairwise phase-field energies as
\begin{equation}
\label{eq:energyalphabeta}
\psi^\textrm{elas}_{\alpha,\beta} := \frac{\phi_\alpha}{\phi_\alpha + \phi_\beta} \psi^\textrm{elas}_{\alpha\beta} (\Bvarepsilon_{\alpha\beta}) + \frac{\phi_\beta}{\phi_\alpha + \phi_\beta} \psi^\textrm{elas}_{\beta\alpha} (\Bvarepsilon_{\beta\alpha})\,. 
\end{equation}
The strain jump vector $\ba_{\alpha\beta}$ is obtained by the relaxation of the previous pairwise elastic energy with respect to it 
 \begin{equation}
 \label{eq:pr1}
{\delta}_{ \ba_{\alpha\beta}}{\psi} ^\textrm{elas}_{\alpha,\beta} =  \frac{\phi_\alpha\phi_\beta}{[\phi_\alpha + \phi_\beta]^2}  \left[   -  \underbrace{\frac{\partial\psi^\textrm{elas}_{\alpha\beta}}{\partial \Bvarepsilon_{\alpha\beta}}}_{\Bsigma_{\alpha \beta}} \cdot \bn_{\alpha\beta} + \underbrace{\frac{\partial  \psi^\textrm{elas}_{\beta \alpha}}{\partial \Bvarepsilon_{\beta \alpha}}}_{\Bsigma_{\beta\alpha}} \cdot \bn_{\alpha\beta} \right] \cdot \delta \ba_{\alpha\beta}= 0\,,
\end{equation} 
and the static equilibrium of the dual interface is achieved $\Bsigma_{\beta\alpha} \cdot \bn_{\alpha\beta} =  \Bsigma_{\alpha\beta} \cdot \bn_{\alpha\beta}$. The elastic energy for the pairwise phase-field energy has the standard Hook's law form 
\begin{equation}
\psi^\textrm{elas}_{\alpha\beta} = \frac{1}{2} \left[ \Bvarepsilon_{\alpha\beta} -  \Bvarepsilon_\alpha^\textrm{B} \right] : \IC_{\alpha} :  \left[ \Bvarepsilon_{\alpha\beta} -  \Bvarepsilon_\alpha^\textrm{B} \right]\,. 
\end{equation} 

A closed-form solution is available for the pairwise strain jump vector $\ba_{\alpha\beta}$ by the pairwise energy relaxation in \Cref{eq:energyalphabeta} with respect to the jump vector after lengthy but straightforward calculations     
\begin{equation}
\label{eq:aalphabeta:1}
 \begin{aligned}
\ba_{\alpha\beta} =  & - \Big[ \bn_{\alpha\beta} \cdot [\phi_\beta \IC_\alpha +\phi_\alpha \IC_\beta] \cdot \bn_{\alpha\beta} \Big]^{-1}  \cdot  \Big[  \IC_\beta : \left[ \Bvarepsilon -\sum_{i=1,i\neq\alpha,\beta}^n  \phi_i \Bvarepsilon_i - (\phi_\alpha+\phi_\beta) \Bvarepsilon_\beta^\textrm{B} \right]   \\ & - 
\ \IC_\alpha : \left[ \Bvarepsilon - \sum_{i=1,i\neq\alpha,\beta}^n \phi_i \Bvarepsilon_i - (\phi_\alpha+\phi_\beta)  \boldsymbol{\varepsilon}_\alpha^\textrm{B} \right] \Big]  \cdot \; \bn_{\alpha\beta}  \;, 
 \end{aligned}
\end{equation} 
and the same solution can be obtained  starting from  $\Bsigma_{\beta\alpha} \cdot \bn_{\alpha\beta} =  \Bsigma_{\alpha\beta} \cdot \bn_{\alpha\beta}$. The strain jump vector simplifies as follows if the relevant phase-fields share the same stiffness 
\begin{equation}
\label{eq:aalphabeta:2}
\ba_{\alpha\beta} =  - \Big[ \bn_{\alpha\beta} \cdot \IC \cdot \bn_{\alpha\beta} \Big]^{-1}  \cdot \Big[  \IC : [ \Bvarepsilon_\alpha^\textrm{B} - \Bvarepsilon_\beta^\textrm{B}] \Big]  \cdot \; \bn_{\alpha\beta}  \;.  
\end{equation} 
To calculate the pairwise driving force of the phase-field pair $\alpha$ and $\beta$, we rewrite the effective energy in \Cref{eq:ee} isolating the possible pairwise phase-field energies involving  $\alpha$ or $\beta$ 
\begin{equation}
\label{eq:edmdf}
 \begin{aligned}
\psi^\textrm{elas} =& \frac{\phi_\alpha   \phi_\beta}{1-\phi_\alpha}  \psi^\textrm{elas}_{\alpha\beta} + \frac{\phi_\alpha   \phi_\beta}{1-\phi_\beta}  \psi^\textrm{elas}_{\beta\alpha} + \frac{\phi_\alpha}{1- \phi_\alpha}  \sum_{i\neq\alpha,\beta}^n \phi_i \psi^\textrm{elas}_{\alpha i}  + \frac{\phi_\beta}{1- \phi_\beta}  \sum_{i\neq\alpha,\beta}^n \phi_i \psi^\textrm{elas}_{\beta i} \\ &+ \phi_\alpha \sum_{i \neq \alpha,\beta}^n \frac{\phi_i}{1-\phi_i} \psi^\textrm{elas}_{i \alpha } + \phi_\beta \sum_{i \neq \alpha,\beta}^n \frac{\phi_i}{1-\phi_i} \psi^\textrm{elas}_{i \beta }  + \sum_{i \neq \alpha,\beta}^n  \frac{\phi_i}{1-\phi_i} \sum_{j \neq \alpha,\beta,i}^n \phi_j \psi^\textrm{elas}_{i j} \,. 
 \end{aligned}
\end{equation}
The pairwise phase-field strain in \Cref{eq:bdg} depends on all non-vanishing order parameters resulting in a significant complexity to calculate the partial derivatives in \Cref{eq:gdf} where all pairwise phase-field energies and strains are involved. Seeking simplification, we assume when we calculate the pairwise driving force of a phase-field pair that other phase-fields have a strain equal to the bulk strain.  Under this assumption, the pairwise strain of a phase-field $i$ with respect to a phase-field $j$ can be reformulated from  \Cref{eq:bdg} as
 \begin{equation}
 \label{eq:sim}
 \Bvarepsilon_{i j} = \Bvarepsilon - \frac{\phi_j}{\phi_i + \phi_j} \left[ \ba_{ij} \otimes \bn_{ij} \right]^\textrm{sym}  \,.
 \end{equation}
The mechanical driving force of the phase-field pair $\alpha$ and $\beta$ is obtained by the partial derivatives of the elastic energy density in \Cref{eq:edmdf} with respect to the order parameters $\alpha$ and $\beta$ taking into consideration the previous simplified pairwise phase-field strain in \Cref{eq:sim}
 \begin{equation}
 \label{eg:DVgeneral}
 \begin{aligned}
\Delta G^\textrm{elas}_{\alpha\beta} = & - \Delta G^\textrm{elas}_{\beta\alpha}  = (\frac{\partial }{\partial \phi_\beta} - \frac{\partial }{\partial \phi_\alpha})  \psi^{\textrm{elas}} 
\\ = &  \,  \frac{\phi_\alpha[1-\phi_\alpha] -\phi_\beta}{[1-\phi_\alpha]^2} \psi_{\alpha\beta}^\textrm{elas} + \frac{\phi_\alpha - [1-\phi_\beta]\phi_\beta}{[1-\phi_\beta]^2} \psi_{\beta\alpha}^\textrm{elas} \\
 & + \sum_{i\neq\alpha,\beta}^{n} \phi_i \left[ \frac{\psi_{\beta i}^\textrm{elas}}{[1-\phi_\beta]^2} - \frac{\psi_{\alpha i}^\textrm{elas}}{[1-\phi_\alpha]^2} \right] 
+ \sum_{i\neq\alpha,\beta}^n \frac{\phi_i}{1-\phi_i} [ \psi_{i\beta}^\textrm{elas} - \psi_{i\alpha}^\textrm{elas}] \\ &
-   \frac{\phi_\alpha \phi_\beta}{\phi_\alpha + \phi_\beta}  \Big[  \frac{\Bsigma_{\alpha\beta}}{1 - \phi_\alpha}  + \frac{\Bsigma_{\beta\alpha}}{1 - \phi_\beta}   \Big] : [\ba_{\alpha \beta} \otimes \bn_{\alpha \beta}] \\ & +   \,
 \sum_{i \neq \alpha,\beta}^n \frac{\phi_i^2\phi_\beta}{[\phi_i + \phi_\beta]^2}  \Big[  \frac{\Bsigma_{\beta i}}{1 - \phi_\beta}  + \frac{\Bsigma_{i \beta}}{1 - \phi_i}   \Big] : [\ba_{\beta i} \otimes \bn_{\beta i}]\\ & 
- \sum_{i \neq\alpha, \beta}^n \frac{\phi_i^2\phi_\alpha}{[\phi_i + \phi_\alpha]^2}  \Big[  \frac{\Bsigma_{\alpha i}}{1 - \phi_\alpha}  + \frac{\Bsigma_{i \alpha}}{1 - \phi_i}   \Big] :  [\ba_{\alpha i} \otimes \bn_{\alpha i}] \,, 
\end{aligned}
\end{equation}
which is reduced for dual-phase-field case to 
\begin{equation}
\Delta G^\textrm{elas}_{\alpha\beta} = \psi_\beta - \psi_\alpha - (\phi_\alpha \Bsigma_\alpha + \phi_\beta \Bsigma_\beta)  :  [\ba_{\alpha \beta} \otimes \bn_{\alpha \beta}]\,.
\end{equation}
Here, the previous mechanical driving force does not contain the term $(\frac{\partial }{\partial \nabla \phi_\beta} - \frac{\partial }{\partial \nabla \phi_\alpha})  \psi^{\textrm{elas}}$ which accounts for the difference of the micro-stress vectors between phase-fields $\alpha$ and $\beta$. This term represents the surface tension with a direction parallel to the normal of the interface for the isotropic materials \cite{FriGur:1994:dsst} which is expected to have a role for martensitic transformation during the martensitic nucleation in the austenite but not at planar interfaces between martensitic variants. However, due to the high complexity of considering this term in the framework of multi-phase-field we omit this term for this work.  For dual-phase-field case, the  mechanical driving force with surface tension effects reads 

\begin{equation}
\begin{aligned}
\label{eq:surfacetension}
\Delta G^\textrm{elas}_{\alpha\beta} = & \psi_\beta - \psi_\alpha - (\phi_\alpha \Bsigma_\alpha + \phi_\beta \Bsigma_\beta)  :  [\ba_{\alpha \beta} \otimes \bn_{\alpha \beta}]  \\ & - \Div \left(\ba_{\alpha \beta} \cdot [\Bsigma_\beta - \Bsigma_\alpha ] \cdot \dfrac{\bI - \bn_{\alpha \beta} \otimes \bn_{\alpha \beta}}{\|  \nabla \phi_\beta \|}\right) \,.
\end{aligned}
\end{equation}
 Notice that the surface tension term does not exist under equal-strain or equal-stress assumptions and appears only in the rank-one relaxation model because its formulation depends on the normal vector of the interface.

%% file: Section5.tex
 \section{Numerical examples} 
 \subsection{Plane interface between two variants: validation with the shape interface} 
  For the dual-phase-field case with a planar interface, the mechanical driving force can be reformulated from \Cref{eq:surfacetension} to 
  \begin{equation}
  \label{eq:dpdf}
  \Delta G^\textrm{elas}_{\alpha\beta} = \llbracket \psi^\textrm{elas} \rrbracket_{\alpha\beta} - \left[ \Bsigma  \cdot \bn_{\alpha\beta} \right] \cdot \left[ \llbracket \Bvarepsilon \rrbracket_{\alpha\beta} \cdot \bn_{\alpha\beta} \right],
  \end{equation} 
  and the diffusion equation referring to \Cref{eq:DeltaGtotal} simplifies to the following form in the absence of a chemical driving force 
\begin{equation} 
\dot{\phi}_\alpha = - \dot{\phi}_\beta =  \mu_{\alpha\beta}\left[ \gamma_{\alpha \beta} \left[ \nabla^2\phi_\alpha + \frac{\pi^2}{\eta^2} (\phi_\alpha - \frac{1}{2}) \right] + \frac{\pi^2 }{8 \eta} \left[ \llbracket \psi^\textrm{elas} \rrbracket_{\alpha\beta} - \left[ \Bsigma  \cdot \bn_{\alpha\beta} \right] \cdot \left[ \llbracket \Bvarepsilon \rrbracket_{\alpha\beta} \cdot \bn_{\alpha\beta} \right]  \right] \right] \,.
\end{equation} 
  We rewrite the previous equation, omitting the phase-field indices, as 
   \begin{equation}
b  \, \dot{\phi}=    \llbracket \psi^\textrm{elas} \rrbracket - \Bsigma  :  \llbracket \Bvarepsilon \rrbracket    + \gamma \kappa = \llbracket \psi^\textrm{elas} \rrbracket - \left[ \Bsigma  \cdot \bn \right] \cdot \left[ \llbracket \Bvarepsilon \rrbracket \cdot \bn \right]  + \gamma \kappa\,,
   \end{equation}
   where $\kappa =  \frac{8\eta}{\pi^2} \left[ \nabla^2\phi + \frac{\pi^2}{\eta^2} (\phi - \frac{1}{2}) \right] $ and  $b= \frac{8 \eta}{\pi^2   \mu}$ represents a generalized curvature and a kinematic parameter, respectively \cite{SteApe:2006:mfp}. For the sharp interface, the normal configurational force balance was introduced in \cite{Gur:2000:cfa,Gur:1995:tno} as 
      \begin{equation}
      \label{eq:SI}
b^*  \, \bv \cdot \bn =   \llbracket \psi^\textrm{elas} \rrbracket^* - \left< \Bsigma \right>^*  :  \llbracket \Bvarepsilon \rrbracket^*  + \gamma \kappa^* =  \llbracket \psi^\textrm{elas} \rrbracket^*- \left[ \left< \Bsigma \right>^* \cdot \bn \right] \cdot \left[ \llbracket \Bvarepsilon \rrbracket^* \cdot \bn \right]  + \gamma \kappa^* \,,
   \end{equation}
where $b^*$ is a kinematic modulus, $\kappa^*$ is the total curvature (twice the mean curvature), $\bv$ is the interfacial velocity vector with a projected normal component $\bv \cdot \bn$, $\llbracket \psi^\textrm{elas} \rrbracket^* = \psi^\textrm{elas,+} - \psi^\textrm{elas,-} $ is the elastic energy difference between the two sides of the interface,  $\llbracket \Bvarepsilon \rrbracket^* = \Bvarepsilon^+ - \Bvarepsilon^-$ is the strain jump across the sharp interface and  $\left< \Bsigma \right>^* = \left[\Bsigma^+ + \Bsigma^- \right]/2$ is the average stress acting on both sides of the sharp interface. We notice that the phase-field diffusion equation, coupled with the mechanical driving  of the MPFR1 model, follows the structure of the sharp interface mobility equation.

 As a three-dimensional validation numerical example, we set two martensitic variants with a planar interface whose outward normal vector is in the $z$-direction.  The phase-fields $\alpha$ and $\beta$ have the same stiffness but different eigenstrains corresponding to cubic-to-tetragonal transformation. The geometry is shown in Figure \ref{fig:geometryfirstexample:1} and the relevant parameters are presented in Table \ref{tab:example1}. The middle phase-field ($\beta$) has a thickness of half of the edge length so that both phase-fields have the same total volume. Under periodic boundary conditions, the elastic energies of both variants are equal and the system is at equilibrium. To induce the movement of the interface, we add an overall strain $ \Bvarepsilon^\textrm{external} = \textrm{diag} \left(0.01,0.0075,0.005\right)$, favoring one of the phase-fields.

\begin{table}[h]
    \centering
    \begin{tabular}{|c|c|} 
        \hline
        \textbf{parameter} & \textbf{value} \\ 
        \hline
        system size & $100  \times 100 \times 100$ grid points \\
        grid spacing:  $\Delta x$ & $10^{-7}$ m \\ 
        interface width: $\eta$ & $5$  $\Delta x$ \\ 
      Lam\'e constants & $\lambda = 80$ GPa and $\mu = 120$ GPa \\
   Bain strains & $  \Bvarepsilon_\alpha^\textrm{B} = \textrm{diag}\left(0.02,-0,01,-0,01\right)$,   $\Bvarepsilon_\beta^\textrm{B} = \textrm{diag}\left(-0.01,-0,01,0,02\right)$ \\
   \hline
       time step size: $\Delta t$  & $10^{-8}$ sec \\
      interfacial energy:  $\gamma$  & $0.1 \, \textrm{J}/\textrm{m}^2$ \\
       interfacial mobility: $M$  &  $3 \cdot 10^{-7} \, \textrm{m}^4 \,  / ( \textrm{J}\cdot \textrm{sec})$  \\
        \hline
    \end{tabular}
    \caption{Parameters used for plane interface example}
    \label{tab:example1}
\end{table}
 
The mechanical solution of each phase-field is homogeneous in the bulk region  ($\phi_\alpha  \phi_\beta  = 0$)  but not necessarily in the diffuse interface regions ($\phi_\alpha  \phi_\beta  \neq 0$). Only the MPFR1 model provides a homogeneous solution for both phase-fields, including in the diffuse interface region, while the equal-strain or equal-stress assumptions lead to varying fields along the interface normal vector, which is unphysical in the context of sharp interface mechanics. We demonstrate this behavior in   \Cref{fig:geometryfirstexample:2}. The elastic energy of phase-field $\alpha$ obtained by the MPFR1 model is constant in the diffuse interface regions, whereas it varies for equal-strain and equal-stress models. As expected, the effective energy of the MPFR1 model lies between the upper and lower bounds obtained by the equal-strain and equal-stress assumptions, respectively.

\begin{figure}[ht]
	\centering
	\unitlength=1mm
\subfigure[]{\begin{picture}(65,55)
	\graphicspath{{figures/}}
	 \put(0,3){\def\svgwidth{0.35 \textwidth}{\small\input{figures/Example_1_geo.eps_tex}}}
	 \label{fig:geometryfirstexample:1}
	\end{picture}} 
	\subfigure[]{\begin{picture}(90,46)
	\graphicspath{{figures/}}
	 \put(-5,0){\def\svgwidth{0.5\textwidth}{\small\input{figures/psi_df_1.eps_tex}}}
	  \label{fig:geometryfirstexample:2}
	\end{picture}} 
	\caption{Example 1: two martensitic variants with plane interfaces. a) Problem dimensions  in grid points with cut-through showing the line of inspection in white. b) The elastic energy of phase-field $\alpha$ and the effective energy along inspection line using different mode. MPFR1 model give homogeneous solution for both phase-fields.  }
	\label{fig:geometryfirstexample}
\end{figure}

Now we focus on the mechanical driving force. We use the homogeneous fields in the bulk, far from the diffuse interface region, to calculate the driving force as if a sharp interface were present. The jumps of the strain and energy fields between are given by
\begin{equation}
\llbracket \psi^\textrm{elas} \rrbracket^*  = 7.0288547 \cdot 10^{7}  \, \textrm{J}/\textrm{m}^3\,,  \quad  \llbracket \Bvarepsilon \rrbracket^* =  \textrm{diag}\left(0,0,0.01714278\right)\,, \\
\end{equation}
and the stress fields equal 
\begin{equation}
\begin{aligned}
  \Bsigma_{\alpha} =  \textrm{diag}\left(0.071435,4.47143,2.7\right)\cdot10^9 \, \textrm{N}/\textrm{m}^2 \,, \\  \Bsigma_{\beta} = \textrm{diag}\left(6.92857,6.52857,2.7\right)\cdot10^9 \, \textrm{N}/\textrm{m}^2\,.
\end{aligned}
\end{equation}
The jump in the strain field projected on the interface plane ($xy$-plane) vanishes, satisfying kinematic compatibility. The stresses projection onto the interface normal vector ($z$-direction) are equal, satisfying static equilibrium. The mechanical driving force of the sharp interface following  \Cref{eq:SI} equals
\begin{equation}
\begin{aligned}
\Delta G^\textrm{SI} & =     \llbracket \psi^\textrm{elas} \rrbracket^*- \left[ \left< \Bsigma \right>^* \cdot \bn \right] \cdot \left[ \llbracket \Bvarepsilon \rrbracket^* \cdot \bn \right]  + \gamma \kappa^*  \\ & =   7.0288547
 \cdot 10^{7} - 2.7 \cdot 10^9  \cdot (0.01714278) + 0 = 2.4 \cdot 10^7 \, \textrm{J}/\textrm{m}^3\,, 
\end{aligned}
\end{equation} 
which is indeed obtained using the MPFR1 model as a constant value along the interface normal, see \Cref{fig:DG:1:PR1} at $t_0$. The driving force of equal-strain and equal-stress models varies across the interface.  We show  in Figure \ref{fig;1;DG} the profile of $\phi_\beta$ with the calculated driving force  at $t_0$ and at $t>t_0$.  For the equal-strain model, the driving force tends to narrow (sharpening) the interface, see \Cref{fig:DG:1:Voigt}, whereas for the equal-stress model, it widens (smearing) the interface as shown in \Cref{fig:DG:1:Reuss}.

\begin{figure}[ht]
	\centering
	\unitlength=1mm
\subfigure[MPFR1]{\begin{picture}(53,55)
	\graphicspath{{figures/}}
	 \put(0,3){\def\svgwidth{0.32 \textwidth}{\small\input{figures/Example_1_DG_PR1.pdf_tex}}}
	 \label{fig:DG:1:PR1}
	\end{picture}} 
\subfigure[equal-stress]{\begin{picture}(53,55)
	\graphicspath{{figures/}}
	 \put(0,3){\def\svgwidth{0.32 \textwidth}{\small\input{figures/Example_1_DG_Reuss.pdf_tex}}}
	 \label{fig:DG:1:Reuss}
	 \end{picture}} 
\subfigure[equal-strain]{\begin{picture}(53,55)
	\graphicspath{{figures/}}
	 \put(0,3){\def\svgwidth{0.32 \textwidth}{\small\input{figures/Example_1_DG_Voigt.pdf_tex}}}
	 \label{fig:DG:1:Voigt}
	 \end{picture}} 
	\caption{The value of the order-parameter of phase-field $\beta$ and the elastic driving force for the implemented models over the inspection line. The MPFR1 model is the most stable numerically and maintains the phase-field profile, whereas the equal-strain and equal-stress models destabilize it. }
	\label{fig;1;DG}
\end{figure}

To preserve the phase-field profiles, different actions are required. First, even a constant driving force interferes with the steady state wave solution, this is addressed by scaling the driving force with a contour function $h = (\textrm{constant}*\sqrt{\phi_\alpha \phi_\beta})$  motivated by the traveling wave solution for the double obstacle potential \cite{Ste:2009:mas}. For the remainder of this work, we will use the so-called antisymmetric approximation

\begin{equation} 
\dot{\phi}_\alpha = \sum_{\beta=1}^n M_{\alpha\beta}\left[ \gamma_{\alpha\beta} \left[ \phi_\beta \nabla^2 \phi_\alpha - \phi_\alpha \nabla^2 \phi_\beta  
+ \frac{\pi^2}{2 \eta^2} (\phi_\alpha - \phi_\beta) \right] + \frac{\pi}{\eta} \sqrt{\phi_\alpha \phi_\beta} \left[ \Delta G^\textrm{elas}_{\alpha\beta} + \Delta G^\textrm{chem}_{\alpha\beta}  \right] \right] \,. 
\end{equation} 

The second action is averaging the driving force to obtain a smooth "semi-constant" profile. The driving force for a point belonging to the diffuse interface will be averaged over nearby points that are within a distance less than a certain limit, which is chosen to be the interface's width for the work. Seeking comparison, these actions will be applied to all mechanical models. We show in \Cref{fig;2;DG} that these two actions indeed lead to a stable phase-field profile.

\begin{figure}[ht]
	\centering
	\unitlength=1mm
\subfigure[MPFR1]{\begin{picture}(53,55)
	\graphicspath{{figures/}}
	 \put(0,3){\def\svgwidth{0.32 \textwidth}{\small\input{figures/Example_1_DG_PR1_enhanced.pdf_tex}}}
	\end{picture}} 
\subfigure[equal-stress]{\begin{picture}(53,55)
	\graphicspath{{figures/}}
	 \put(0,3){\def\svgwidth{0.32 \textwidth}{\small\input{figures/Example_1_DG_Reuss_enhanced.pdf_tex}}}
	 \end{picture}} 
\subfigure[equal-strain]{\begin{picture}(53,55)
	\graphicspath{{figures/}}
	 \put(0,3){\def\svgwidth{0.32 \textwidth}{\small\input{figures/Example_1_DG_Voigt_enhanced.pdf_tex}}}
	 \end{picture}} 
	\caption{The value of the order-parameter of phase-field $\beta$ and the elastic driving force for the implemented models over the inspection line.  The stabilization actions success to preserve the phase-field profile. Dashed and solid green lines indicate to the mechanical driving force before and after averaging, respectively. }
\label{fig;2;DG}
\end{figure}

\FloatBarrier

\subsection{Growing martensitic nucleus in an austenitic matrix} 

\begin{figure}[ht!]
	\centering
	\unitlength=1mm
\subfigure[start geometry]{\begin{picture}(50,45)
	\graphicspath{{figures/}}
    \put(0.0,0.0){\includegraphics[width=0.32 \textwidth]{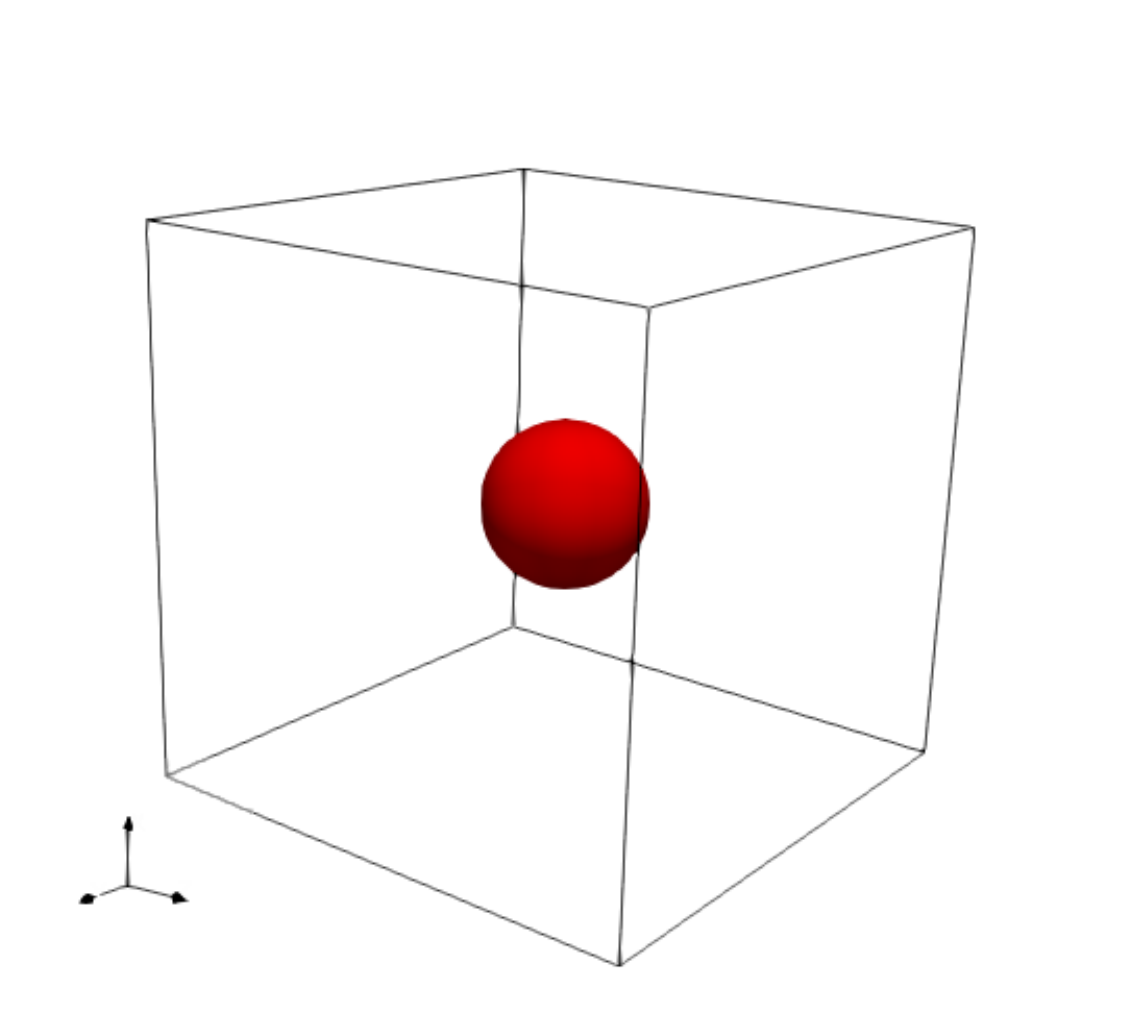}}
	\end{picture}\label{fig;2:start}} 
\subfigure[no mechanics]{\begin{picture}(50,45)
	\graphicspath{{figures/}}
    \put(0.0,0.0){\includegraphics[width=0.32 \textwidth]{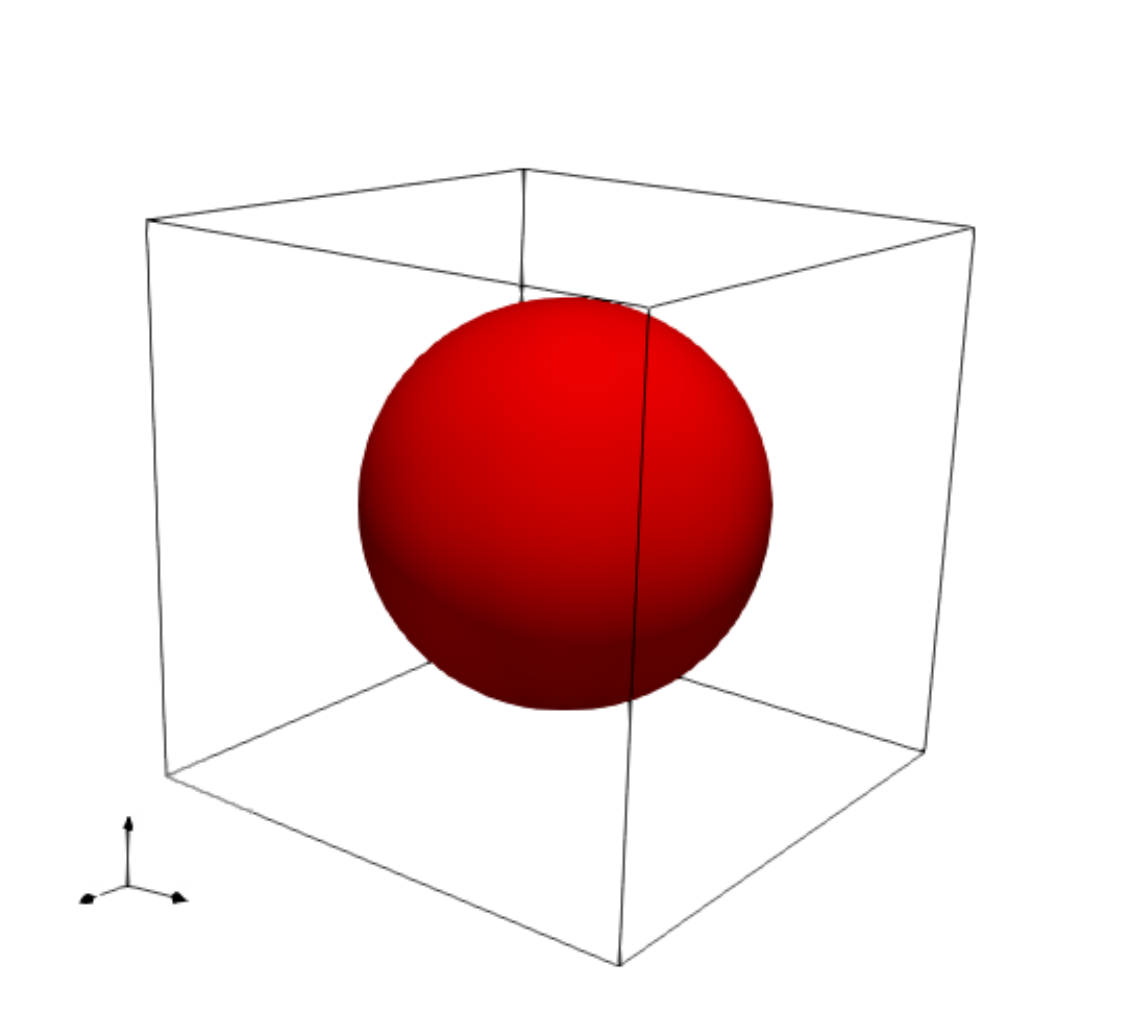}}
	\end{picture}\label{fig;2:nomechanics}} 
\subfigure[equal stress]{\begin{picture}(50,45)
	\graphicspath{{figures/}}
    \put(0.0,0.0){\includegraphics[width=0.32 \textwidth]{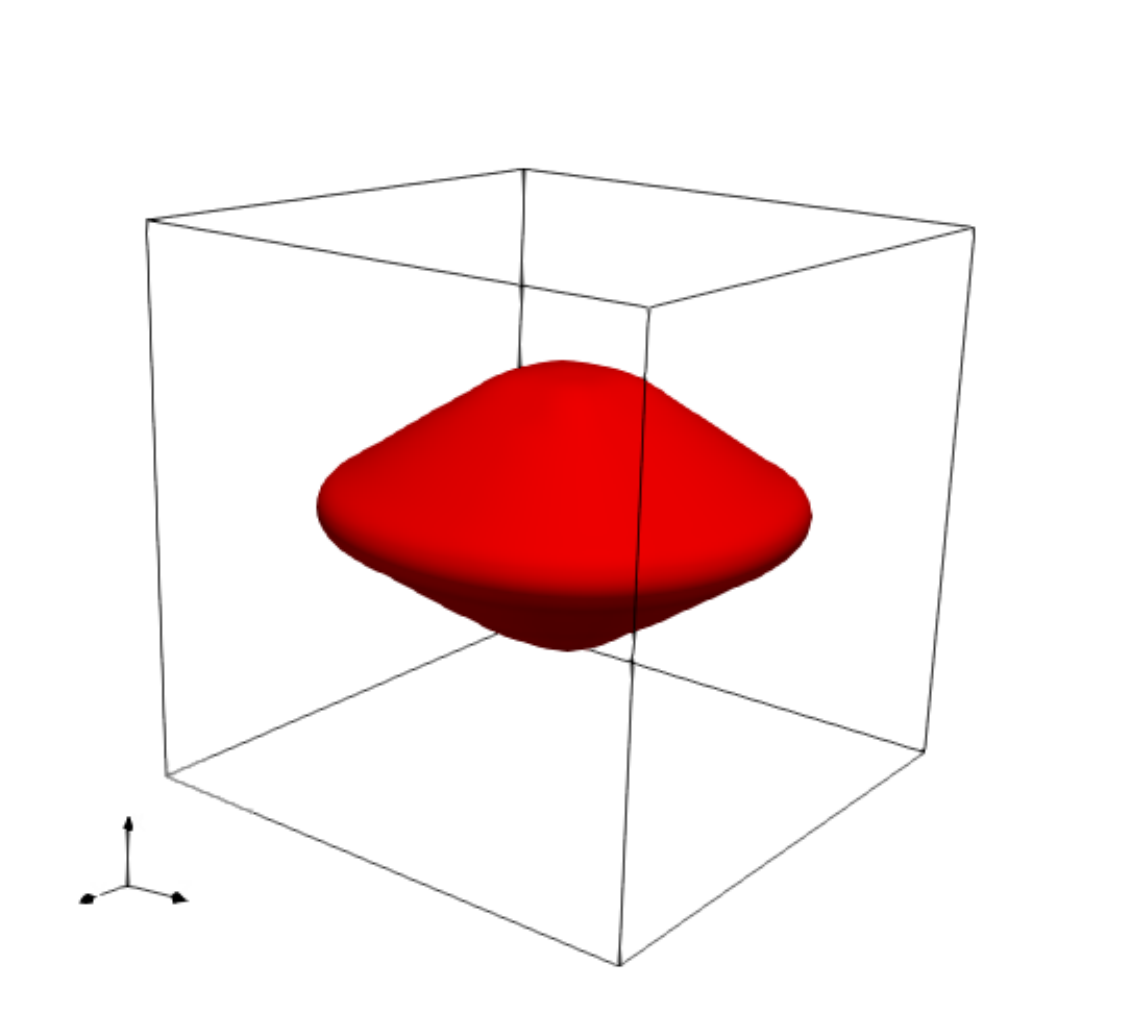}}
	\end{picture}\label{fig;2:Reuss}} 
\subfigure[equal strain]{\begin{picture}(50,45)
	\graphicspath{{figures/}}
    \put(0.0,0.0){\includegraphics[width=0.32 \textwidth]{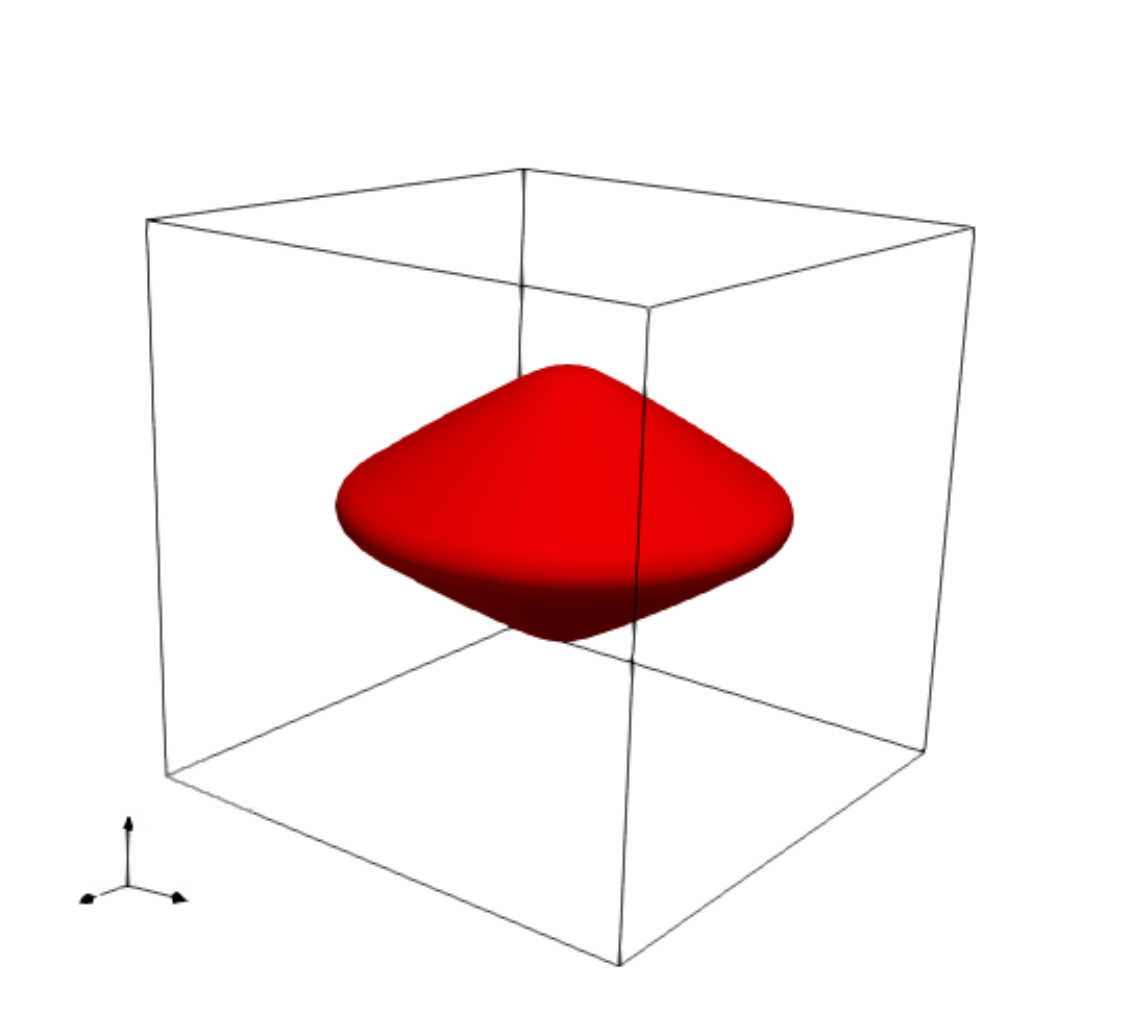}}
	\end{picture}\label{fig;2:Voigt}} 
\subfigure[MPFR1]{\begin{picture}(50,45)
	\graphicspath{{figures/}}
    \put(0.0,0.0){\includegraphics[width=0.32 \textwidth]{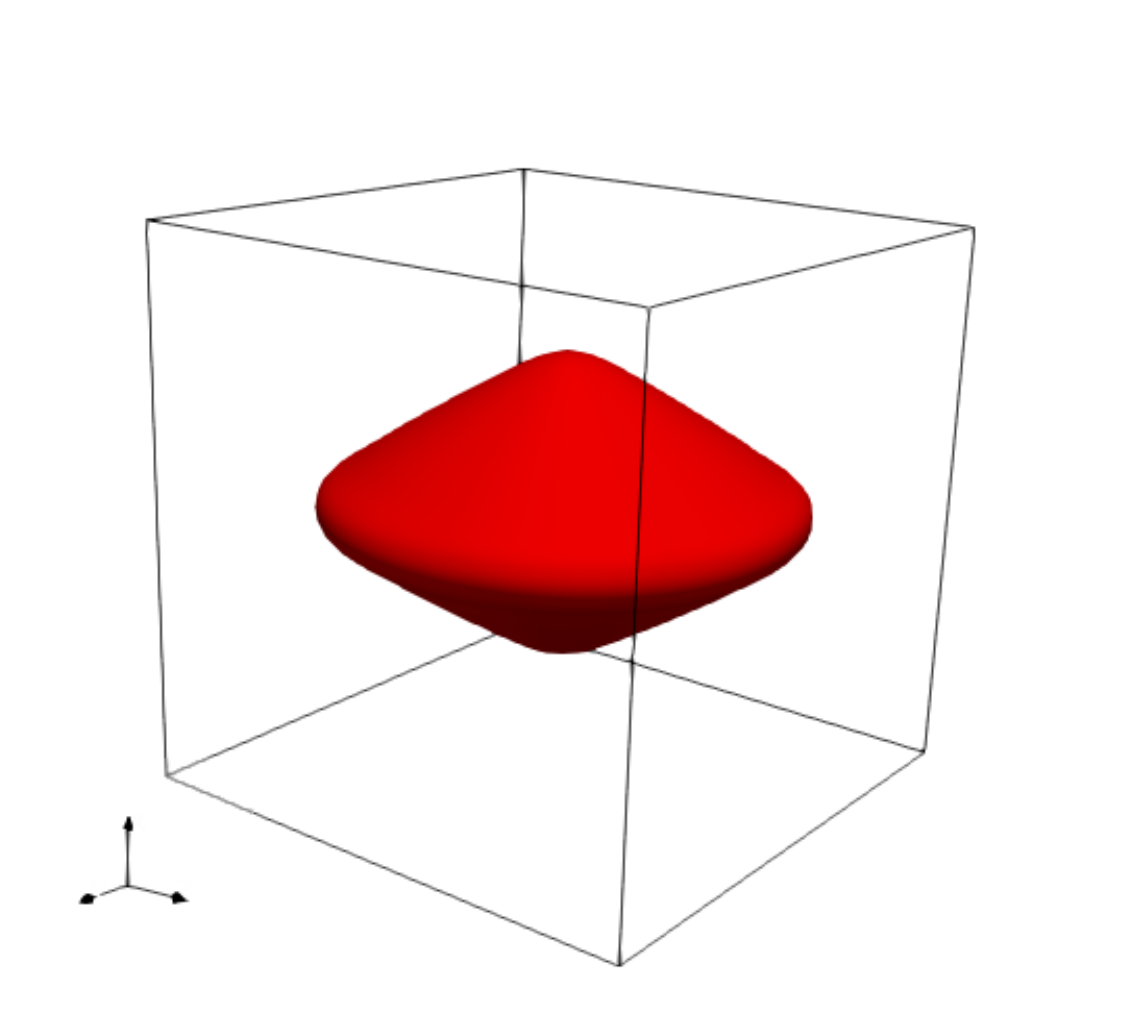}}
	\end{picture}\label{fig;2:MPFR1_without}} 
\subfigure[MPFR1 with surface tension]{\begin{picture}(50,40)
	\graphicspath{{figures/}}
    \put(0.0,0.0){\includegraphics[width=0.32 \textwidth]{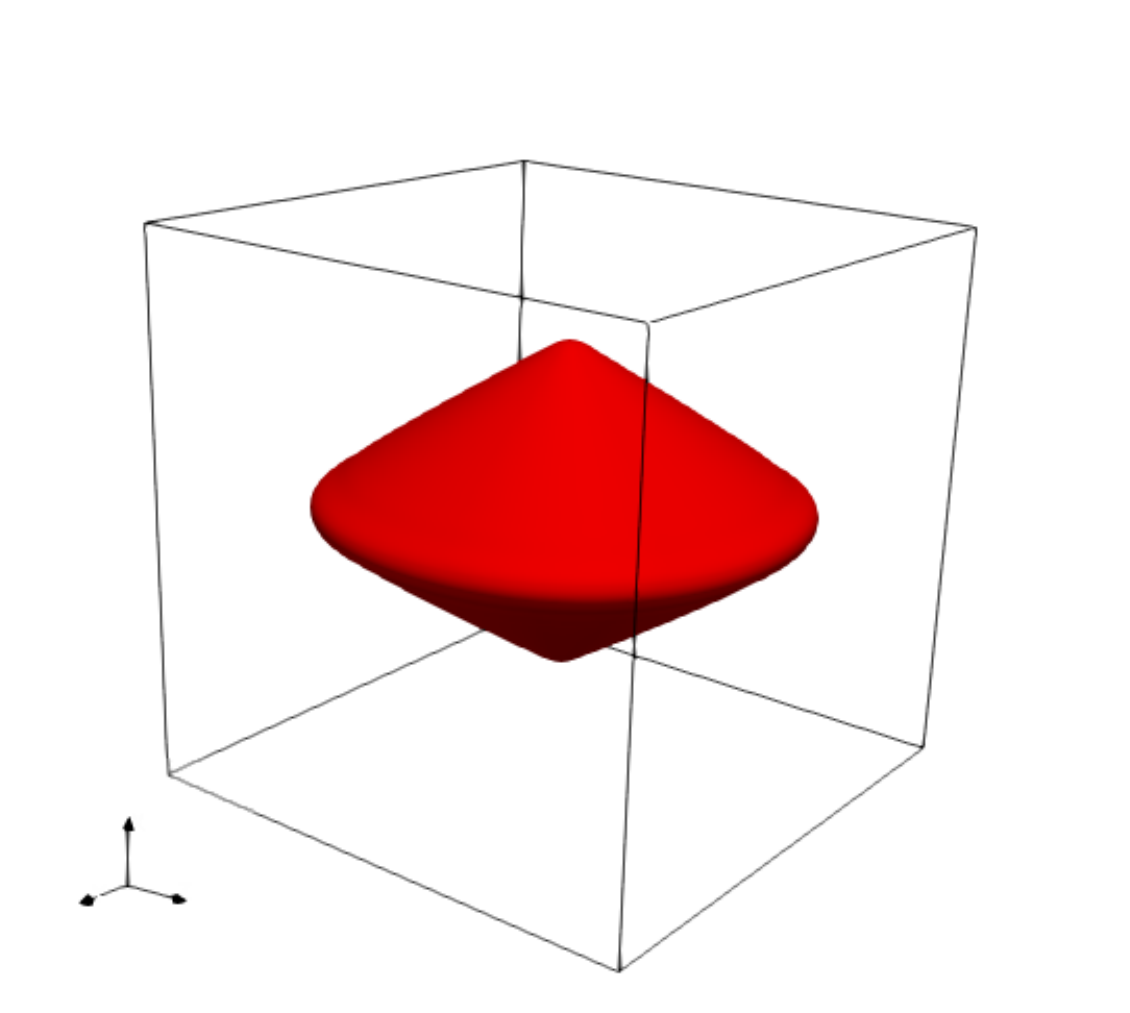}}
	\end{picture}\label{fig;2:MPFR1_with}} 
	\caption{The anisotropic growth of a martensitic nucleus using different mechanical models}
	\label{fig;2}
\end{figure}

A single nucleus (phase-field $\beta$) is initialized as a sphere with a radius $r$ at the center of the simulation box filled with phase-field $\alpha$, see Figure \ref{fig;2:start}. The used parameters are shown in Table \ref{tab:example2}. Periodic boundary conditions are applied to the mechanical solution. The nucleus has an eigenstrain corresponding to a cubic-to-tetragonal transformation leading to a mechanical stressing.  Furthermore, we introduce a constant chemical driving force to promote the growth of the nucleus.

\begin{table}[h]
    \centering
    \begin{tabular}{|c|c|} 
        \hline
        \textbf{parameter} & \textbf{value} \\ 
        \hline
        system size: & $100  \times 100 \times 100$ grid points \\
        grid spacing: $\Delta x$ & $10^{-8}$ m \\ 
        radius: $r$ & $15$  $\Delta x$ \\        
       interface width: $\eta$ & $5$  $\Delta x$ \\ 
       time step size: $\Delta t$  & $10^{-10}$ sec \\
      interfacial energy:  $\gamma$  & $0.1 \, \textrm{J}/\textrm{m}^2$ \\
       interfacial mobility: $M$  &  $3 \cdot 10^{-7} \, \textrm{m}^4 \,  / ( \textrm{J}\cdot \textrm{sec})$  \\
      Lam\'e constants & $\lambda = 80$ GPa and $\mu = 120$ GPa \\
   Bain strains &  $\Bvarepsilon_\beta^\textrm{B} = \textrm{diag}\left(-0.01,-0,01,0,02\right)$ \\
chemical driving force: $\Delta G_{\alpha \beta}^\textrm{chem}$ & $5.33 \cdot 10^7 \, \textrm{J}/\textrm{m}^3$ \\ 
        \hline
    \end{tabular}
    \caption{Parameters used for growing nucleus example.}
    \label{tab:example2}
\end{table}
 
We focus here on the shape of the growing martensitic nucleus for the different mechanical models shown in \Cref{fig;2}. The nucleus grows isotropically when the mechanical driving force is not considered, as shown in \Cref{fig;2:nomechanics}. For the MPFR1 model,  two cases are investigated: with and without the surface tension effects in the mechanical driving force (divergence term). The equal-stress  model leads to a different curvilinear shape, see \Cref{fig;2:Reuss}.  The equal-strain and MPFR1 models result in a similar  nucleus shape, as shown in \Cref{fig;2:Voigt,fig;2:MPFR1_without}, respectively. However, adding the divergence term in the mechanical driving force to MPFR1 model, which is corresponding to the surface tension, leads to a more angular and edged shape, see \Cref{fig;2:MPFR1_with}. Thus, the surface tension has a critical role in the nucleation and growth of martensite phase and the evolution of the microstructure.  

\subsection{Triple junction}
\begin{figure}[H]
	\centering
	\unitlength=1mm
\subfigure[Geometry with zoom in triple junction]{\begin{picture}(100,50)
	\graphicspath{{figures/}}
    \put(0,2){\includegraphics[width=0.6\textwidth]{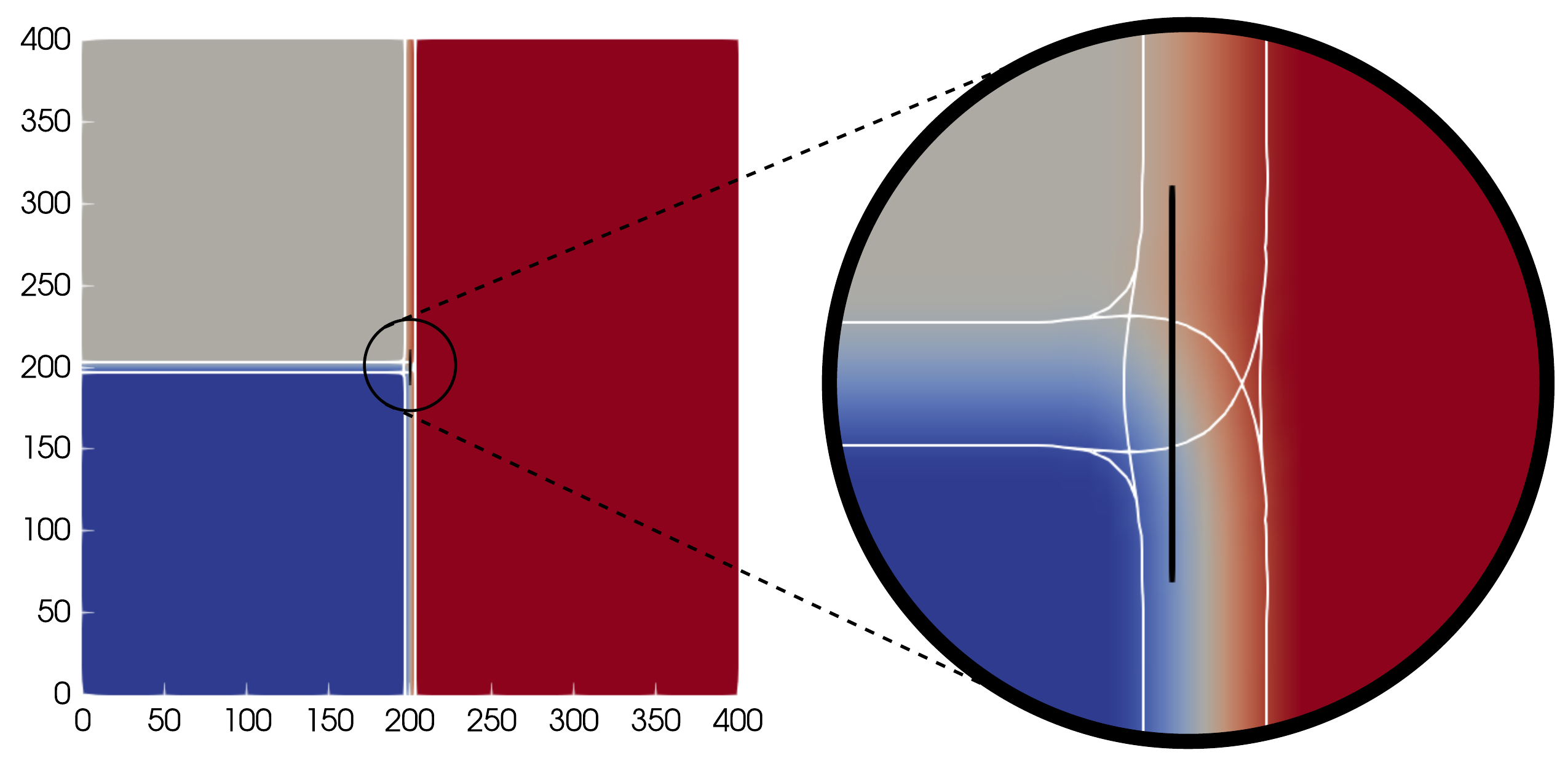}}
    \put(20.405323271,-0.055607354){\color[rgb]{0,0,0}\makebox(0,0)[lb]{\scriptsize \smash{$x \textrm{ in } 10^{-4} \textrm{ m}$}}}%
 \put(0.005323271,20.055607354){\color[rgb]{0,0,0}\rotatebox{90}{\makebox(0,0)[lb]{\scriptsize \smash{$y \textrm{ in } 10^{-4} \textrm{ m}$}}}}%
    \put(13.405323271,15.055607354){\color[rgb]{0,0,0}\makebox(0,0)[lb]{ \color{white} \smash{$\alpha$}}}%
    \put(13.405323271,35.055607354){\color[rgb]{0,0,0}\makebox(0,0)[lb]{ \smash{$\beta$}}}%
    \put(35.405323271,24.055607354){\color[rgb]{0,0,0}\makebox(0,0)[lb]{ \color{white} \smash{$\gamma$}}}%
	\end{picture}}
\subfigure[Order parameters over the inspection line]{\begin{picture}(60,50)
	\graphicspath{{figures/}}
    \put(0,3){\includegraphics[width=0.35\textwidth]{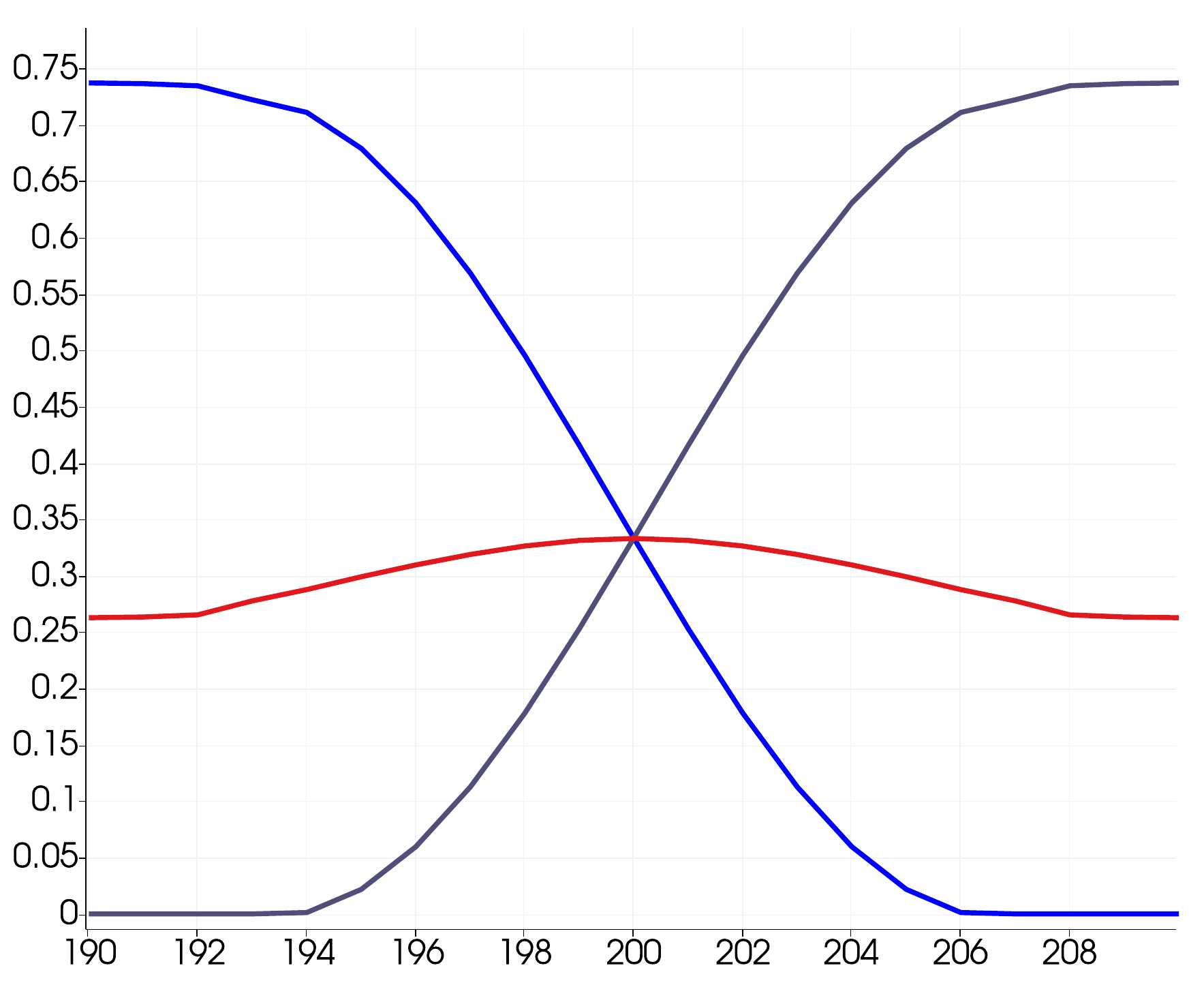}}
    \put(20.405323271,-0.455607354){\color[rgb]{0,0,0}\makebox(0,0)[lb]{\scriptsize \smash{$y \textrm{ in } 10^{-4} \textrm{ m}$}}}%
    \put(6.405323271,48.055607354){\color[rgb]{0,0,0}\makebox(0,0)[lb]{ \smash{$\phi_\alpha$}}}%
    \put(50.405323271,48.055607354){\color[rgb]{0,0,0}\makebox(0,0)[lb]{ \smash{$\phi_\beta$}}}%
    \put(12.405323271,24.055607354){\color[rgb]{0,0,0}\makebox(0,0)[lb]{ \smash{$\phi_\gamma$}}}%
	\end{picture}}	
	\caption{Triple junction example; geometry with zoom in the triple junction zone. The interfacial regions are marked by white lines. The inspection line can be seen in black color which is passing through the center of the triple junction  }
	\label{fig:geometrythirdexample}
\end{figure}

We set up a two-dimensional triple junction consisting of three phase-fields $\alpha,\beta$, and $\gamma$ as shown in \Cref{fig:geometrythirdexample}. The main purpose of this example is to investigate the elastic energy and the related driving force in the triple junction. The used parameters are shown in  Table \ref{tab:example3}.  We impose periodic boundary conditions in all directions.  Bain strains for phase-field $\alpha$ and $\beta$ are enforced whereas phase-field $\gamma$ has no Bain strain. Phase-fields $\alpha$ and $\beta$ can be interpreted as two martensitic variants while phase-field $\gamma$ is an austenitic phase. The inspection line begins at the dual interfacial region between $\alpha$ and $\gamma$ through the triple junction center ($\phi_\alpha = \phi_\beta = \phi_\gamma$) and ends at the dual interfacial region between $\beta$ and $\gamma$.

\begin{table}[H]
    \centering
    \begin{tabular}{|c|c|} 
        \hline
        \textbf{parameter} & \textbf{value} \\ 
        \hline
        system size: & $400  \times 400 \times 1$ grid points \\
        grid spacing: $\Delta x$ & $10^{-4}$ m \\   
       interface width: $\eta$ & $10$  $\Delta x$ \\ 
       Lam\'e constants & $\lambda = 80$ GPa and $\mu = 120$ GPa \\
   Bain strains &  $\Bvarepsilon_\alpha^\textrm{B} = \textrm{diag}\left(0.01,-0,01,0\right)$ \\
				 &  $\Bvarepsilon_\beta^\textrm{B} = \textrm{diag}\left(-0.01,0,01,0\right)$ \\           \hline
    \end{tabular}
    \caption{Parameters used for triple junction example.}
    \label{tab:example3}
\end{table}

\begin{figure}[h!]
	\centering
	\unitlength=1mm
\subfigure[Elastic energy contour of MPFR1 model]{\begin{picture}(60,75)
	\graphicspath{{figures/}}
    \put(0,4){\includegraphics[width=0.37\textwidth]{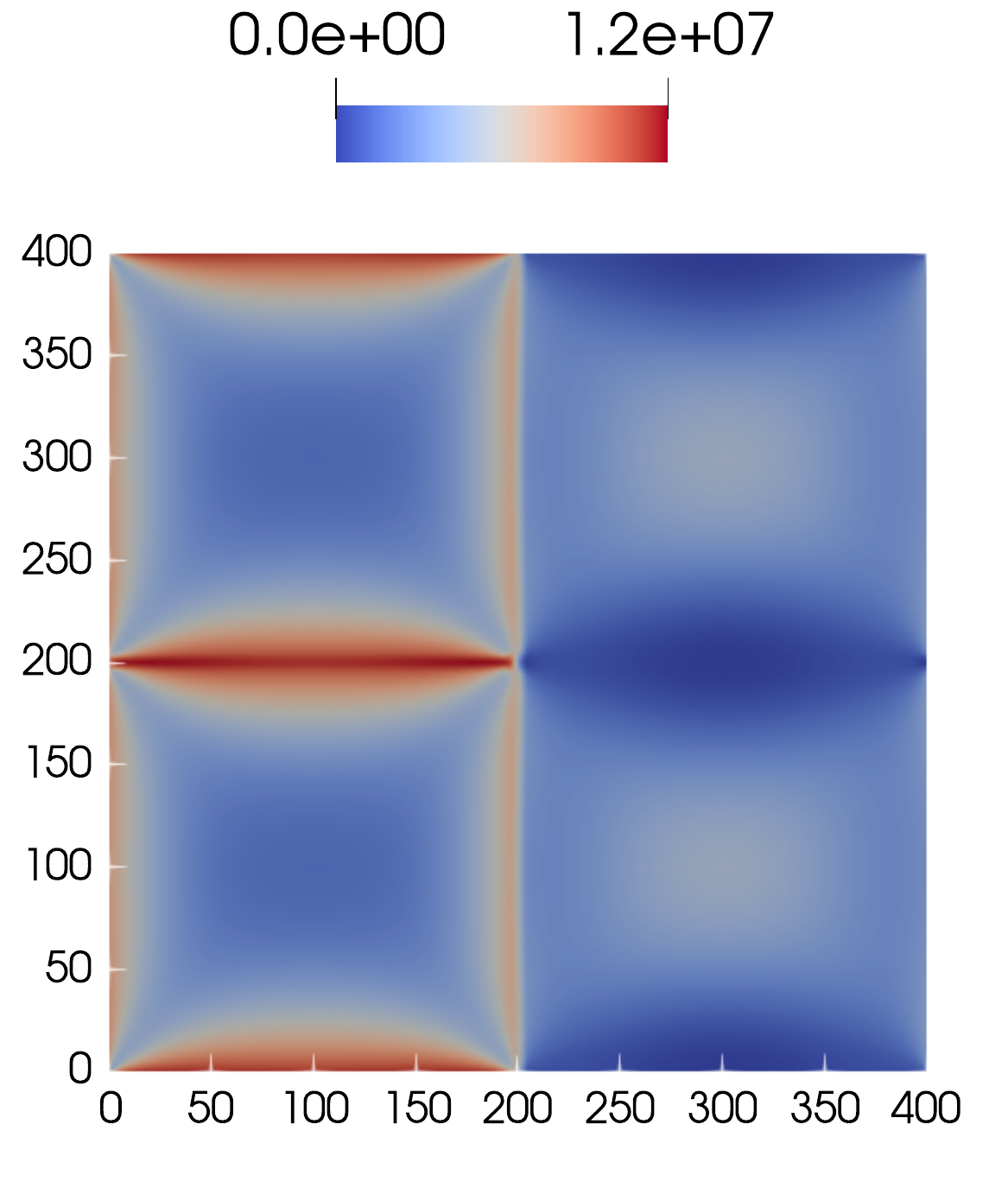}}
    \put(20.405323271,3.055607354){\color[rgb]{0,0,0}\makebox(0,0)[lb]{ \smash{$x \textrm{ in } 10^{-4} \textrm{ m}$}}}%
 \put(0.005323271,25.055607354){\color[rgb]{0,0,0}\rotatebox{90}{\makebox(0,0)[lb]{ \smash{$y \textrm{ in } 10^{-4} \textrm{ m}$}}}}%
     \put(20.405323271,62.055607354){\color[rgb]{0,0,0}\makebox(0,0)[lb]{ \smash{$\psi^\textrm{elas} \quad  \textrm{J}/\textrm{m}^3$}}}%
	\end{picture}}
\subfigure[Effective elastic energy and associated driving force of ($\alpha$,$\beta$) pair.]{\begin{picture}(100,70)
	\graphicspath{{figures/}}
    \put(0,4){\includegraphics[width=0.6 \textwidth]{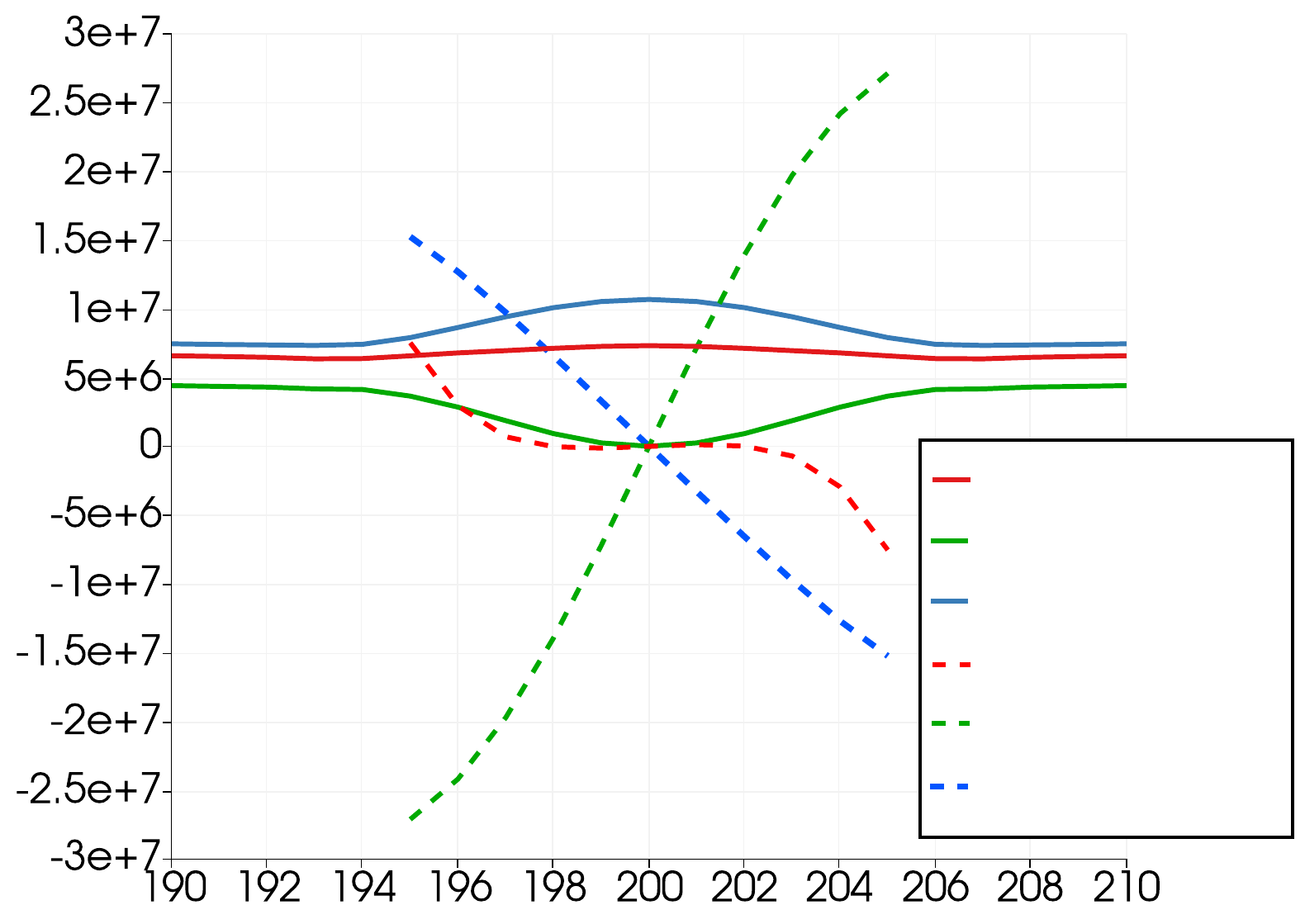}}
    \put(75.200885861,27.520908704){\color[rgb]{0,0,0}\makebox(0,0)[lb]{\smash{\scriptsize $\psi^\textrm{elas,iso-strain}$}}}%
    \put(75.200944127,32.573168089){\color[rgb]{0,0,0}\makebox(0,0)[lb]{\smash{\scriptsize $\psi^\textrm{elas,iso-stress}$}}}%
    \put(75.20063731,37.61700973){\color[rgb]{0,0,0}\makebox(0,0)[lb]{\smash{\scriptsize $\psi^\textrm{elas}$}}}%
    \put(75.200885861,23.4820908704){\color[rgb]{0,0,0}\makebox(0,0)[lb]{\smash{\scriptsize $\Delta G^\textrm{ela}_{\alpha \beta}$}}}%
    \put(75.200885861,18.4820908704){\color[rgb]{0,0,0}\makebox(0,0)[lb]{\smash{\scriptsize $\Delta G^\textrm{ela,iso-stress}_{\alpha \beta}$}}}%
    \put(75.200885861,13.4820908704){\color[rgb]{0,0,0}\makebox(0,0)[lb]{\smash{\scriptsize $\Delta G^\textrm{ela,iso-strain}_{\alpha \beta}$}}}%
    \put(40.47543573,0.00070059){\color[rgb]{0,0,0}\makebox(0,0)[lb]{\smash{$y \textrm{ in } 10^{-4} \textrm{ m}$}}}%
    \put(2.04,35.2723396){\color[rgb]{0,0,0}\rotatebox{90}{\makebox(0,0)[lb]{\smash{$\textrm{J}/\textrm{m}^3$}}}}%
     \put(17.47543573,65.00070059){\color[rgb]{0,0,0}\makebox(0,0)[lb]{\smash{$\alpha-\gamma$}}}%
     \put(40.47543573,65.00070059){\color[rgb]{0,0,0}\makebox(0,0)[lb]{\smash{$\alpha-\beta-\gamma$}}}%
     \put(70.47543573,65.00070059){\color[rgb]{0,0,0}\makebox(0,0)[lb]{\smash{$\beta-\gamma$}}}%
    
	\end{picture}}	
	\caption{Triple junction example; elastic energy and the associated driving force.   }
	\label{fig:geometrysecondexample}
\end{figure}

In \Cref{fig:geometrysecondexample}(a), we show the elastic energy contour of the MPFR1 model, the energy profile is symmetric with respect to the interface between phase-fields $\alpha$ and $\beta$. In other words, phase-fields $\alpha$ and $\beta$ are energetically equivalent so the interface between this pair should not move and the driving force should vanish.  In  \Cref{fig:geometrysecondexample}(b), we show the elastic energy and the driving force of the phase-field pair $\alpha$ and $\beta$ along the inspection line for the implemented modes. The elastic energy of MPFR1 model is between the ones of equal- strain and stress models. Under the equal-stress assumption, the elastic energy drops to zero at the center of triple junction. The behavior of the equal-stress model aligns with that of interpolation models that use an effective elastic energy of the form  $\psi^{\textrm{elas,int}} (\varepsilon)$ since  $\varepsilon = 0$ at the center of triple junction. All the implemented models deliver an elastic driving force with a vanishing average, however, the one of MPFR1 show a vanish driving force over two-thirds of the interface width.

\FloatBarrier

\subsection{Cubic-to-tetragonal transformation in a three-dimensional single grain}

 We initialize an austenitic matrix with three martensitic nuclei near the center of the simulation box. Each nucleus has a different Bain strain corresponding to a cubic-to-tetragonal variant. It is known that the mathematical solution leads to a twin laminate whose outward interface's normal is oriented at $45^\circ$ with respect to two basic vectors of the original austenitic crystal and is perpendicular to the third axis. Since the kinematic compatibility is achievable and periodic boundary conditions are enforced, the formed twin laminate of any two variants should be stress-free. The parameters used in the example are shown in  \Cref{tab:example4}. A chemical driving force, that favors the growth of the martensitic nuclei, is added. Furthermore, we enforce periodic boundary conditions  which should allow building a stress-free laminate.

\begin{table}[h!]
    \centering
    \begin{tabular}{|c|c|} 
        \hline
        \textbf{parameter} & \textbf{value} \\ 
        \hline
        system size: & $100  \times 100 \times 100$ grid points \\
        grid spacing: $\Delta x$ & $10^{-7}$ m \\ 
       interface width: $\eta$ & $5$  $\Delta x$ \\ 
       time step size: $\Delta t$  &  $5\,\cdot\,10^{-9}$ sec \\
      interfacial energy:  $\gamma$  & $0.1 \, \textrm{J}/\textrm{m}^2$ \\
       interfacial mobility: $M$  &  $3 \cdot 10^{-7} \, \textrm{m}^4 \,  / ( \textrm{J}\cdot \textrm{sec})$  \\
      Lam\'e constants & $\lambda = 80$ GPa and $\mu = 120$ GPa \\
   Bain strains &  $ \Bvarepsilon^\textrm{B}_{\textrm{M},1} = \textrm{diag}\left(0.02,-0.01,-0.01\right)$ \\
 &  $\Bvarepsilon^\textrm{B}_{\textrm{M},2} = \textrm{diag}\left(-0.01,-0.01,0.02\right)$  \\
 &  $\Bvarepsilon^\textrm{B}_{\textrm{M},3} = \textrm{diag}\left(-0.01,0.02,-0.01\right)$  \\
chemical driving force: $\Delta G_{A M}^\textrm{chem}$ & $7.5 \cdot 10^7 \, \textrm{J}/\textrm{m}^3$ \\ 
        \hline
    \end{tabular}
    \caption{Parameters of the martensitic transformation in a single grain example.}
    \label{tab:example4}
\end{table}

\begin{figure}[ht] 
 	\centering
	\unitlength=1mm
\subfigure[MPFR1 (t = 0)]{\begin{picture}(40,35)
	\graphicspath{{figures/}}
    \put(0.0,0.0){\includegraphics[width=0.24 \textwidth]{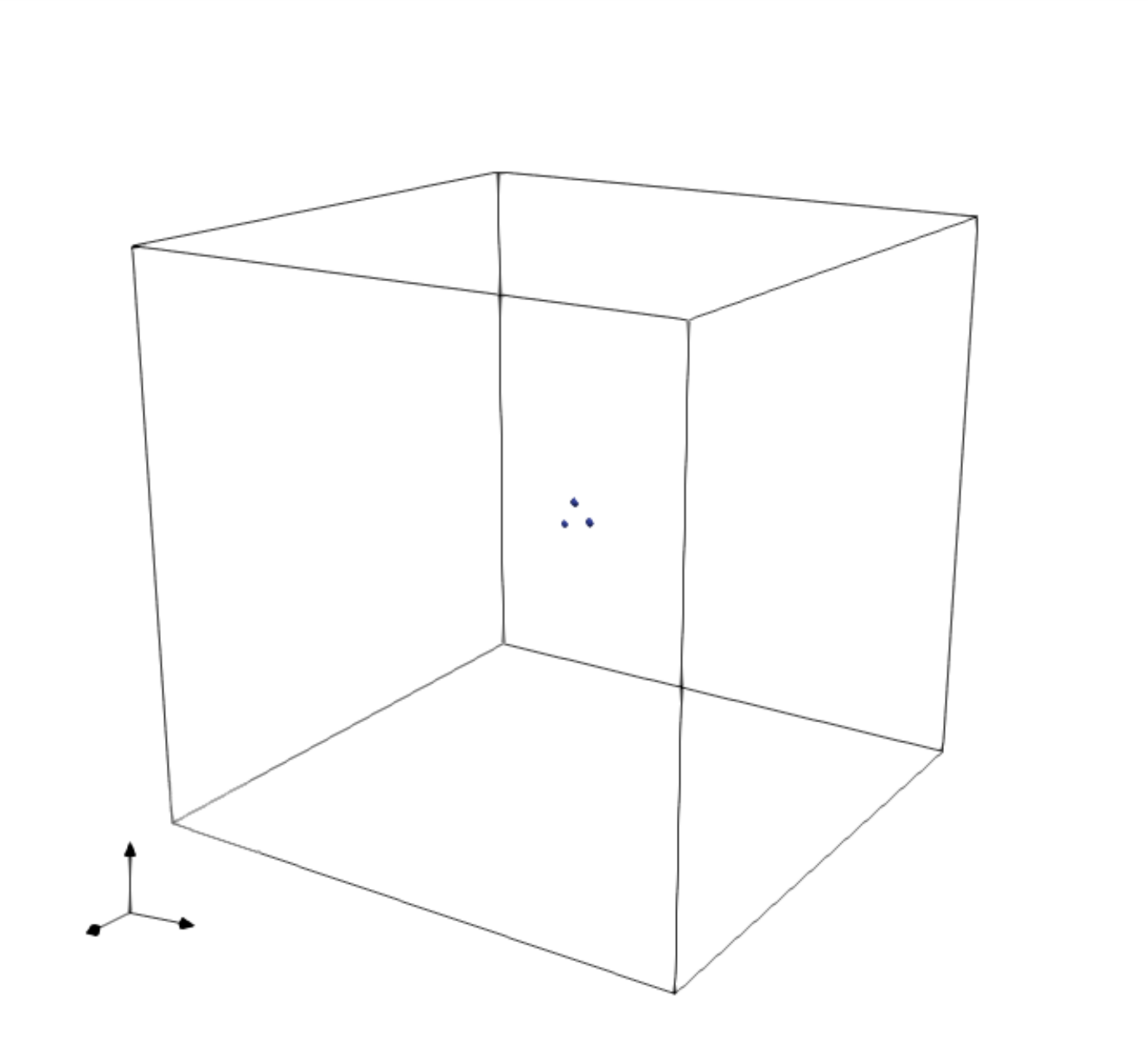}}
	\end{picture}} 
\subfigure[MPFR1 ($t = 200 \Delta t$)]{\begin{picture}(40,35)
	\graphicspath{{figures/}}
    \put(0.0,0.0){\includegraphics[width=0.24 \textwidth]{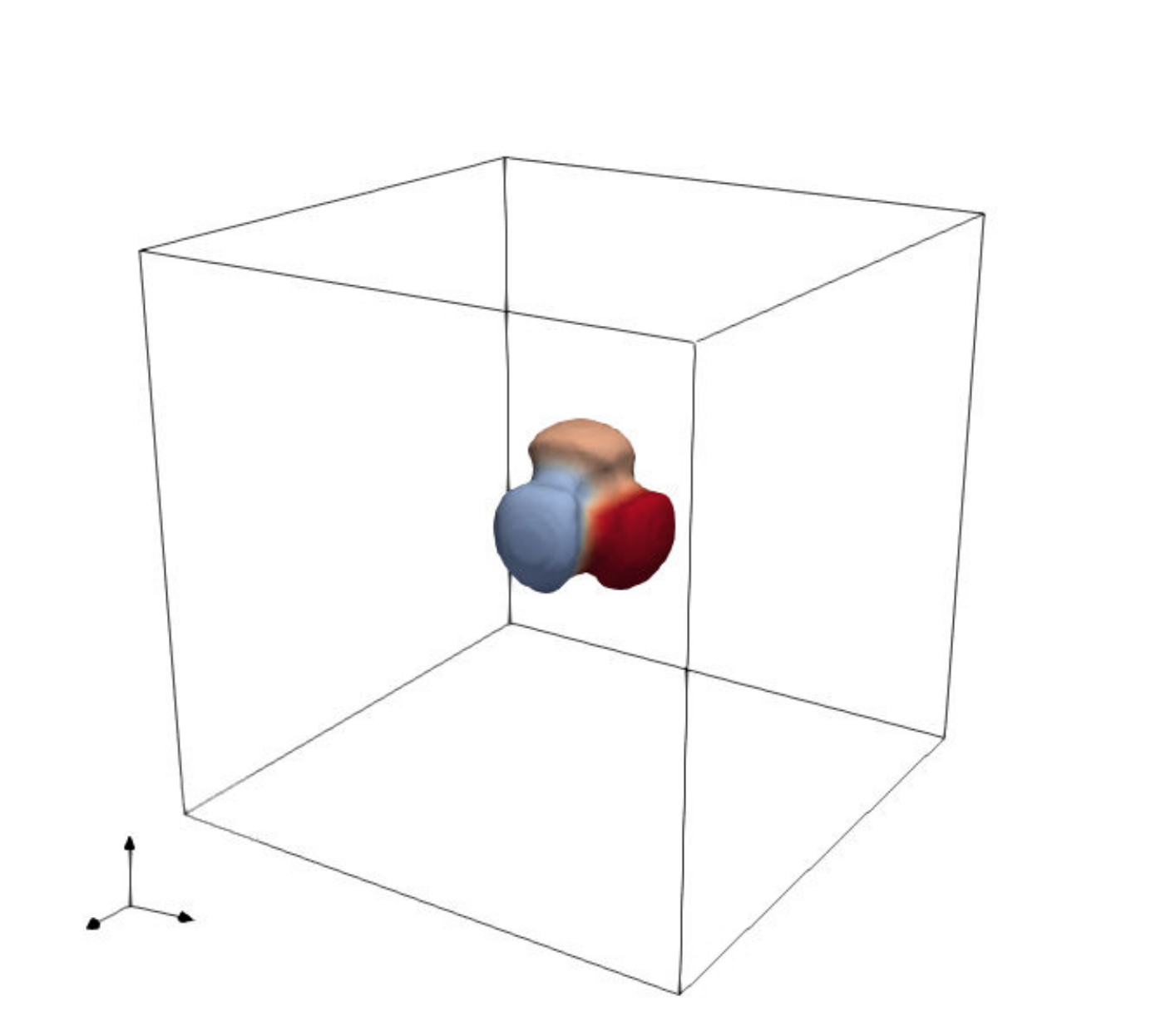}}
	\end{picture}} 
\subfigure[MPFR1 ($t = 800 \Delta t$)]{\begin{picture}(40,35)
	\graphicspath{{figures/}}
    \put(0.0,0.0){\includegraphics[width=0.24 \textwidth]{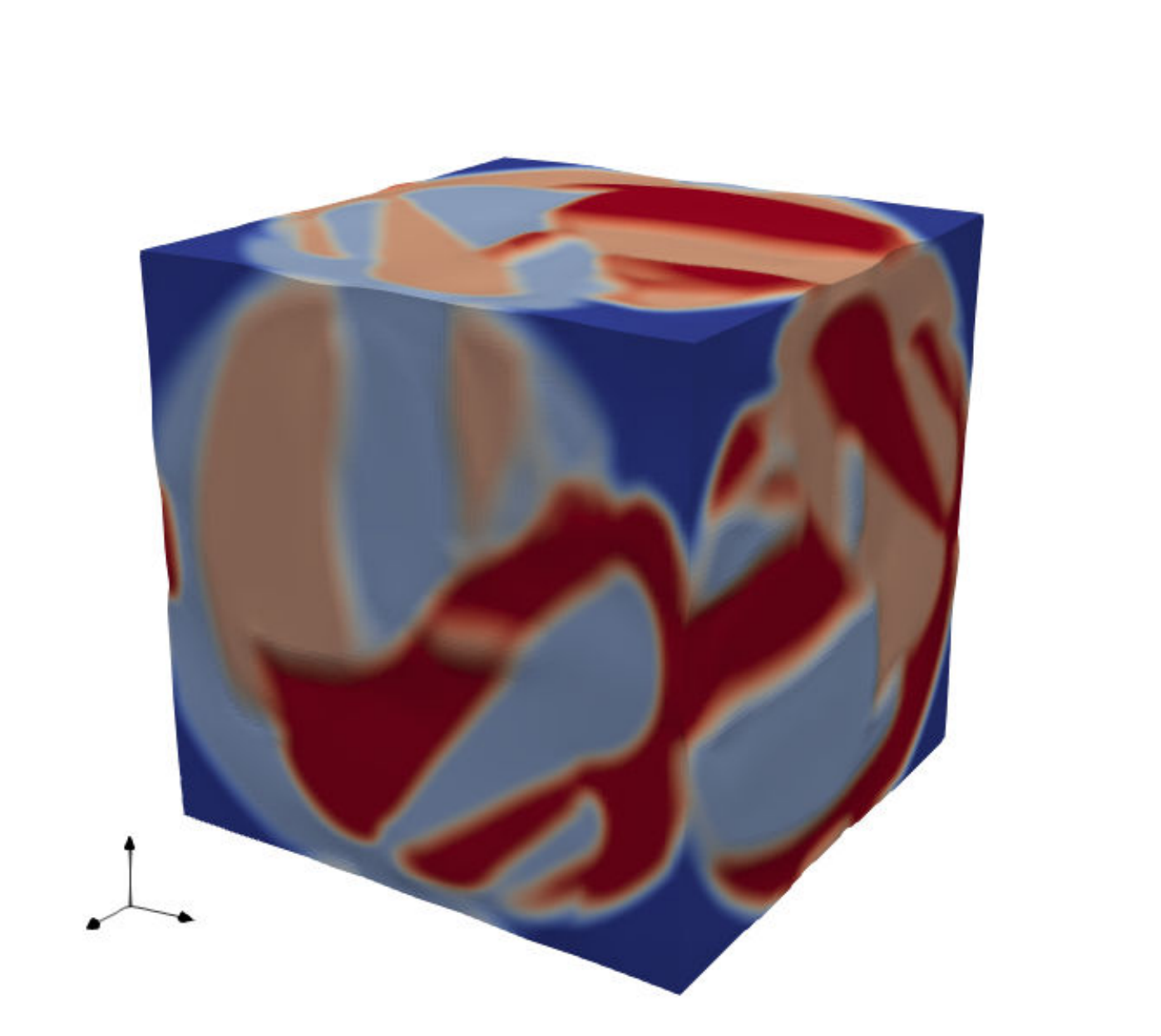}}
	\end{picture}} 
\subfigure[MPFR1 (final)]{\begin{picture}(40,35)
	\graphicspath{{figures/}}
    \put(0.0,0.0){\includegraphics[width=0.24 \textwidth]{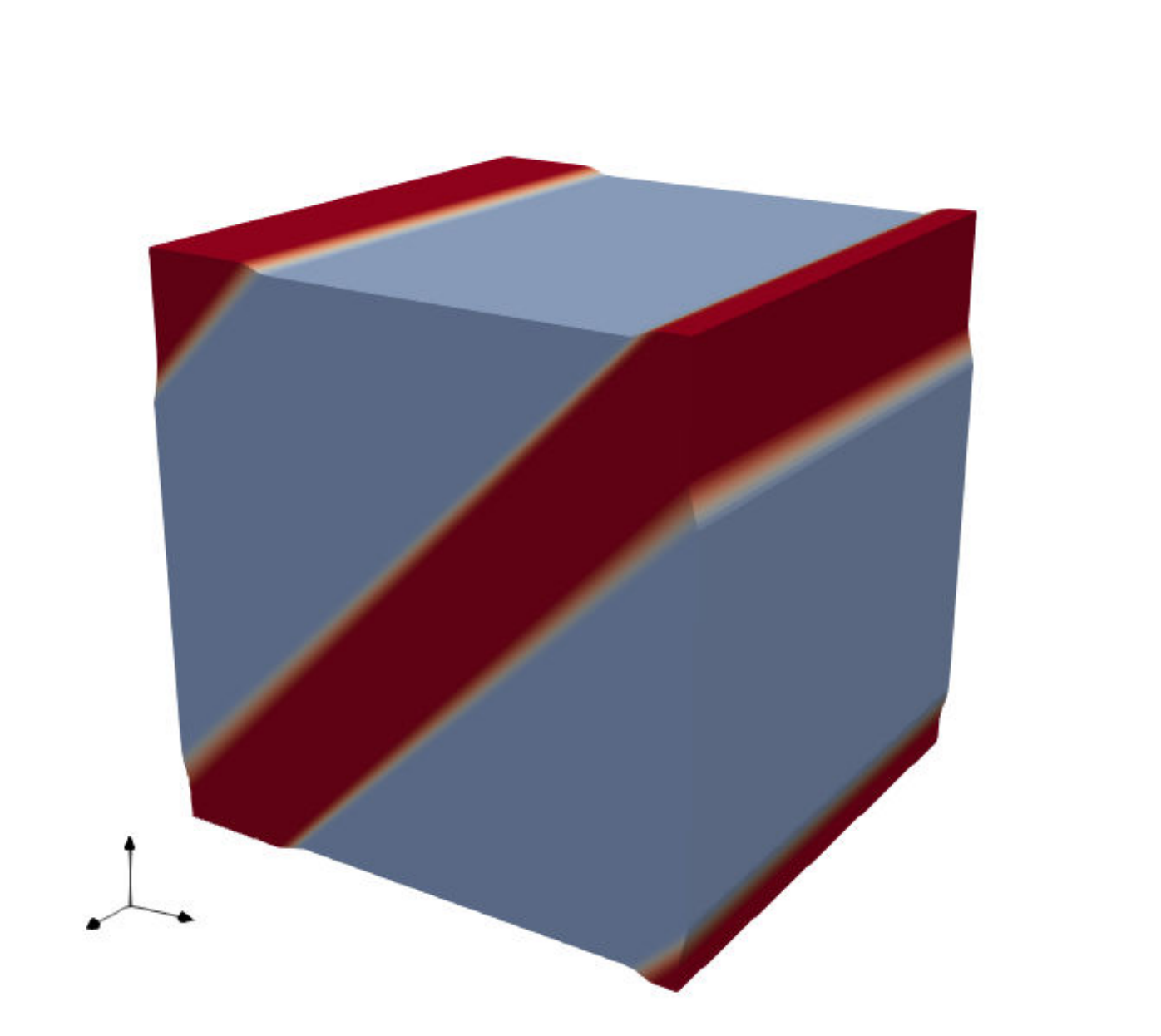}}
	\end{picture}} 
\subfigure[equal-stress (final)]{\begin{picture}(40,35)
	\graphicspath{{figures/}}
    \put(0.0,0.0){\includegraphics[width=0.24 \textwidth]{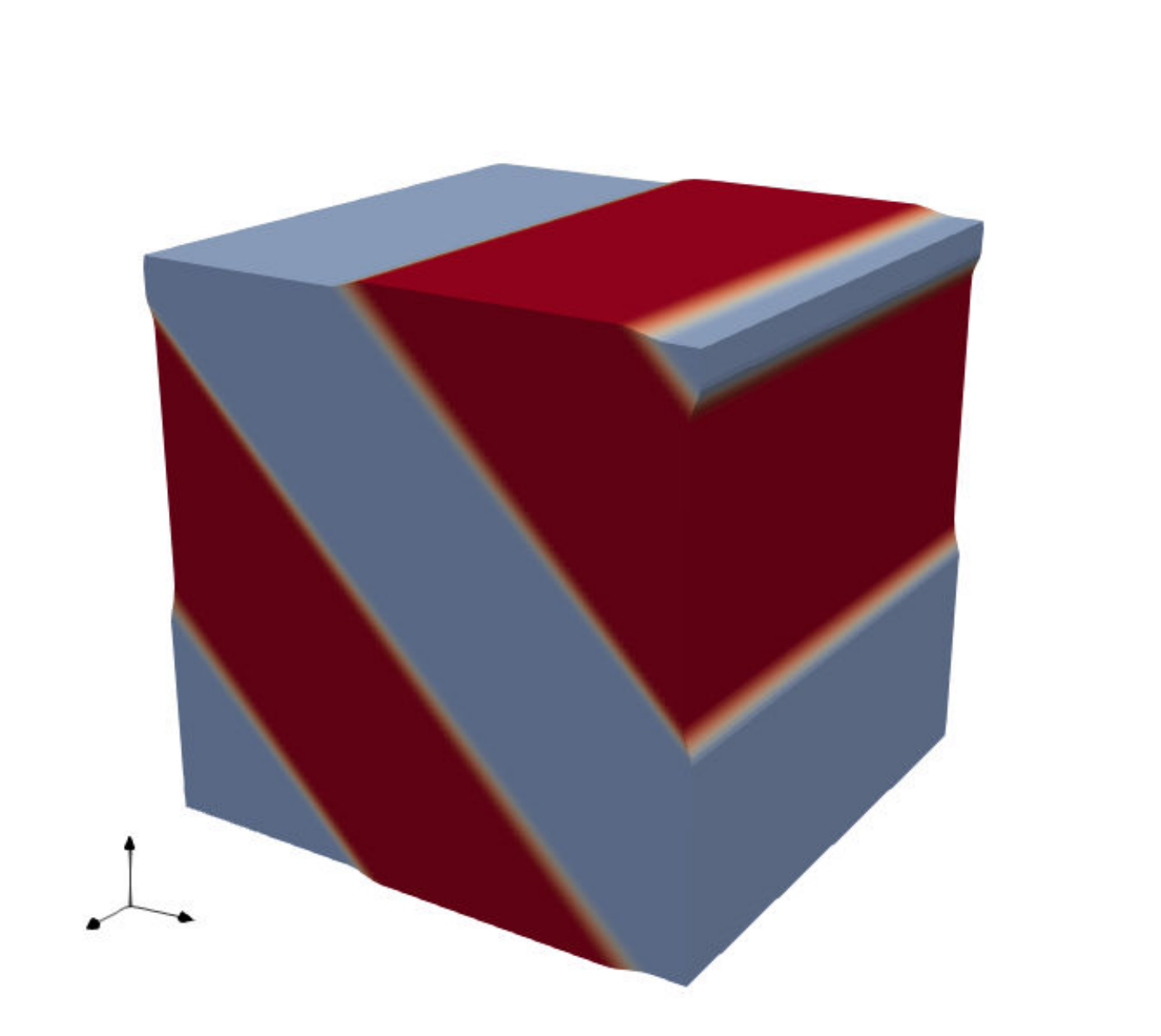}}
	\end{picture}} 
\subfigure[equal-strain (final)]{\begin{picture}(40,35)
	\graphicspath{{figures/}}
    \put(0.0,0.0){\includegraphics[width=0.24 \textwidth]{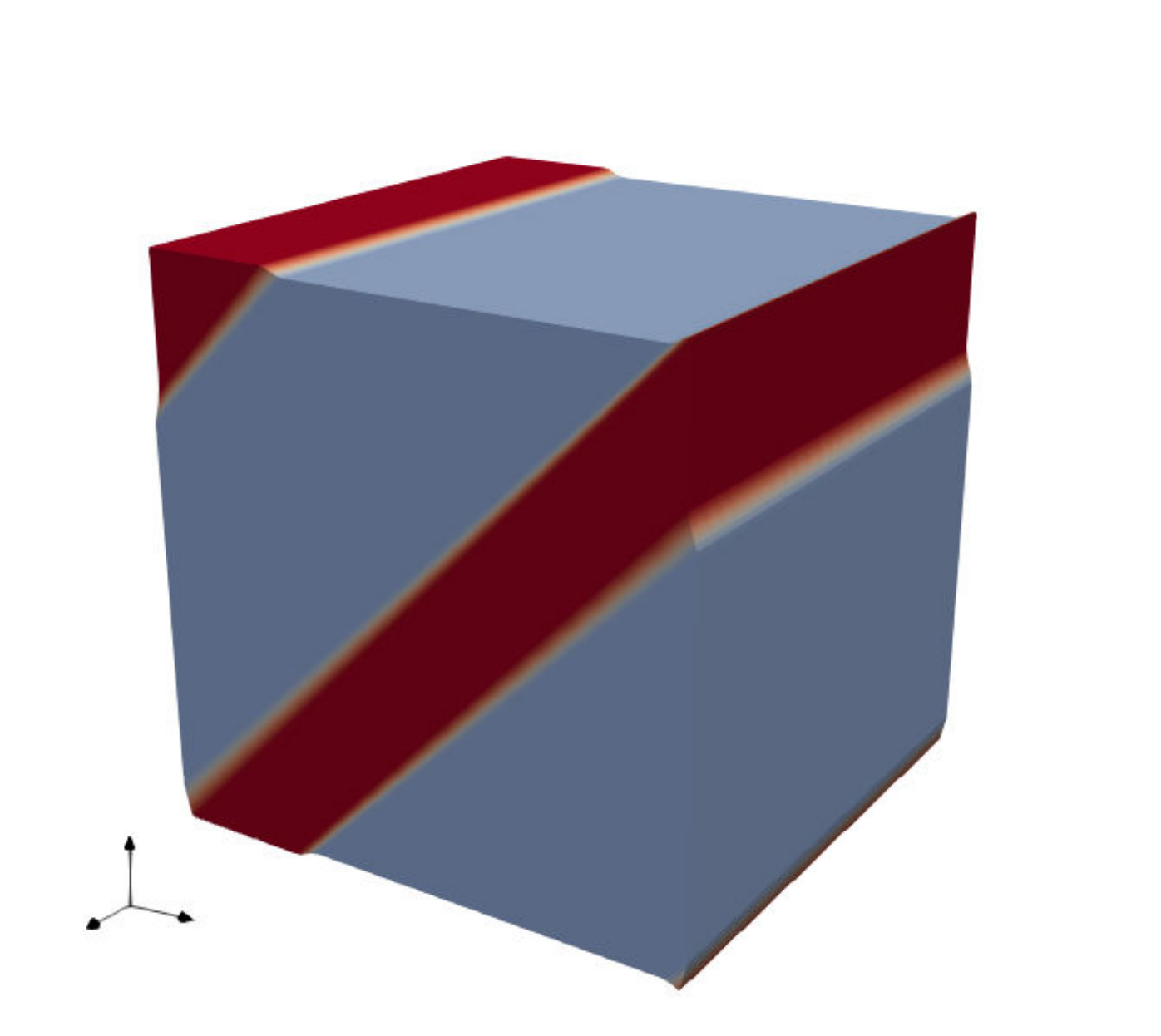}}
	\end{picture}} 
\subfigure[mechanical energy of equal-strain]{\begin{picture}(60,40)
	\graphicspath{{figures/}}
    \put(6.5,0.0){\includegraphics[width=0.24 \textwidth]{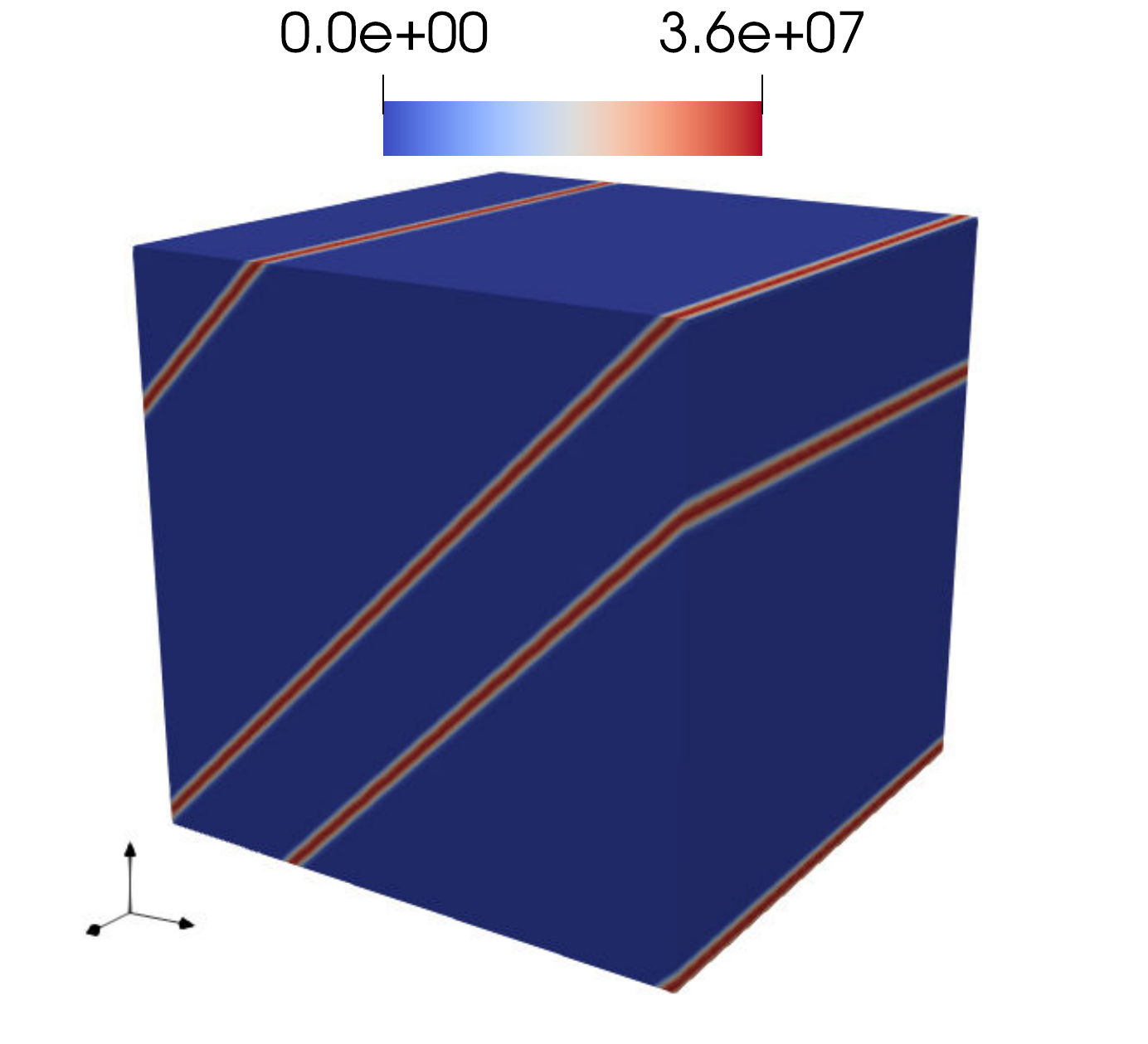}}
     \put(18.005323271,36.055607354){\color[rgb]{0,0,0}\rotatebox{0}{\makebox(0,0)[lb]{\tiny $\psi^\textrm{elas} \, (\textrm{J}/\textrm{m}^3)$}}} 
	\end{picture}} 
	\caption{Results of the martensitic transformation in a single grain. The microstructure evolution using the developed model is shown in (a-d) while the final microstructures using Voigt/Taylor and Reuss/Sachs approaches appear in (e-f). The mechanical energy of the final microstructure obtained by equal-strain model is shown in (g).  }
	\label{fig:resultsfourthexample}
\end{figure}

 We show in \Cref{fig:resultsfourthexample} the microstructure evolution for the MPFR1 model along with the final microstructures obtained by the different models. All models deliver final microstructures  that match the analytical solution of a twinned laminate. The twinned laminate consists of variants 1 and 3, with an outward normal vector in the direction $[110]$. Different variants for the twinned laminate may be selected based on the initialization of the nuclei, but the characteristics of the final microstructure remain the same. We notice that the equal-strain model is unable to produce a stress-free twinned laminate as enforcing the equal-strain assumption leads to stress concentration in the interfacial regions, see \Cref{fig:resultsfourthexample}(g). This is because of not allowing any stress relaxation along the normal of the twinned laminate.

\begin{figure}[h!]
	\centering
	\unitlength=1mm
\subfigure[MPFR1]{\begin{picture}(100,50)
	\graphicspath{{figures/}}
    \put(10.0,5.0){\includegraphics[width=0.67 \textwidth]{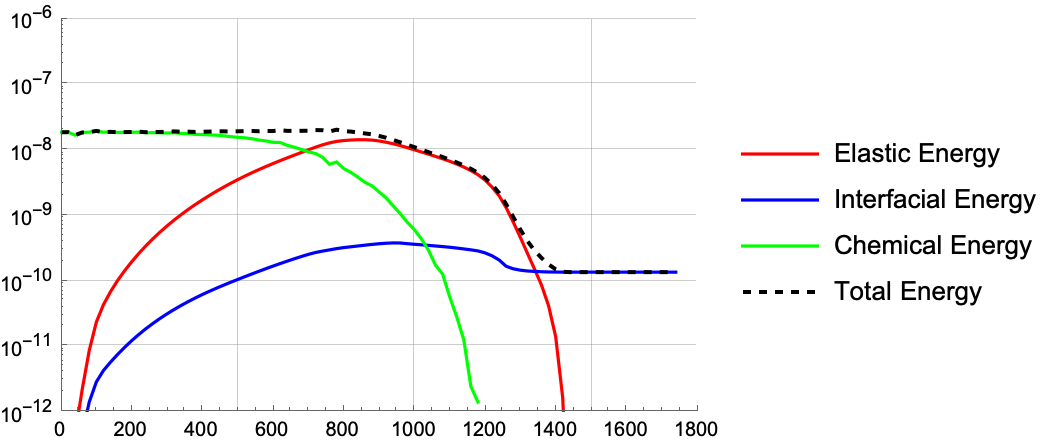}}
     \put(9.005323271,25.055607354){\color[rgb]{0,0,0}\rotatebox{90}{\makebox(0,0)[lb]{ $\textrm{Joule}$}}}
     \put(40.005323271,0.055607354){\color[rgb]{0,0,0}\rotatebox{0}{\makebox(0,0)[lb]{ time step}}} 
	\end{picture}} 
\subfigure[equal-stress]{\begin{picture}(80,50)
	\graphicspath{{figures/}}
    \put(0.0,5.0){\includegraphics[width=0.45 \textwidth]{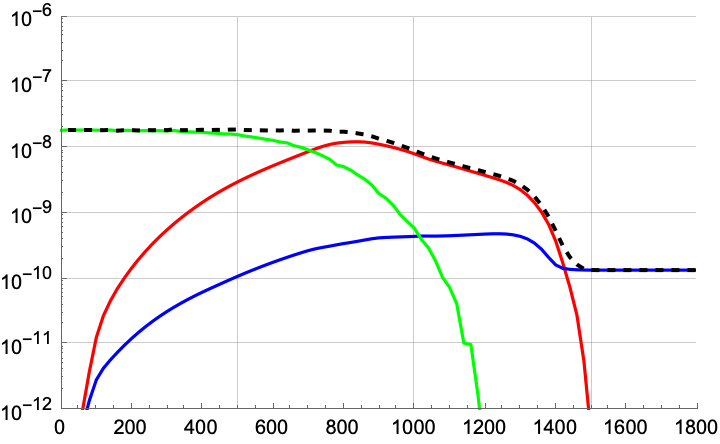}}
         \put(0.005323271,25.055607354){\color[rgb]{0,0,0}\rotatebox{90}{\makebox(0,0)[lb]{ $\textrm{Joule}$}}}
     \put(30.005323271,0.055607354){\color[rgb]{0,0,0}\rotatebox{0}{\makebox(0,0)[lb]{ time step}}} 
	\end{picture}} 
\subfigure[equal-strain]{\begin{picture}(80,50)
	\graphicspath{{figures/}}
    \put(0.0,5.0){\includegraphics[width=0.45 \textwidth]{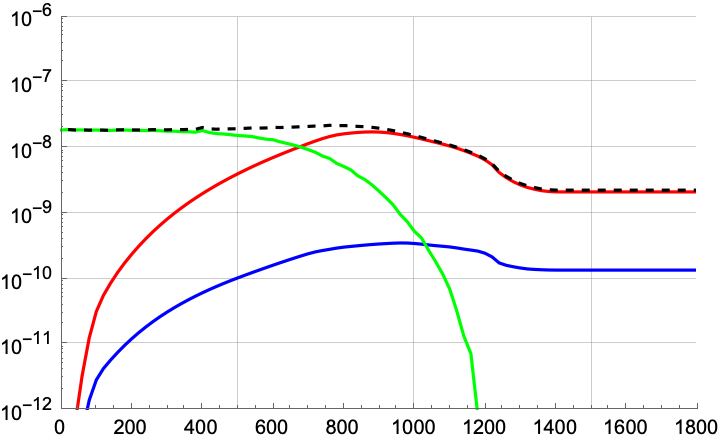}}
         \put(0.005323271,25.055607354){\color[rgb]{0,0,0}\rotatebox{90}{\makebox(0,0)[lb]{ $\textrm{Joule}$}}}
     \put(30.005323271,0.055607354){\color[rgb]{0,0,0}\rotatebox{0}{\makebox(0,0)[lb]{ time step}}} 
	\end{picture}} 
	\caption{Energy evolution over time for cubic-to-tetragonal transformation in a single grain}
	\label{fig:resultsfourthexample_2}
\end{figure}

\begin{figure}[h!]
 \center
 \begin{picture}(160,160)
	\graphicspath{{figures/}}
    \put(-40.0,0.0){\includegraphics[width=0.65 \textwidth]{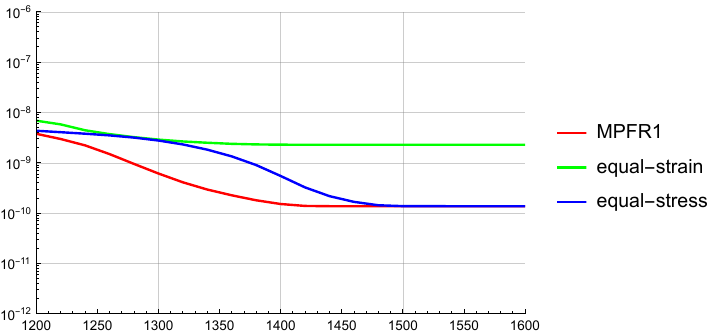}}
         \put(-45.005323271,50.055607354){\color[rgb]{0,0,0}\rotatebox{90}{\makebox(0,0)[lb]{ $\textrm{Joule}$}}}
     \put(70.005323271,-10.055607354){\color[rgb]{0,0,0}\rotatebox{0}{\makebox(0,0)[lb]{ time step}}} 
	\end{picture} 
	\caption{Total energy evolution for cubic-to-tetragonal transformation in a single grain. Equal-strain model reaches the final microstructure first, MPFR1 model is second and equal-stress is last. The final total energy for MPFR1 and equal-stress models include only interfacial energy whereas equal-strain model contains elastic energy as well.  }
		\label{fig:resultsfourthexample_3}
\end{figure}

We investigate in \Cref{fig:resultsfourthexample_2} the evolution of the total energy with the individual energy contributions for the implemented models. The chemical energy becomes zero when austenite phase vanishes from the simulation domain. The elastic energy becomes the dominating driving force, which vanishes when a twinned laminate is formed except for equal-strain model whose elastic energy remains non-zero. The final interfacial energy of all models is nearly the same. We observe that the transformation speeds differ: the equal-strain model achieves the final microstructure the fastest, the equal-stress model is the slowest, and MPFR1 lies in between, see \Cref{fig:resultsfourthexample_3}. 
 
\FloatBarrier

\subsection{Cubic-to-tetragonal transformation in a two-dimensional polycrystal} 

\begin{figure}[h!]
 \center
 \begin{picture}(200,200)
	\graphicspath{{figures/}}
    \put(0.0,0.0){\includegraphics[width=0.4 \textwidth]{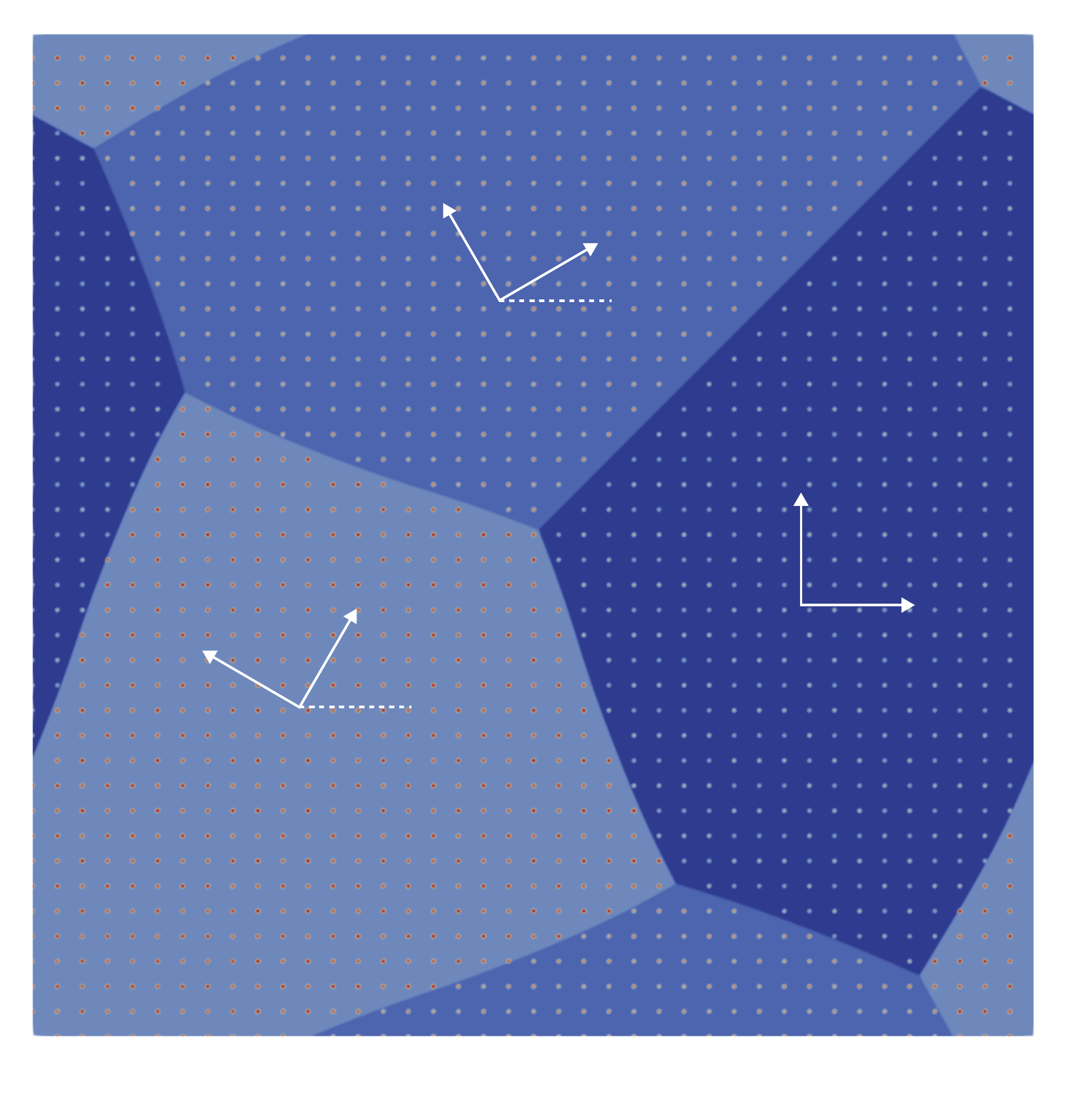}}
      \put(58.005323271,70.055607354){\color[rgb]{1,1,1}\rotatebox{0}{\makebox(0,0)[lb]{\tiny $60^\textrm{o}$}}} 
      \put(96.005323271,139.055607354){\color[rgb]{1,1,1}\rotatebox{0}{\makebox(0,0)[lb]{\tiny $30^\textrm{o}$}}}
	\end{picture} 
	\caption{ The two-dimensional polycrystal at $t=0$ consisting of three grains with different orientation of the parent cubic crystals. The nuclei of the martensite tetragonal varanits are the white dots.   }
		\label{fig:geometryfourthexample_5}
\end{figure}

We consider here a two-dimensional polycrystal with three grains.  Different orientations are assumed in each grain for the parent austenite phase, see \Cref{fig:geometryfourthexample_5}.   We plant martensitic nuclei at equal distances in the grains. A chemical driving force is added to enforce the growth of the nuclei.  The parameters used in this simulation are shown in    \Cref{tab:example5}.

\begin{table}[h!]
    \centering
    \begin{tabular}{|c|c|} 
        \hline
        \textbf{parameter} & \textbf{value} \\ 
        \hline
        system size: & $1000  \times 1000 $ grid points \\
        grid spacing: $\Delta x$ & $10^{-6}$ m \\ 
       interface width: $\eta$ & $3$  $\Delta x$ \\ 
       time step size: $\Delta t$  &  $10^{-6}$ sec \\
      interfacial energy:  $\gamma$  & $0.1 \, \textrm{J}/\textrm{m}^2$ \\
       interfacial mobility: $M$  &  $3 \cdot 10^{-7} \, \textrm{m}^4 \,  / ( \textrm{J}\cdot \textrm{sec})$  \\
      Lam\'e constants & $\lambda = 80$ GPa and $\mu = 120$ GPa \\
   Bain strains &  $ \Bvarepsilon^\textrm{B}_{\textrm{M},1} = \textrm{diag}\left(0.02,-0.01\right)$ \\
 &  $\Bvarepsilon^\textrm{B}_{\textrm{M},2} = \textrm{diag}\left(-0.01,0.02\right)$  \\
chemical driving force: $\Delta G_{A M}^\textrm{chem}$ & $7.76\cdot 10^7 \, \textrm{J}/\textrm{m}^3$ \\ 
        \hline
    \end{tabular}
    \caption{Parameters of the martensitic transformation in a two-dimensional polycrystal example.}
    \label{tab:example5}
\end{table}

\begin{figure}[h!]
	\centering
	\unitlength=1mm
\subfigure[MPFR1, $t = 400 \,\Delta t$]{\begin{picture}(50,55)
	\graphicspath{{figures/}}
	\put(0,0){\includegraphics[width=0.32\textwidth]{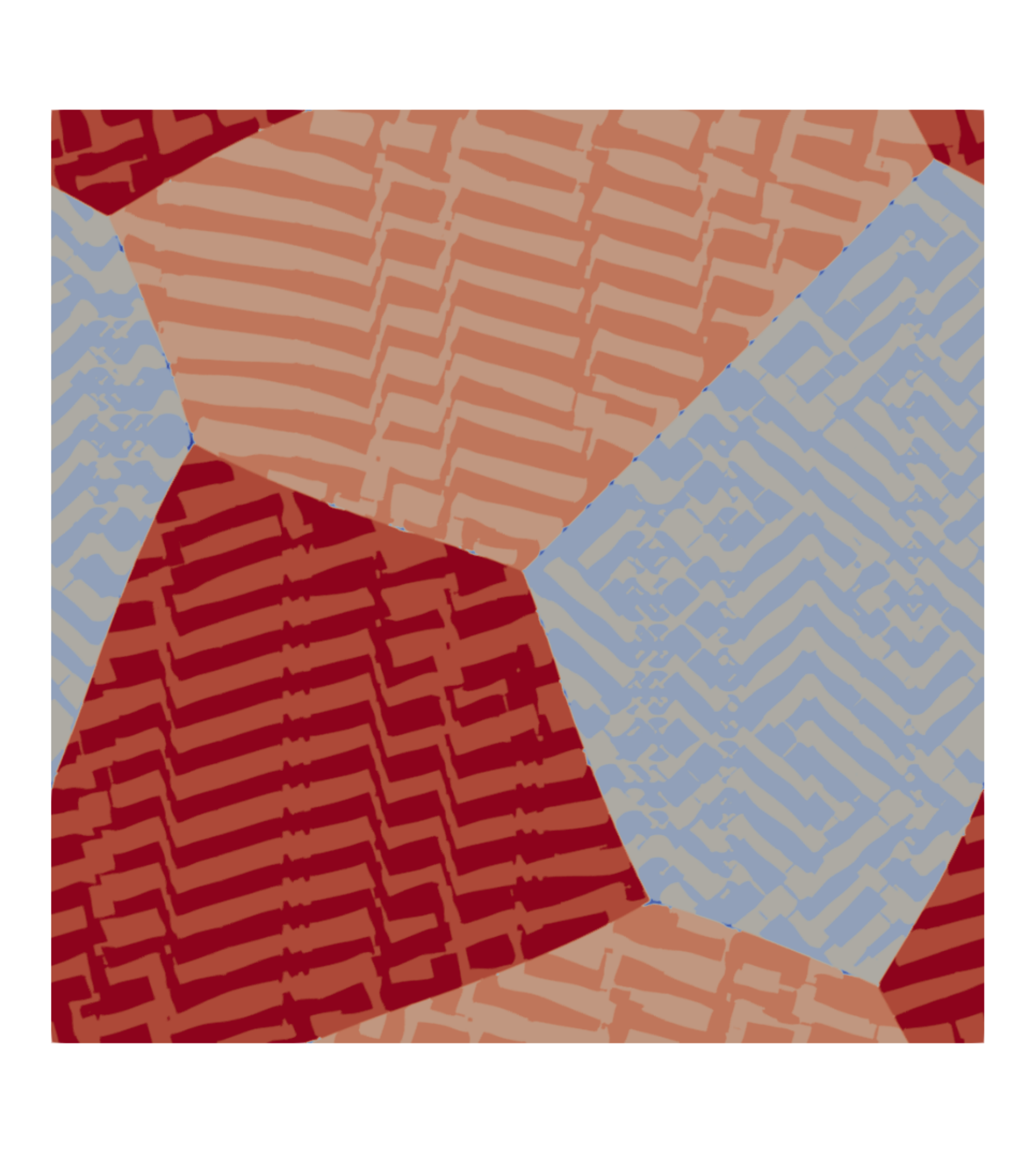}}
	\end{picture}} 
\subfigure[MPFR1, $t = 800 \,\Delta t$]{\begin{picture}(50,55)
	\graphicspath{{figures/}}
	\put(0,0){\includegraphics[width=0.32\textwidth]{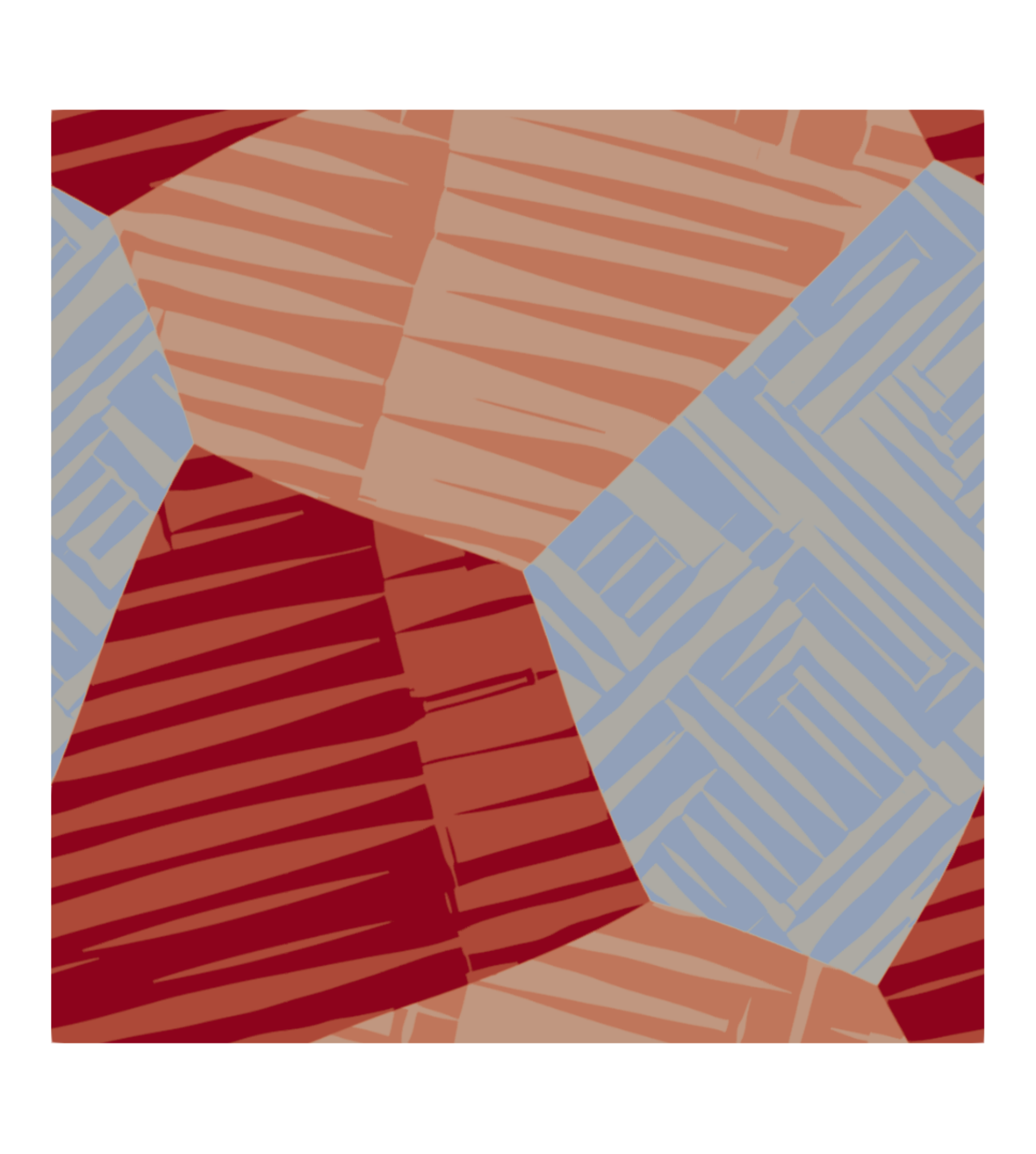}}
	\end{picture}} 
\subfigure[MPFR1, $t = 10000 \,\Delta t$]{\begin{picture}(50,5)
	\graphicspath{{figures/}}
	\put(0,0){\includegraphics[width=0.32\textwidth]{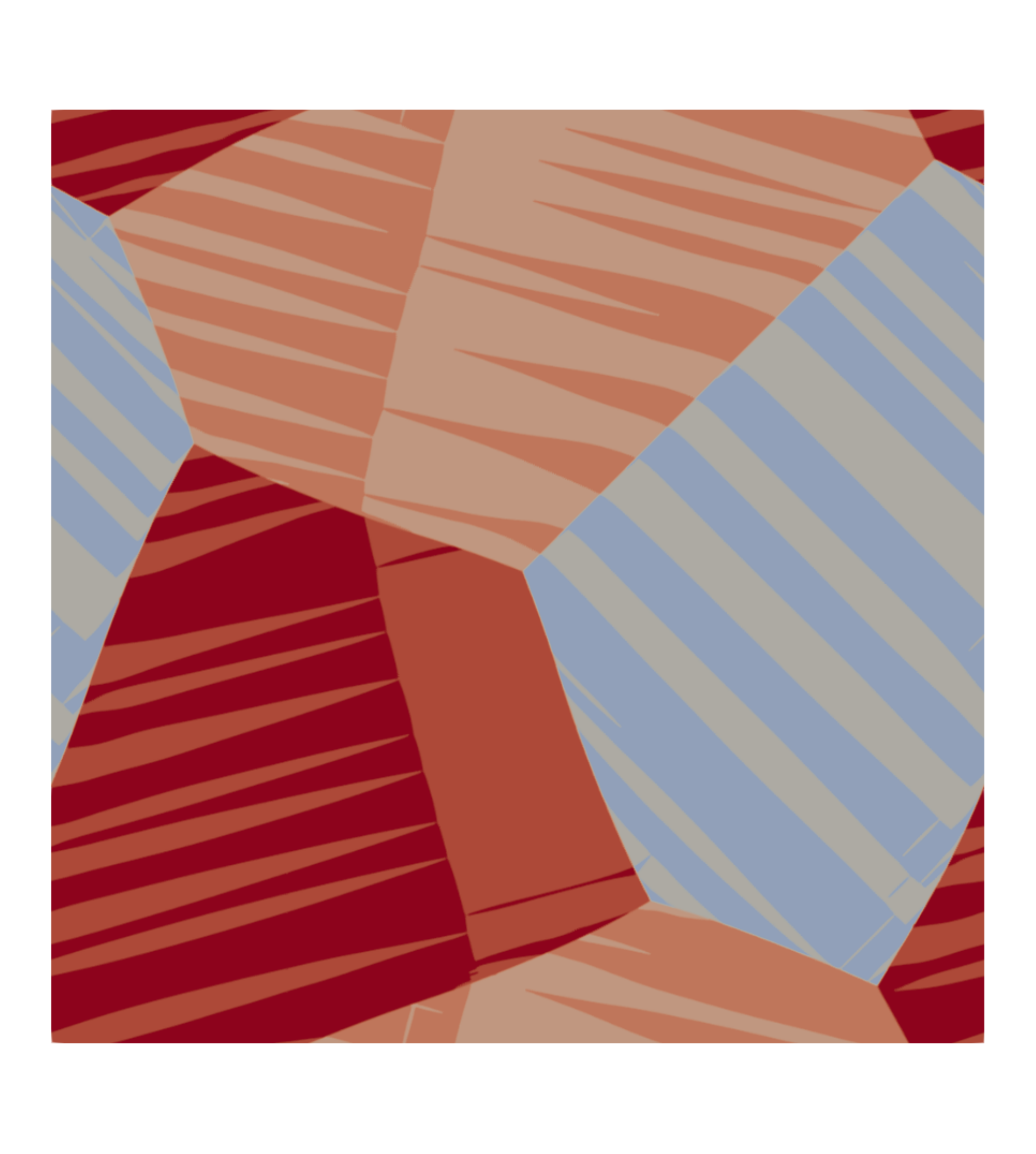}}
	\end{picture}} 
\subfigure[equal-stress, $t = 400 \,\Delta t$]{\begin{picture}(50,55)
	\graphicspath{{figures/}}
	\put(0,0){\includegraphics[width=0.32\textwidth]{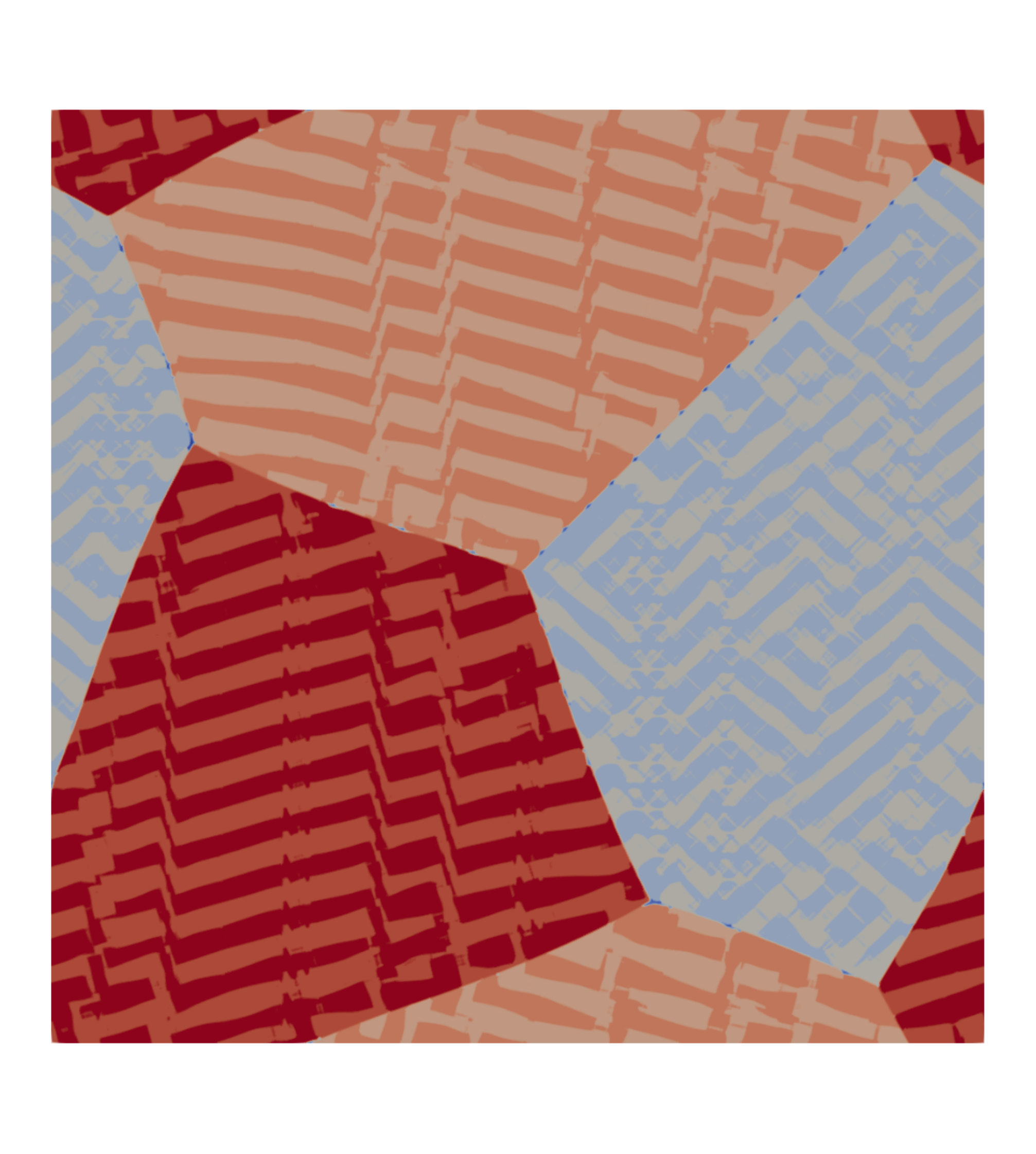}}
	\end{picture}} 
\subfigure[equal-stress, $t = 800 \,\Delta t$]{\begin{picture}(50,55)
	\graphicspath{{figures/}}
	\put(0,0){\includegraphics[width=0.32\textwidth]{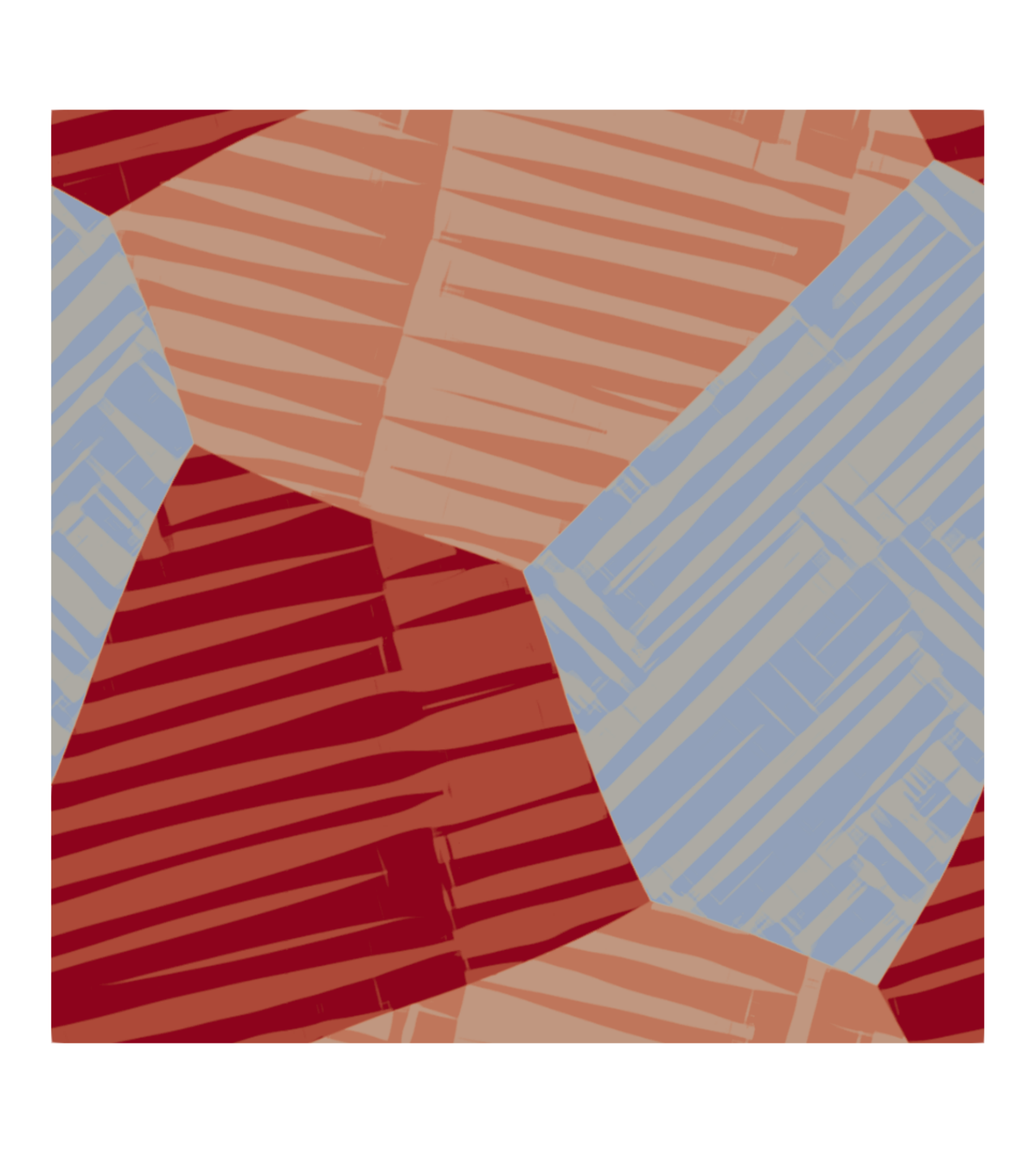}}
	\end{picture}} 
\subfigure[equal-stress, $t = 10000 \,\Delta t$]{\begin{picture}(50,55)
	\graphicspath{{figures/}}
	\put(0,0){\includegraphics[width=0.32\textwidth]{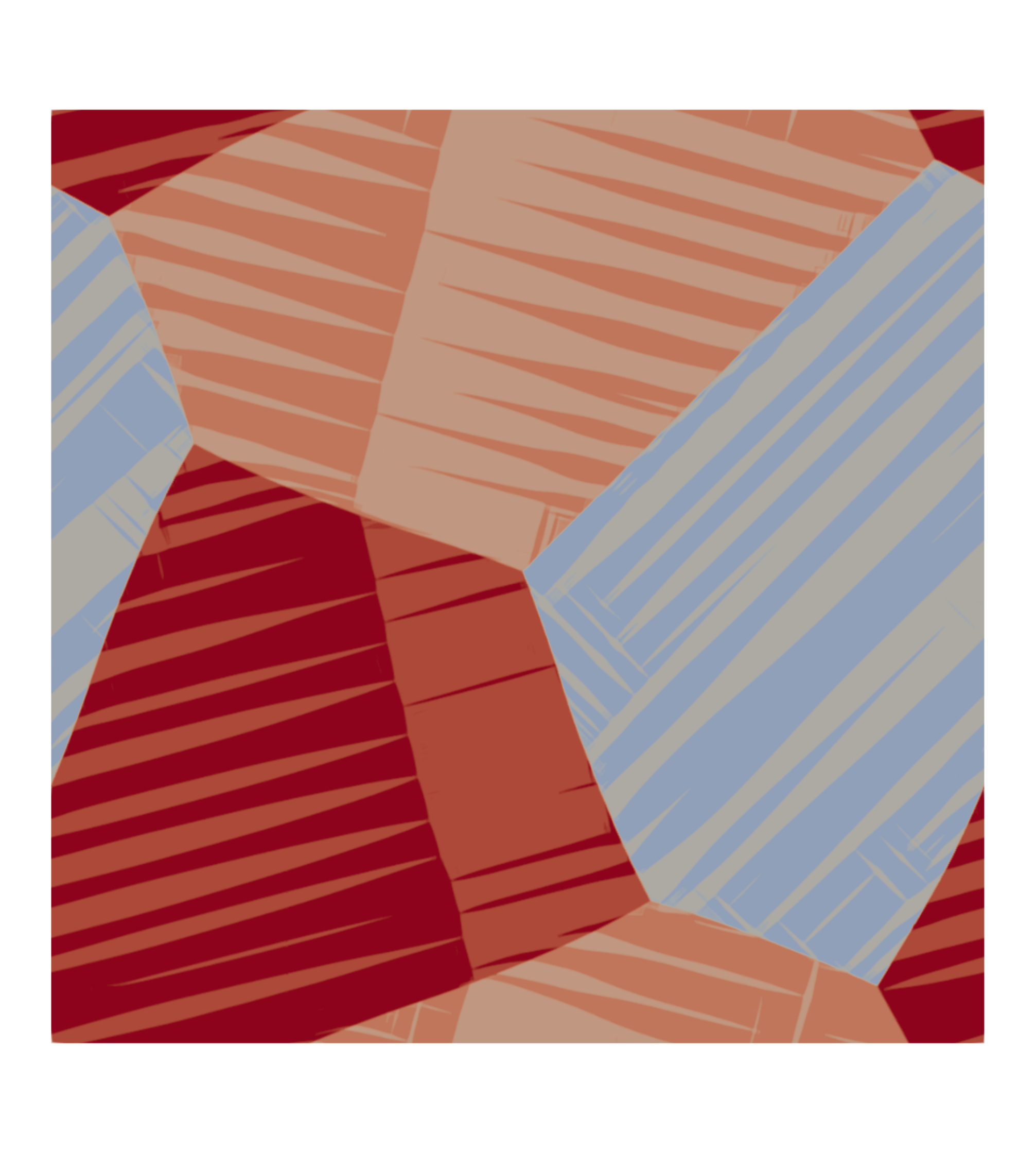}}
	\end{picture}} 
\subfigure[equal-strain, $t = 400 \,\Delta t$]{\begin{picture}(50,55)
	\graphicspath{{figures/}}
	\put(0,0){\includegraphics[width=0.32\textwidth]{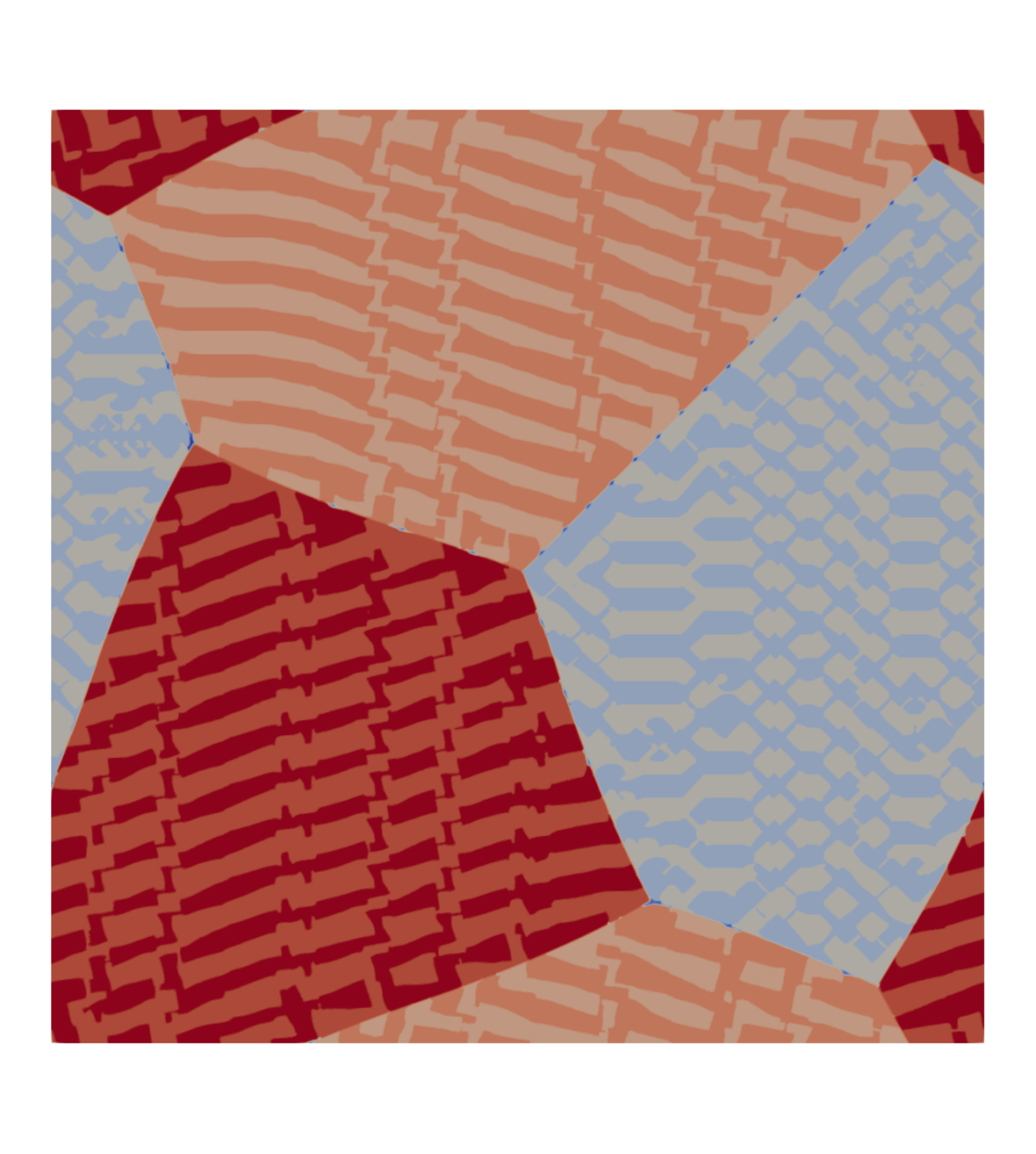}}
	\end{picture}} 
\subfigure[equal-strain, $t = 800 \,\Delta t$]{\begin{picture}(50,55)
	\graphicspath{{figures/}}
	\put(0,0){\includegraphics[width=0.32\textwidth]{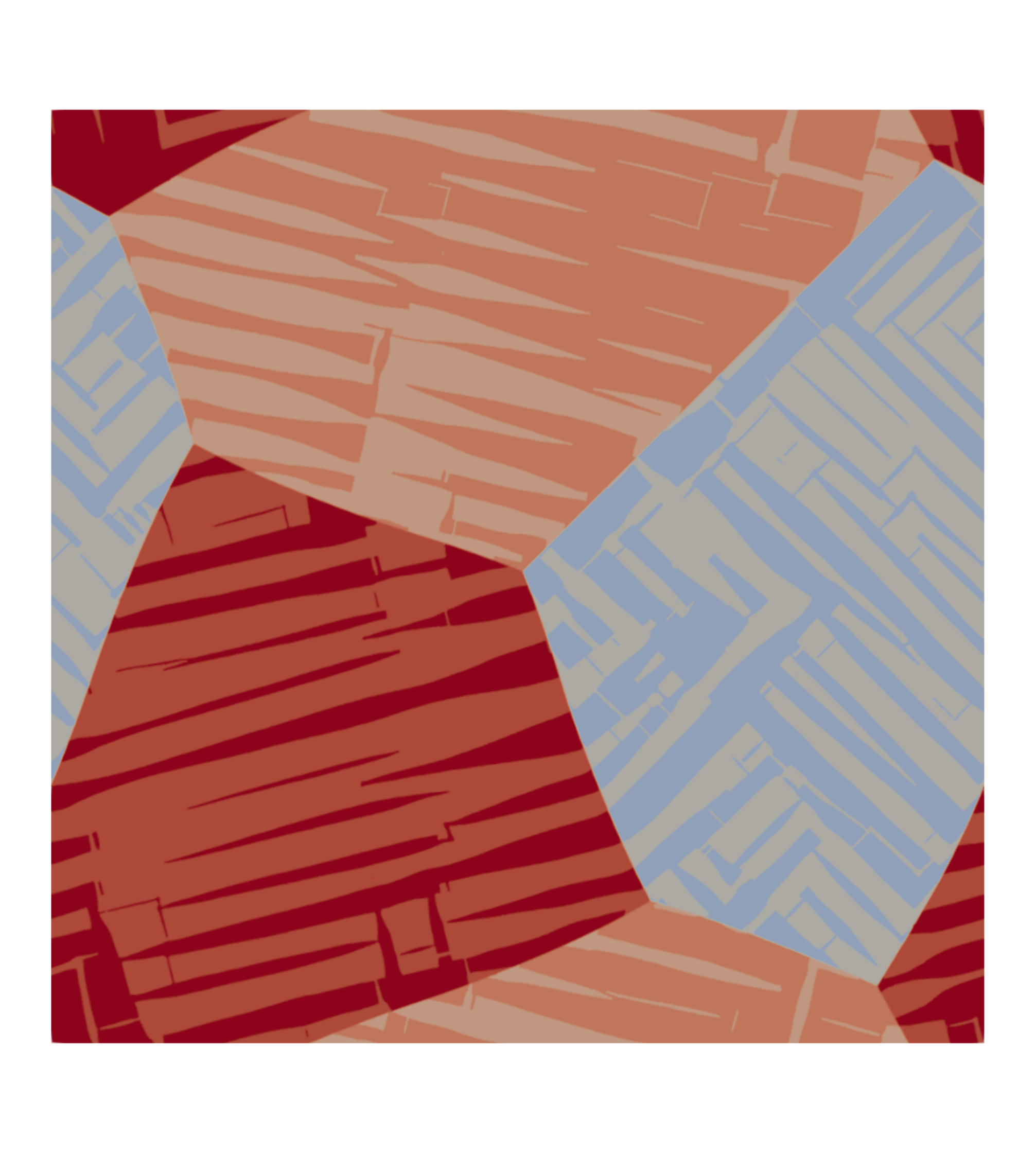}}
	\end{picture}} 
\subfigure[equal-strain, $t = 10000 \, \Delta t$]{\begin{picture}(50,55)
	\graphicspath{{figures/}}
	\put(0,0){\includegraphics[width=0.32\textwidth]{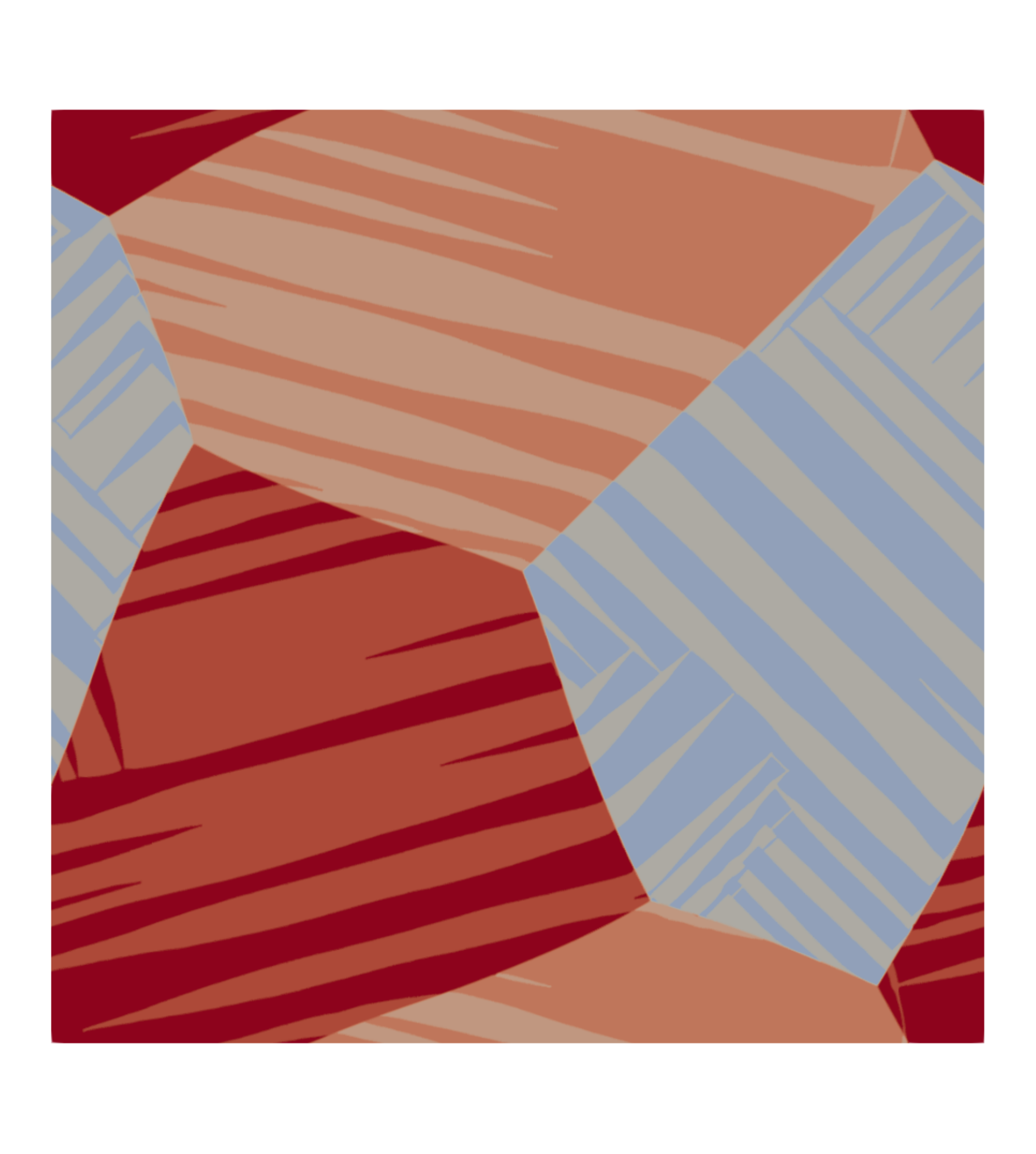}}
	\end{picture}} 
	\caption{Evolution of the microstructure for the two-dimensional Cubic to tetragonal martensitic transfomration.  }
	\label{fig:resultsfifththexample:1}
\end{figure}

The microstructure evolution using the presented models are shown in \Cref{fig:resultsfifththexample:1}. The final microstructures of all models are showing the same characteristic properties with martensitic  laminate with an interfaces perpendicular to the cubic lattice of the parent austenite.  The equal-strain approach delivers coarser laminates, since the interfacial regions are associated with higher elastic energy density. The equal stress approach ends with a closer microstructure to the one obtained by the MPFR1 model.

\begin{figure}[h!]
	\centering
	\unitlength=1mm
\subfigure[MPFR1, $t = 10000 \, \Delta t$]{\begin{picture}(50,70)
	\graphicspath{{figures/}}
	\put(0,0){\includegraphics[width=0.32\textwidth]{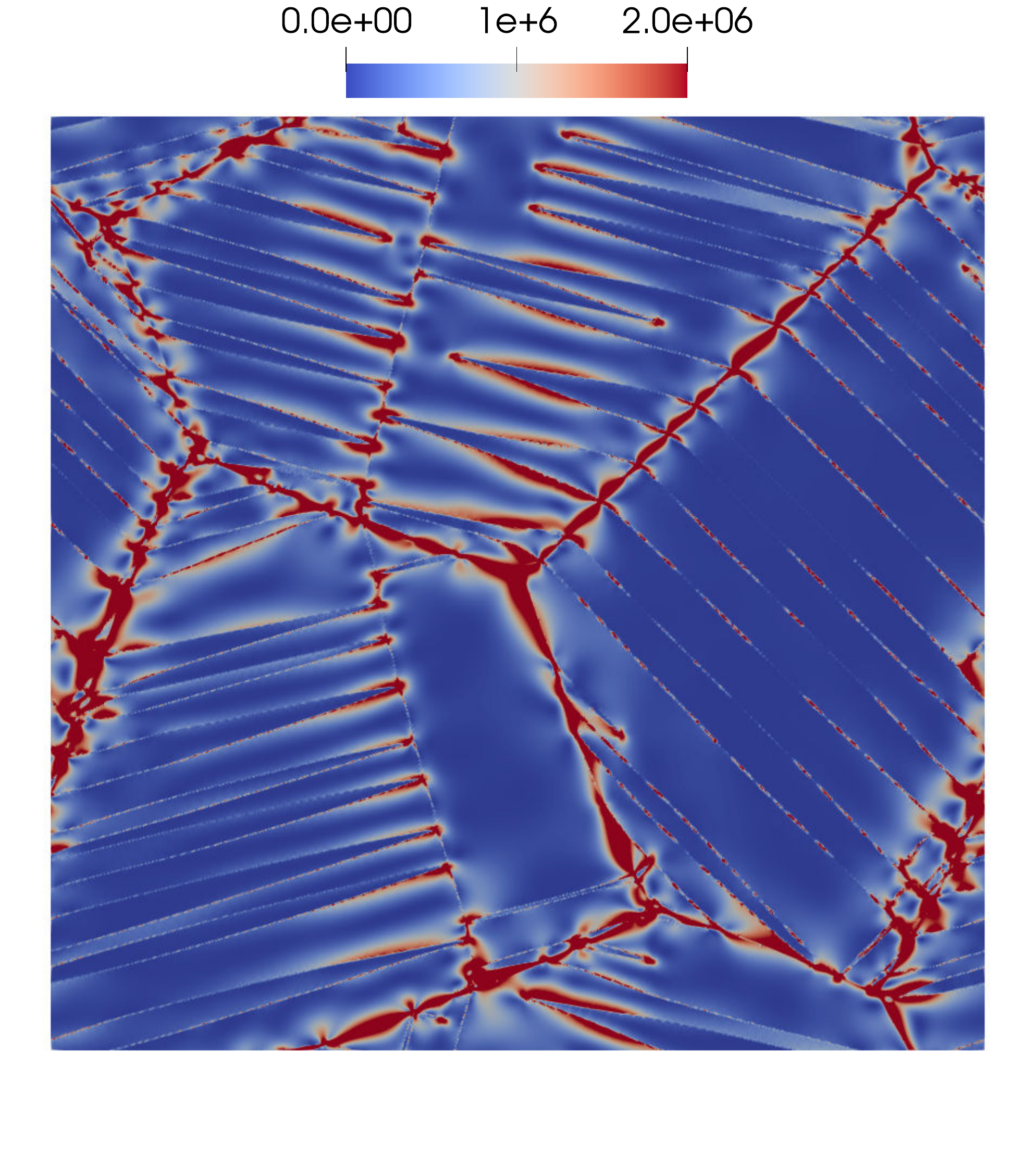}}
     \put(18.005323271,60.055607354){\color[rgb]{0,0,0}\rotatebox{0}{\makebox(0,0)[lb]{$\psi^\textrm{elas} \, (\textrm{J}/\textrm{m}^3)$}}} 
	\end{picture}} 
\subfigure[equal-stress, $t = 10000 \, \Delta t$]{\begin{picture}(50,50)
	\graphicspath{{figures/}}
	\put(0,0){\includegraphics[width=0.32\textwidth]{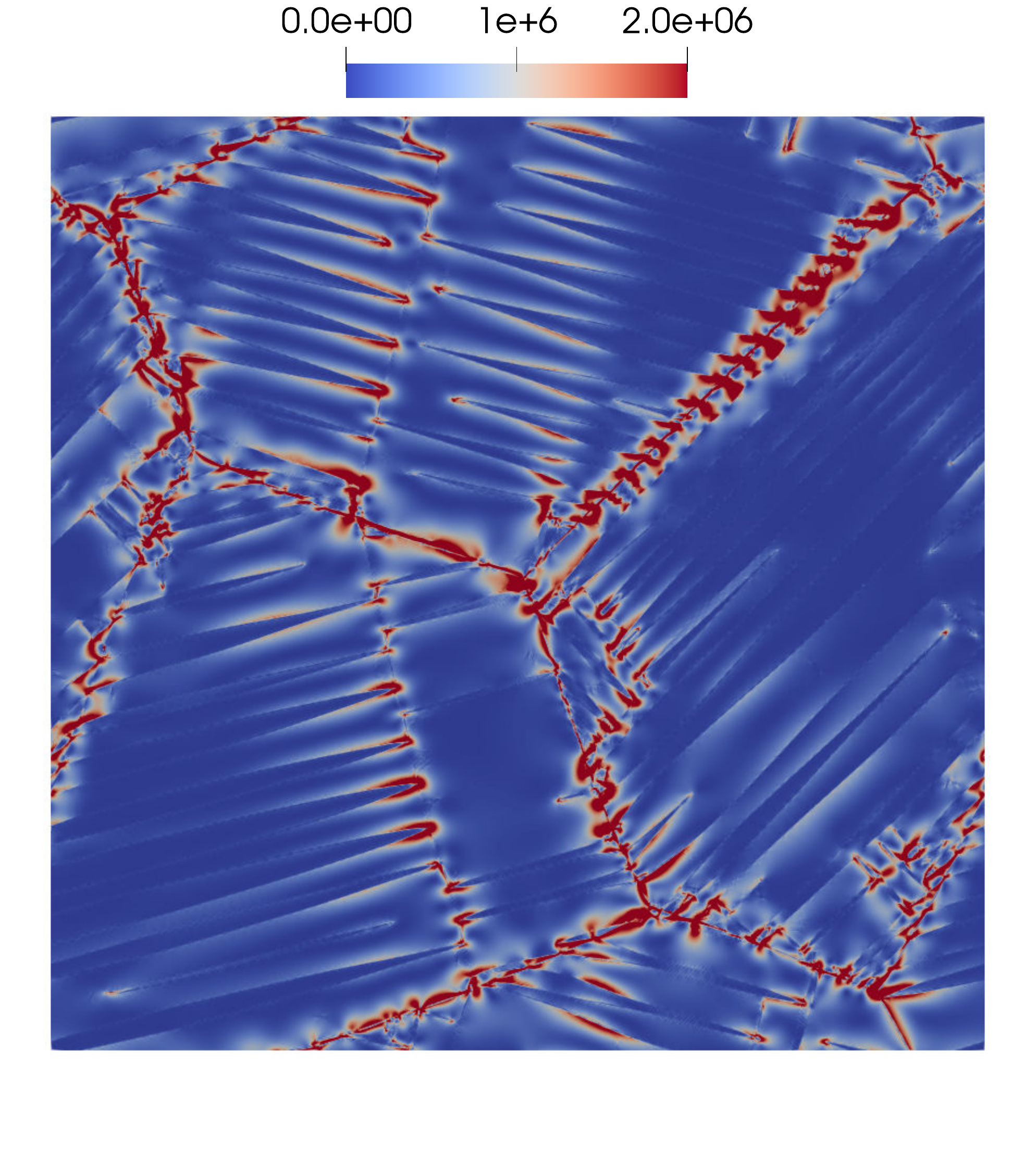}}
     \put(18.005323271,60.055607354){\color[rgb]{0,0,0}\rotatebox{0}{\makebox(0,0)[lb]{$\psi^\textrm{elas} \, (\textrm{J}/\textrm{m}^3)$}}} 
	\end{picture}} 
\subfigure[equal-strain, $t = 10000 \, \Delta t$]{\begin{picture}(50,50)
	\graphicspath{{figures/}}
	\put(0,0){\includegraphics[width=0.32\textwidth]{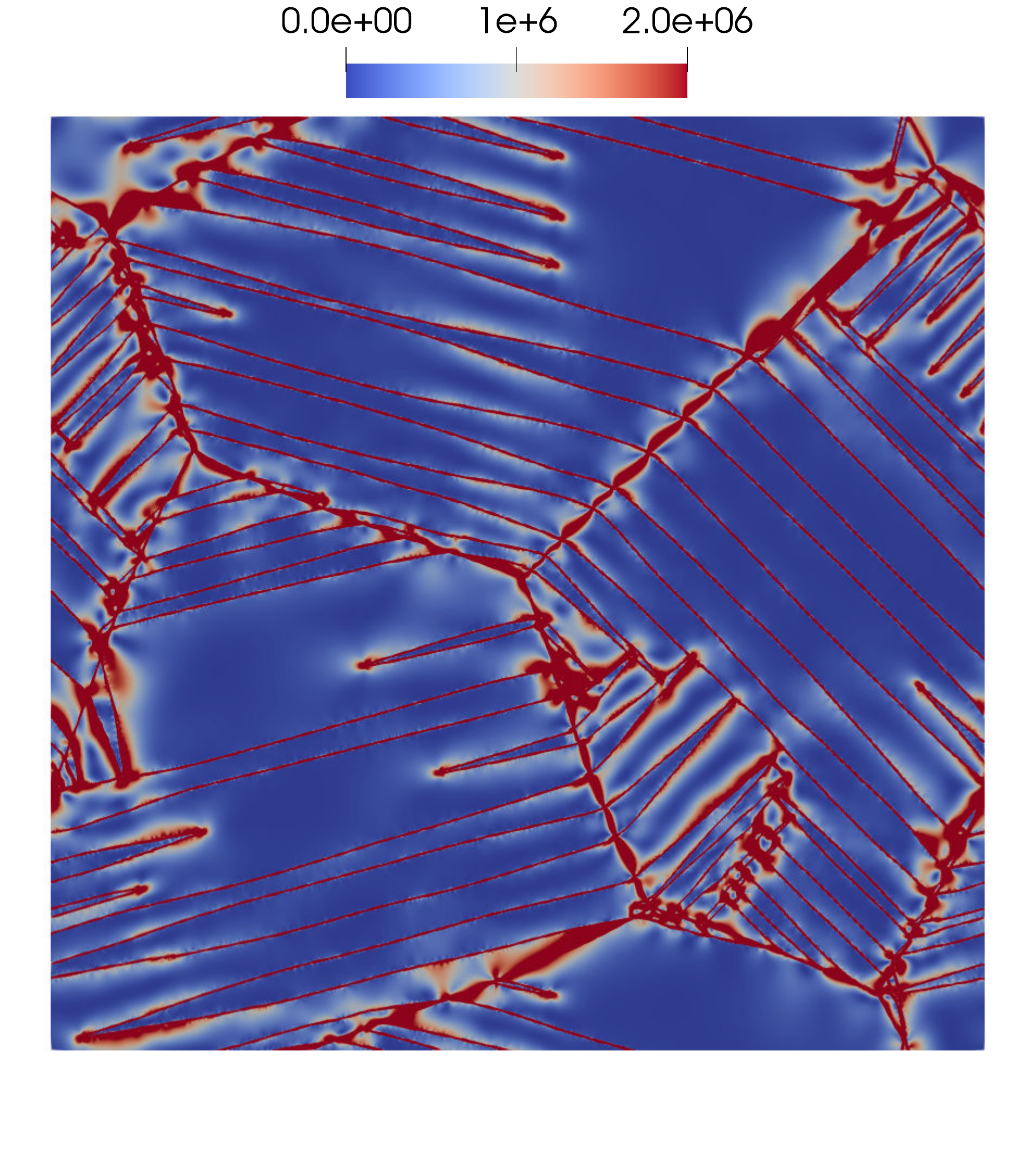}}
     \put(18.005323271,60.055607354){\color[rgb]{0,0,0}\rotatebox{0}{\makebox(0,0)[lb]{$\psi^\textrm{elas} \, (\textrm{J}/\textrm{m}^3)$}}} 
	\end{picture}} 
\subfigure[energy evolution]{\begin{picture}(150,75)
	\graphicspath{{figures/}}
	\put(0,5){\includegraphics[width=0.99\textwidth]{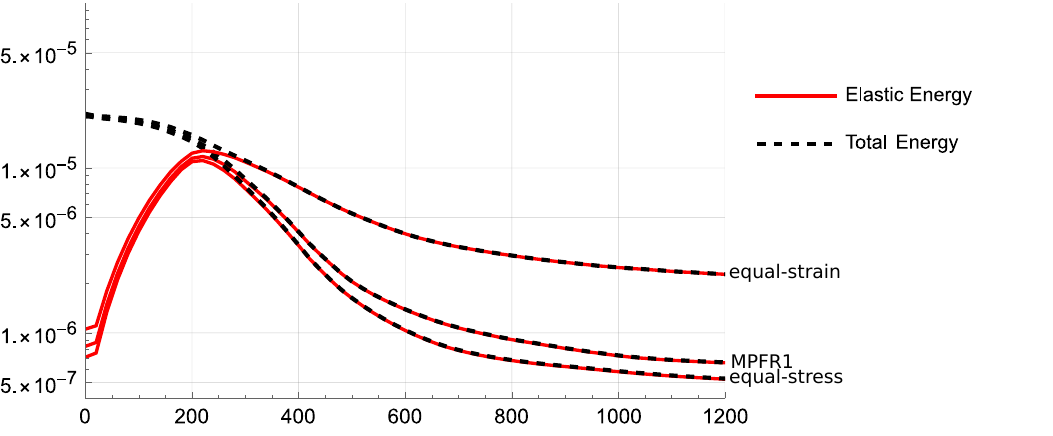}}
	\put(-3.005323271,30.055607354){\color[rgb]{0,0,0}\rotatebox{90}{\makebox(0,0)[lb]{ $\textrm{Joule}$}}}
     \put(60.005323271,0.055607354){\color[rgb]{0,0,0}\rotatebox{0}{\makebox(0,0)[lb]{ time step}}} 
	\end{picture}} 
	\caption{The energy evolution for the two-dimensional polycrystal example using the different models with the elastic energy density contours of the final microstructures (a-c).}
  \label{fig:resultsfifththexample:2}
\end{figure}

We demonstrate in  \Cref{fig:resultsfifththexample:2}(d) the evolution of the total and elastic energies for the first 1200 time steps.  Notice that we allow full energy relaxation to reach a final microstructure at  $t = 10000 \, \Delta t$.  The total energy is $1.99299 \cdot 10^{-05}$ Joule at the start $t=0$. When austenite vanishes from the simulation domain, the elastic energy become the dominating energy for all the implemented model.  The final total energy for equal-strain model reads $1.64382 \cdot 10^{-06}$ Joule which is four times higher than equal-stress assumption with  final total energy of $3.76084 \cdot 10^{-07}$ Joule.  The developed MPFR1 model delivers a total energy of $4.98237 \cdot 10^{-07}$ Joule which is $30 \%$ higher than equal-stress model.    We show the elastic energy density contours over the final microstructures in \Cref{fig:resultsfifththexample:2}(a-c). All the models exhibit elastic energy concentration on the grain boundaries. However, equal-strain model shows elastic energy concentration on all interfacial regions similar to single grain example. 

\FloatBarrier

%% file: figures/Example_1_geo.eps_tex
\begingroup%
  \makeatletter%
  \providecommand\color[2][]{%
    \errmessage{(Inkscape) Color is used for the text in Inkscape, but the package 'color.sty' is not loaded}%
    \renewcommand\color[2][]{}%
  }%
  \providecommand\transparent[1]{%
    \errmessage{(Inkscape) Transparency is used (non-zero) for the text in Inkscape, but the package 'transparent.sty' is not loaded}%
    \renewcommand\transparent[1]{}%
  }%
  \providecommand\rotatebox[2]{#2}%
  \ifx\svgwidth\undefined%
    \setlength{\unitlength}{721.61524751bp}%
    \ifx\svgscale\undefined%
      \relax%
    \else%
      \setlength{\unitlength}{\unitlength * \real{\svgscale}}%
    \fi%
  \else%
    \setlength{\unitlength}{\svgwidth}%
  \fi%
  \global\let\svgwidth\undefined%
  \global\let\svgscale\undefined%
  \makeatother%
  \begin{picture}(1,0.71182333)%
    \put(0,0){\includegraphics[width=\unitlength]{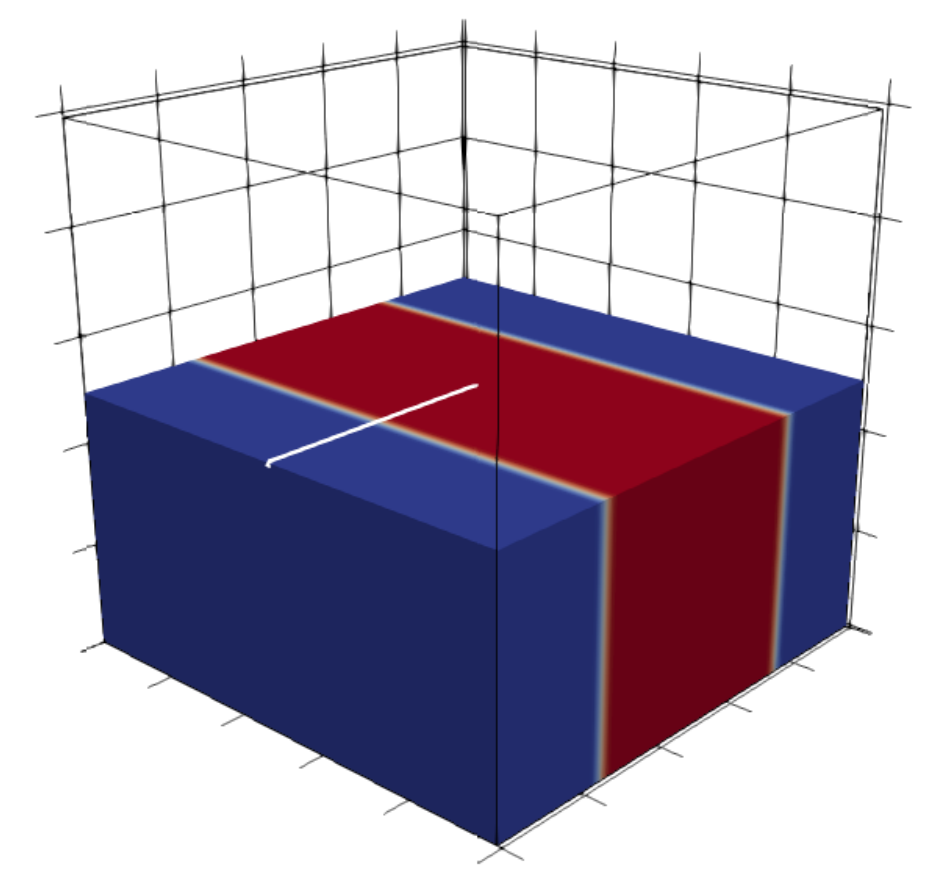}}
    \put(0.10323271,0.15607354){\color[rgb]{0,0,0}\rotatebox{-34}{\makebox(0,0)[lb]{\smash{$y$}}}}%
	\put(0.002456942,0.39876464){\color[rgb]{0,0,0}\rotatebox{90}{\makebox(0,0)[lb]{\smash{$x$}}}}%
    \put(0.693271,0.02607354){\color[rgb]{0,0,0}\rotatebox{30}{\makebox(0,0)[lb]{\smash{$z$}}}}%
    \put(0.04001797939,0.198073634){\color[rgb]{0,0,0}\makebox(0,0)[lb]{\smash{$0$}}}%
    \put(0.14917334,0.26086553){\color[rgb]{0,0,0}\makebox(0,0)[lb]{\smash{$0$}}}%
    \put(0.41285139,-0.0238698){\color[rgb]{0,0,0}\makebox(0,0)[lb]{\smash{$100$}}}%
    \put(-0.050557186,0.78193466){\color[rgb]{0,0,0}\makebox(0,0)[lb]{\smash{$100$}}}%
    \put(0.56285139,-0.0238698){\color[rgb]{0,0,0}\makebox(0,0)[lb]{\smash{$0$}}}%
    \put(0.934001797939,0.198073634){\color[rgb]{0,0,0}\makebox(0,0)[lb]{\smash{$100$}}}%
    \put(0.57001797939,0.21073634){\color[rgb]{1,1,1}\makebox(0,0)[lb]{\smash{$\alpha$}}}%
    \put(0.72001797939,0.28073634){\color[rgb]{1,1,1}\makebox(0,0)[lb]{\smash{$\beta$}}}%
    \put(0.86001797939,0.371073634){\color[rgb]{1,1,1}\makebox(0,0)[lb]{\smash{$\alpha$}}}%
  \end{picture}%
\endgroup%

%% file: figures/psi_df_1.eps_tex
\begingroup%
  \makeatletter%
  \providecommand\color[2][]{%
    \errmessage{(Inkscape) Color is used for the text in Inkscape, but the package 'color.sty' is not loaded}%
    \renewcommand\color[2][]{}%
  }%
  \providecommand\transparent[1]{%
    \errmessage{(Inkscape) Transparency is used (non-zero) for the text in Inkscape, but the package 'transparent.sty' is not loaded}%
    \renewcommand\transparent[1]{}%
  }%
  \providecommand\rotatebox[2]{#2}%
  \ifx\svgwidth\undefined%
    \setlength{\unitlength}{624.31111961bp}%
    \ifx\svgscale\undefined%
      \relax%
    \else%
      \setlength{\unitlength}{\unitlength * \real{\svgscale}}%
    \fi%
  \else%
    \setlength{\unitlength}{\svgwidth}%
  \fi%
  \global\let\svgwidth\undefined%
  \global\let\svgscale\undefined%
  \makeatother%
  \begin{picture}(1,0.52302736)%
    \put(0,0){\includegraphics[width=\unitlength]{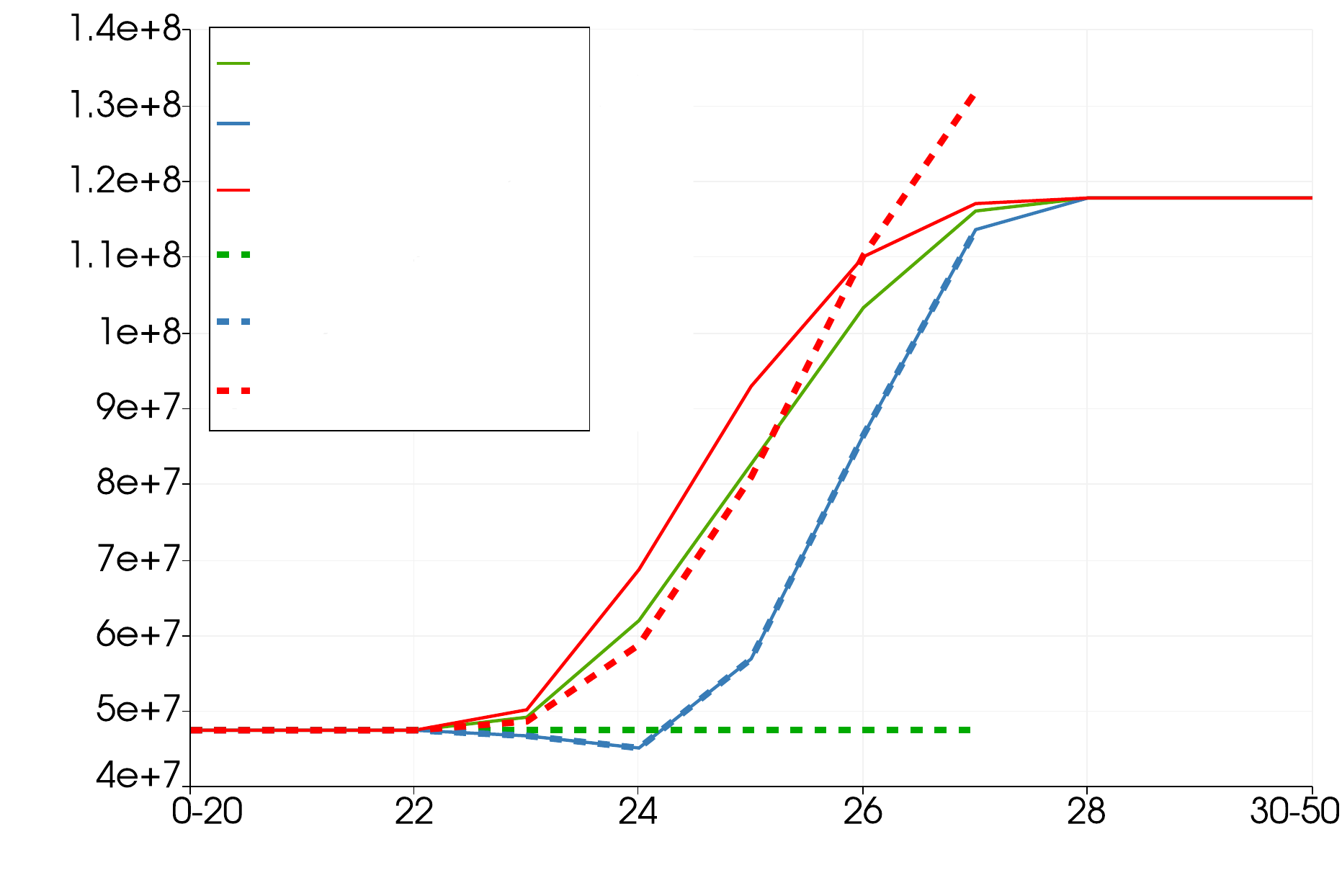}}%
    \put(0.200885861,0.520908704){\color[rgb]{0,0,0}\makebox(0,0)[lb]{\smash{\scriptsize $\psi^\textrm{elas,iso-strain}$}}}%
    \put(0.200944127,0.573168089){\color[rgb]{0,0,0}\makebox(0,0)[lb]{\smash{\scriptsize $\psi^\textrm{elas,iso-stress}$}}}%
    \put(0.20063731,0.61700973){\color[rgb]{0,0,0}\makebox(0,0)[lb]{\smash{\scriptsize $\psi^\textrm{elas}$}}}%
    \put(0.200885861,0.4820908704){\color[rgb]{0,0,0}\makebox(0,0)[lb]{\smash{\scriptsize $\psi^\textrm{ela}_\alpha$}}}%
    \put(0.200944127,0.423168089){\color[rgb]{0,0,0}\makebox(0,0)[lb]{\smash{\scriptsize $\psi^\textrm{elas,iso-stress}_\alpha$}}}%
    \put(0.20063731,0.36700973){\color[rgb]{0,0,0}\makebox(0,0)[lb]{\smash{\scriptsize $\psi^\textrm{elas,iso-strain}_\alpha$}}}%
    \put(0.47543573,0.00070059){\color[rgb]{0,0,0}\makebox(0,0)[lb]{\smash{$z \textrm{ in } 10^{-7} \textrm{ m}$}}}%
    \put(0.04,0.2723396){\color[rgb]{0,0,0}\rotatebox{90}{\makebox(0,0)[lb]{\smash{$\textrm{J}/\textrm{m}^3$}}}}%
    \put(0.15824291,0.21479315){\color[rgb]{0,0,0}\makebox(0,0)[lb]{\smash{phase-field $\alpha$}}}%
    \put(0.76855914,0.347873){\color[rgb]{0,0,0}\makebox(0,0)[lb]{\smash{phase-field $\beta$}}}%
  \end{picture}%
\endgroup%

%% file: figures/Example_1_DG_PR1.pdf_tex
\begingroup%
  \makeatletter%
  \providecommand\color[2][]{%
    \errmessage{(Inkscape) Color is used for the text in Inkscape, but the package 'color.sty' is not loaded}%
    \renewcommand\color[2][]{}%
  }%
  \providecommand\transparent[1]{%
    \errmessage{(Inkscape) Transparency is used (non-zero) for the text in Inkscape, but the package 'transparent.sty' is not loaded}%
    \renewcommand\transparent[1]{}%
  }%
  \providecommand\rotatebox[2]{#2}%
  \ifx\svgwidth\undefined%
    \setlength{\unitlength}{721.61524751bp}%
    \ifx\svgscale\undefined%
      \relax%
    \else%
      \setlength{\unitlength}{\unitlength * \real{\svgscale}}%
    \fi%
  \else%
    \setlength{\unitlength}{\svgwidth}%
  \fi%
  \global\let\svgwidth\undefined%
  \global\let\svgscale\undefined%
  \makeatother%
  \begin{picture}(1,0.71182333)%
    \put(0.0,0.01){\includegraphics[width=\unitlength]{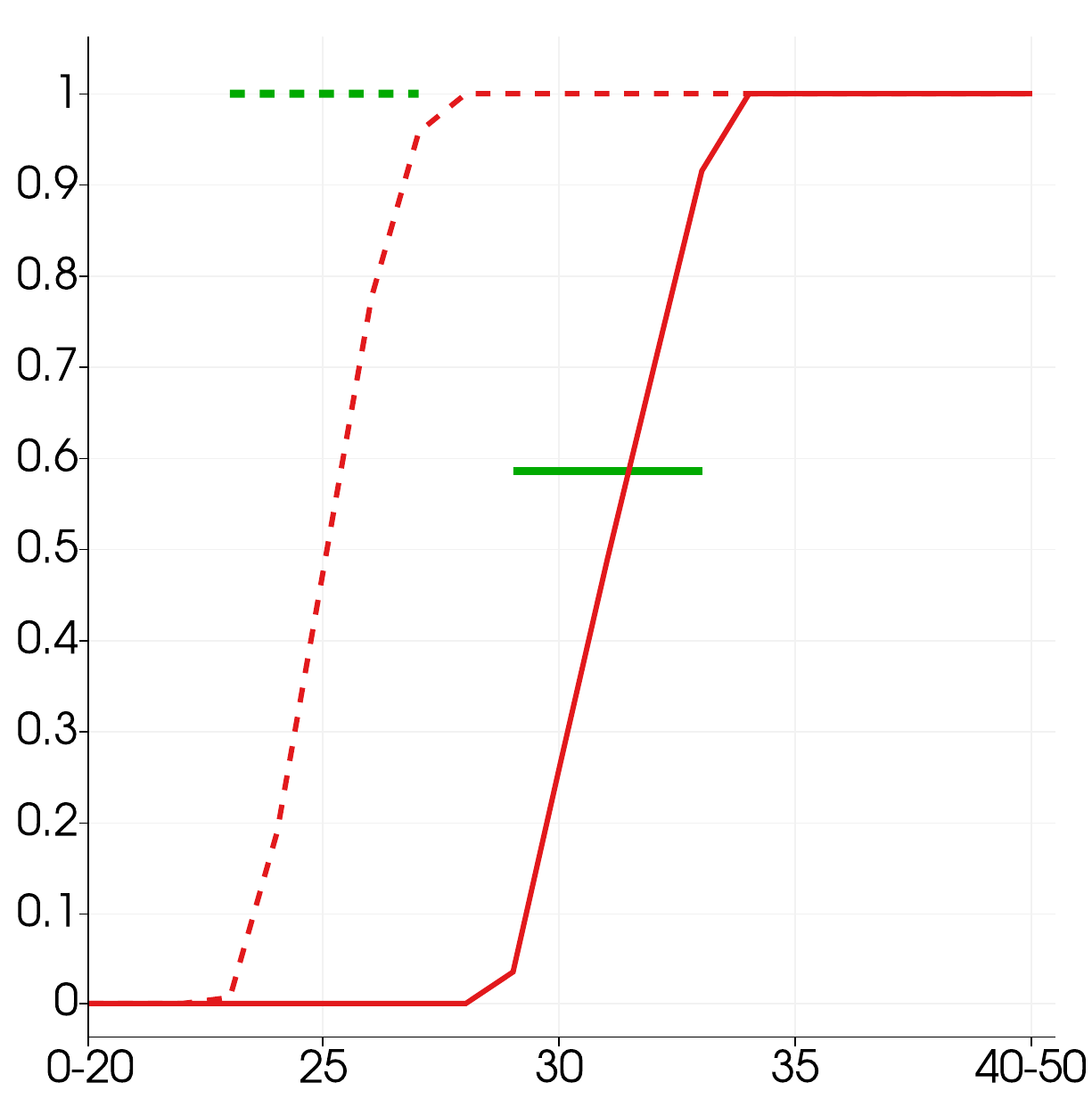}}
    \put(0.405323271,-0.035607354){\color[rgb]{0,0,0}\makebox(0,0)[lb]{\smash{$z \textrm{ in } 10^{-7} \textrm{m}$}}}%
	\put(-0.02456942,0.14876464){\color[rgb]{0,0.5,0}\rotatebox{90}{\makebox(0,0)[lb]{\smash{${\color{red} \phi_\beta} \,, {\Delta G^\textrm{elas}_{\alpha\beta} / (2.4 \, 10^{7}) } \, \textrm{in} \,  \textrm{J}/\textrm{m}^3$}}}}%
    \put(0.245323271,0.635607354){\color[rgb]{0,0,0}\makebox(0,0)[lb]{\smash{$t_0$}}}%
    \put(0.605323271,0.635607354){\color[rgb]{0,0,0}\makebox(0,0)[lb]{\smash{$t>t_0$}}}%
  \end{picture}%
\endgroup%

%% file: figures/Example_1_DG_Reuss.pdf_tex
\begingroup%
  \makeatletter%
  \providecommand\color[2][]{%
    \errmessage{(Inkscape) Color is used for the text in Inkscape, but the package 'color.sty' is not loaded}%
    \renewcommand\color[2][]{}%
  }%
  \providecommand\transparent[1]{%
    \errmessage{(Inkscape) Transparency is used (non-zero) for the text in Inkscape, but the package 'transparent.sty' is not loaded}%
    \renewcommand\transparent[1]{}%
  }%
  \providecommand\rotatebox[2]{#2}%
  \ifx\svgwidth\undefined%
    \setlength{\unitlength}{721.61524751bp}%
    \ifx\svgscale\undefined%
      \relax%
    \else%
      \setlength{\unitlength}{\unitlength * \real{\svgscale}}%
    \fi%
  \else%
    \setlength{\unitlength}{\svgwidth}%
  \fi%
  \global\let\svgwidth\undefined%
  \global\let\svgscale\undefined%
  \makeatother%
  \begin{picture}(1,0.71182333)%
    \put(0.0,0.01){\includegraphics[width=\unitlength]{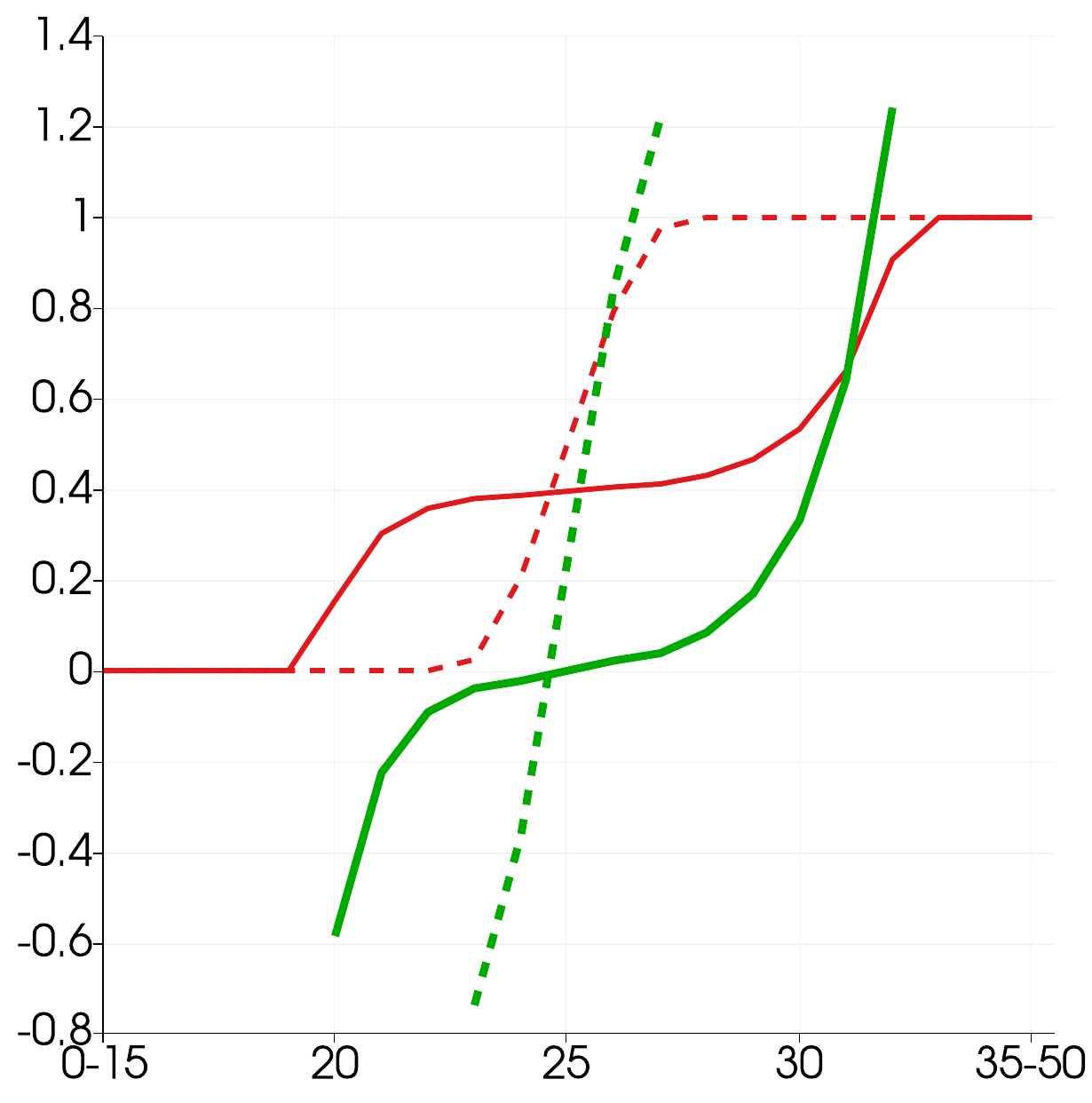}}
    \put(0.405323271,-0.035607354){\color[rgb]{0,0,0}\makebox(0,0)[lb]{\smash{$z \textrm{ in } 10^{-7} \textrm{m}$}}}%
	\put(-0.02456942,0.14876464){\color[rgb]{0,0.5,0}\rotatebox{90}{\makebox(0,0)[lb]{\smash{${\color{red} \phi_\beta} \,, {\Delta G^\textrm{elas}_{\alpha\beta} /  10^{8} } \, \textrm{in} \,  \textrm{J}/\textrm{m}^3$}}}}%
    \put(0.4705323271,0.635607354){\color[rgb]{0,0,0}\makebox(0,0)[lb]{\smash{$t_0$}}}%
    \put(0.7805323271,0.635607354){\color[rgb]{0,0,0}\makebox(0,0)[lb]{\smash{$t>t_0$}}}%
  \end{picture}%
\endgroup%

%% file: figures/Example_1_DG_Voigt.pdf_tex
\begingroup%
  \makeatletter%
  \providecommand\color[2][]{%
    \errmessage{(Inkscape) Color is used for the text in Inkscape, but the package 'color.sty' is not loaded}%
    \renewcommand\color[2][]{}%
  }%
  \providecommand\transparent[1]{%
    \errmessage{(Inkscape) Transparency is used (non-zero) for the text in Inkscape, but the package 'transparent.sty' is not loaded}%
    \renewcommand\transparent[1]{}%
  }%
  \providecommand\rotatebox[2]{#2}%
  \ifx\svgwidth\undefined%
    \setlength{\unitlength}{721.61524751bp}%
    \ifx\svgscale\undefined%
      \relax%
    \else%
      \setlength{\unitlength}{\unitlength * \real{\svgscale}}%
    \fi%
  \else%
    \setlength{\unitlength}{\svgwidth}%
  \fi%
  \global\let\svgwidth\undefined%
  \global\let\svgscale\undefined%
  \makeatother%
  \begin{picture}(1,0.71182333)%
    \put(0.0,0.01){\includegraphics[width=\unitlength]{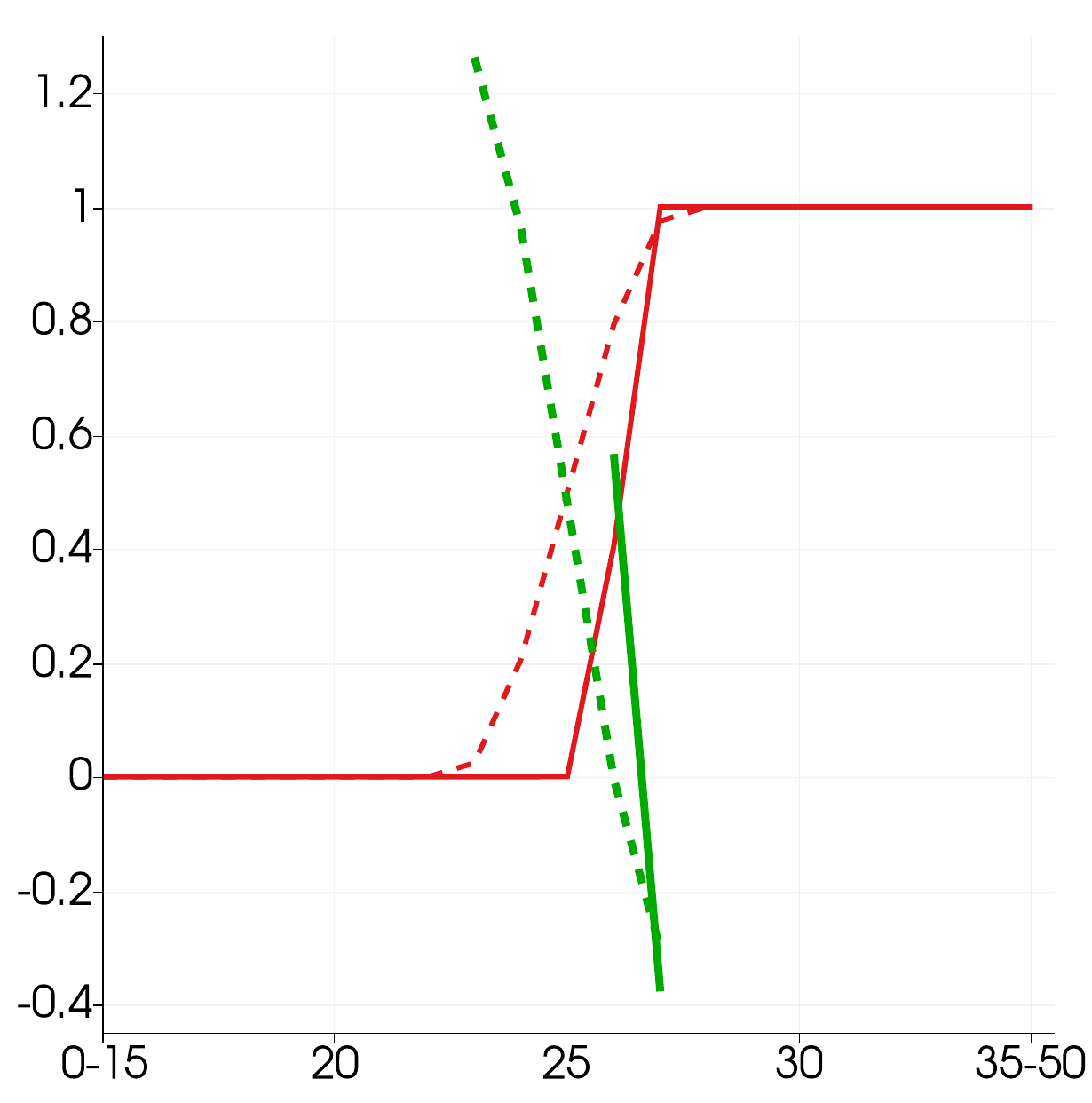}}
    \put(0.405323271,-0.035607354){\color[rgb]{0,0,0}\makebox(0,0)[lb]{\smash{$z \textrm{ in } 10^{-7} \textrm{m}$}}}%
	\put(-0.02456942,0.14876464){\color[rgb]{0,0.5,0}\rotatebox{90}{\makebox(0,0)[lb]{\smash{${\color{red} \phi_\beta} \,, {\Delta G^\textrm{elas}_{\alpha\beta} / (5 \, 10^{7})} \, \textrm{in} \,  \textrm{J}/\textrm{m}^3$}}}}%
    \put(0.425323271,0.435607354){\color[rgb]{0,0,0}\makebox(0,0)[lb]{\smash{$t_0$}}}%
    \put(0.605323271,0.635607354){\color[rgb]{0,0,0}\makebox(0,0)[lb]{\smash{$t>t_0$}}}%
  \end{picture}%
\endgroup%

%% file: figures/Example_1_DG_PR1_enhanced.pdf_tex
\begingroup%
  \makeatletter%
  \providecommand\color[2][]{%
    \errmessage{(Inkscape) Color is used for the text in Inkscape, but the package 'color.sty' is not loaded}%
    \renewcommand\color[2][]{}%
  }%
  \providecommand\transparent[1]{%
    \errmessage{(Inkscape) Transparency is used (non-zero) for the text in Inkscape, but the package 'transparent.sty' is not loaded}%
    \renewcommand\transparent[1]{}%
  }%
  \providecommand\rotatebox[2]{#2}%
  \ifx\svgwidth\undefined%
    \setlength{\unitlength}{721.61524751bp}%
    \ifx\svgscale\undefined%
      \relax%
    \else%
      \setlength{\unitlength}{\unitlength * \real{\svgscale}}%
    \fi%
  \else%
    \setlength{\unitlength}{\svgwidth}%
  \fi%
  \global\let\svgwidth\undefined%
  \global\let\svgscale\undefined%
  \makeatother%
  \begin{picture}(1,0.71182333)%
    \put(0.0,0.01){\includegraphics[width=\unitlength]{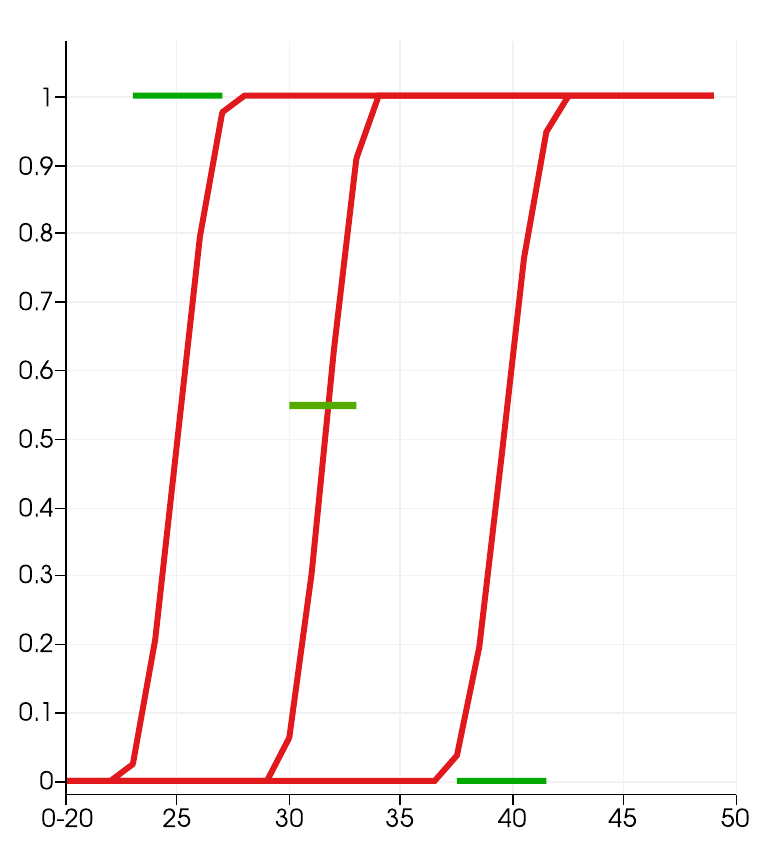}}
    \put(0.405323271,-0.035607354){\color[rgb]{0,0,0}\makebox(0,0)[lb]{\smash{$z \textrm{ in } 10^{-7} \textrm{m}$}}}%
	\put(-0.02456942,0.14876464){\color[rgb]{0,0.5,0}\rotatebox{90}{\makebox(0,0)[lb]{\smash{${\color{red} \phi_\beta} \,, {\Delta G^\textrm{elas}_{\alpha\beta} / (2.4 \, 10^{7}) } \, \textrm{in} \,  \textrm{J}/\textrm{m}^3$}}}}%
    \put(0.245323271,0.635607354){\color[rgb]{0,0,0}\makebox(0,0)[lb]{\smash{$t_0$}}}%
    \put(0.4405323271,0.635607354){\color[rgb]{0,0,0}\makebox(0,0)[lb]{\smash{$t_1$}}}%
    \put(0.680405323271,0.635607354){\color[rgb]{0,0,0}\makebox(0,0)[lb]{\smash{$t_2$}}}%
  \end{picture}%
\endgroup%

%% file: figures/Example_1_DG_Reuss_enhanced.pdf_tex
\begingroup%
  \makeatletter%
  \providecommand\color[2][]{%
    \errmessage{(Inkscape) Color is used for the text in Inkscape, but the package 'color.sty' is not loaded}%
    \renewcommand\color[2][]{}%
  }%
  \providecommand\transparent[1]{%
    \errmessage{(Inkscape) Transparency is used (non-zero) for the text in Inkscape, but the package 'transparent.sty' is not loaded}%
    \renewcommand\transparent[1]{}%
  }%
  \providecommand\rotatebox[2]{#2}%
  \ifx\svgwidth\undefined%
    \setlength{\unitlength}{721.61524751bp}%
    \ifx\svgscale\undefined%
      \relax%
    \else%
      \setlength{\unitlength}{\unitlength * \real{\svgscale}}%
    \fi%
  \else%
    \setlength{\unitlength}{\svgwidth}%
  \fi%
  \global\let\svgwidth\undefined%
  \global\let\svgscale\undefined%
  \makeatother%
  \begin{picture}(1,0.71182333)%
    \put(0.0,0.01){\includegraphics[width=\unitlength]{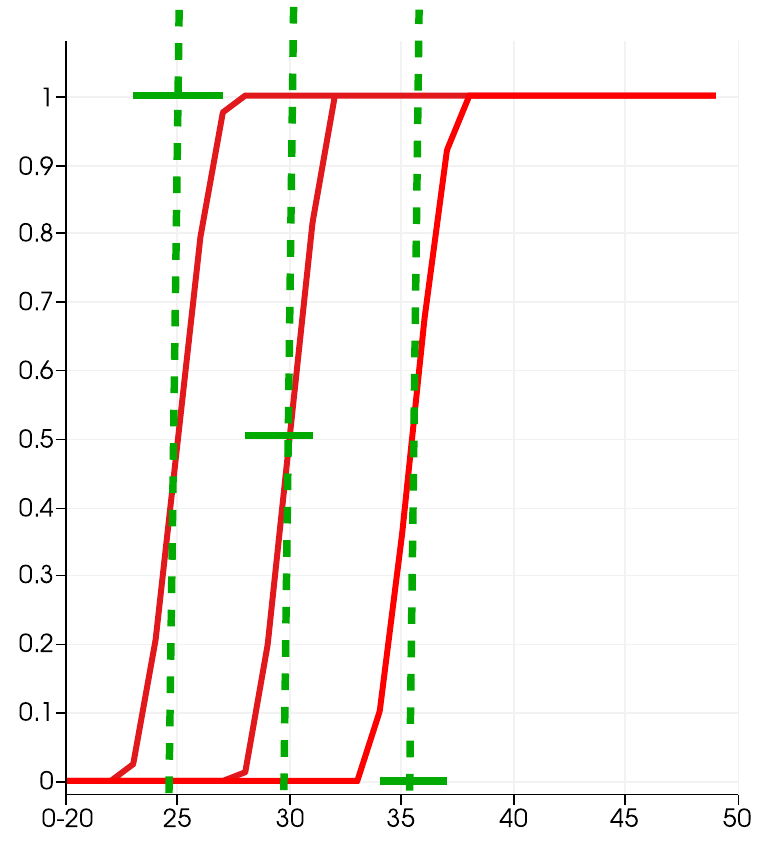}}
    \put(0.405323271,-0.035607354){\color[rgb]{0,0,0}\makebox(0,0)[lb]{\smash{$z \textrm{ in } 10^{-7} \textrm{m}$}}}%
	\put(-0.02456942,0.14876464){\color[rgb]{0,0.5,0}\rotatebox{90}{\makebox(0,0)[lb]{\smash{${\color{red} \phi_\beta} \,, {\Delta G^\textrm{elas}_{\alpha\beta} / (2.4 \, 10^{7}) } \, \textrm{in} \,  \textrm{J}/\textrm{m}^3$}}}}%
    \put(0.2505323271,0.635607354){\color[rgb]{0,0,0}\makebox(0,0)[lb]{\smash{$t_0$}}}%
    \put(0.405323271,0.635607354){\color[rgb]{0,0,0}\makebox(0,0)[lb]{\smash{$t_1$}}}%
    \put(0.565323271,0.635607354){\color[rgb]{0,0,0}\makebox(0,0)[lb]{\smash{$t_2$}}}%
  \end{picture}%
\endgroup%

%% file: figures/Example_1_DG_Voigt_enhanced.pdf_tex
\begingroup%
  \makeatletter%
  \providecommand\color[2][]{%
    \errmessage{(Inkscape) Color is used for the text in Inkscape, but the package 'color.sty' is not loaded}%
    \renewcommand\color[2][]{}%
  }%
  \providecommand\transparent[1]{%
    \errmessage{(Inkscape) Transparency is used (non-zero) for the text in Inkscape, but the package 'transparent.sty' is not loaded}%
    \renewcommand\transparent[1]{}%
  }%
  \providecommand\rotatebox[2]{#2}%
  \ifx\svgwidth\undefined%
    \setlength{\unitlength}{721.61524751bp}%
    \ifx\svgscale\undefined%
      \relax%
    \else%
      \setlength{\unitlength}{\unitlength * \real{\svgscale}}%
    \fi%
  \else%
    \setlength{\unitlength}{\svgwidth}%
  \fi%
  \global\let\svgwidth\undefined%
  \global\let\svgscale\undefined%
  \makeatother%
  \begin{picture}(1,0.71182333)%
    \put(0.0,0.01){\includegraphics[width=\unitlength]{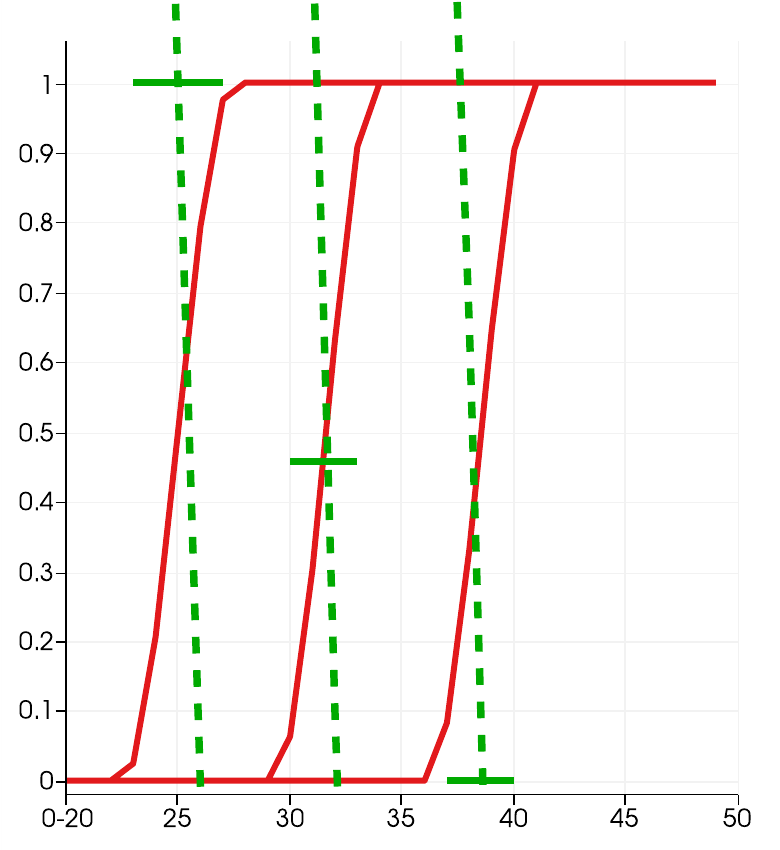}}
    \put(0.405323271,-0.035607354){\color[rgb]{0,0,0}\makebox(0,0)[lb]{\smash{$z \textrm{ in } 10^{-7} \textrm{m}$}}}%
	\put(-0.02456942,0.14876464){\color[rgb]{0,0.5,0}\rotatebox{90}{\makebox(0,0)[lb]{\smash{${\color{red} \phi_\beta} \,, {\Delta G^\textrm{elas}_{\alpha\beta} / (2.4 \, 10^{7})} \, \textrm{in} \,  \textrm{J}/\textrm{m}^3$}}}}%
    \put(0.2505323271,0.635607354){\color[rgb]{0,0,0}\makebox(0,0)[lb]{\smash{$t_0$}}}%
    \put(0.445323271,0.635607354){\color[rgb]{0,0,0}\makebox(0,0)[lb]{\smash{$t_1$}}}%
    \put(0.635323271,0.635607354){\color[rgb]{0,0,0}\makebox(0,0)[lb]{\smash{$t_2$}}}%
  \end{picture}%
\endgroup%

%% file: Conclusion.tex
\section{Final conclusions} 
\label{sect:fin}
In this work, a novel multi-phase-field elasticity model is developed presenting a new definition for the phase-field elastic energy which is linearly interpolated from pairwise phase-field energies. The pairwise phase-field energy is assumed as a function of a pairwise strain after enforcing Hadamard kinematic compatibility  between the two relevant phase-fields. The pairwise static equilibrium is satisfied then by calculating the pairwise strain jump vector through the relaxation of the pairwise elastic energy. Thus, the developed model satisfies the pairwise jump conditions for all active phase-field pairs on their pairwise normals. Different numerical examples were designed to test the developed MPFR1 model against the upper and lower bounds represented by equal-strain and equal-stress assumptions. 

In the first numerical example, the MPFR1 model quantitatively meets the sharp dynamics for a planar interface exhibiting a superior constant driving force compared to the classical equal-strain and equal-stress approaches. We examined in the second numerical example the shape of the growth of a martensitic nucleus in an austenitic matrix where we showed the surface tension effect, which should be added in future work to the MPFR1 model.  We have studied the energy and driving force profile at the triple junction in third example showing a semi-constant energy profile for the MPFR1 model. Finally, we model the cubic-to-tetragonal transformation in a 3D single grain in the fourth numerical example and in a 2D polycrystal in the fifth example. The MPFR1 behaves between the equal-strain and the equal-stress models, but it is able to predict a stress-free laminate in the 3D single grain example as  the equal-stress model does. For the 2D polycrystal example, the final energy of the microstructure of MPFR1 model is 30$\%$ more than the equal-stress assumption but 70$\%$ less than the equal-strain assumption.

Satisfying the jump conditions at the interface resulted in quantitative results with transformation rates and stress level between the nonphysical upper and lower bounds of elasticity presented by equal-strain and equal-stress assumption. Most of the mechanical models used in the phase-field community are interpolation models, which align with the equal-stress assumption and therefore the mechanical energy is underestimated. For this paper, we considered only elasticity at small strains, however, the role of the mechanical models will be more pronounced when plasticity and viscoplasticity are considered which will be tested in future works.